\RequirePackage{fix-cm} 
\documentclass[twoside,reqno,11pt]{amsart}

\usepackage{a4wide,verbatim}
\usepackage[english]{babel}
\usepackage[latin1]{inputenc}
\usepackage{amsmath,amsfonts,amssymb,amsgen,amsbsy}
\usepackage{eucal,esint,dsfont,mathrsfs}
\usepackage{graphicx}
\usepackage{subcaption}
\usepackage[square,numbers]{natbib}
\usepackage{mathtools}
\usepackage{hyperref}
\usepackage{float}
\usepackage{slashbox}
\usepackage{nicematrix}

\newcommand{\eqdef}{\stackrel{\text{\tiny{def}}}{=}} 
\newcommand{\mathsc}[1]{{\normalfont{\textsc{#1}}}}

\newcommand{\ud}{\mathrm{d}}

\newcommand{\ue}{\mathrm{e}}
\newcommand{\ui}{\ii}

\renewcommand{\Re}{\operatorname{Re}}

\newlength{\intwidth}

\renewcommand{\leq}{\leqslant}
\renewcommand{\geq}{\geqslant}

\title[Zigzag finite difference schemes]{\bf A new class of finite difference 
methods:\\ The zigzag schemes}
\author[L.  Poggioni, D. Clamond \& Y. D'Angelo]{Lorenzo Poggioni, Didier Clamond \& Yves D'Angelo}
\address{Universit\'e C\^ote d'Azur, CNRS, LJAD, France.} 
\email{Lorenzo.Poggioni@univ-cotedazur}
\email{Didier.Clamond@univ-cotedazur.fr}
\email{Yves.Dangelo@univ-cotedazur.fr}

\DeclarePairedDelimiter\ceil{\lceil}{\rceil}
\DeclarePairedDelimiter\floor{\lfloor}{\rfloor}
\newcommand{\e}{\text{e}}
\newcommand{\ii}{\text{i}}
\newcommand{\sinc}{\text{sinc}}
\newcommand{\kmax}{k_\text{max}}

\date{\today}

\begin{document}
\maketitle

\begin{abstract}
We introduce a novel class of finite difference approximations, termed \textit{zigzag schemes}, that
employ a hybrid stencil that is neither symmetrical, nor fully one-sided. These zigzag schemes often
enjoy more permissive stability constraints and see their coefficients vanish as the order tends to 
infinity. This property permits the formulation of higher order schemes. An explicit formula is given 
for both collocated and staggered grids for an arbitrary order and a closed-form expression for the 
infinite-order scheme is also provided. A linear stability analysis indicates that the zigzag scheme
offer a broader range of conditional stability compared to the centred and upwind schemes, sometimes
being the only stable scheme.
Additionally, the asymmetrical structure of the stencil of zigzag schemes prevents some issues 
such as the formation of ``ghost solutions''. Moreover, implementing zigzag schemes is relatively 
easy when a code using classical finite differences is available, that is an important feature  
for well-tested legacy codes. 
Overall, zigzag schemes provide a compelling alternative for finite differences methods by 
enabling faster and more stable numerical simulations without sacrificing accuracy or ease of use.

\end{abstract}




\section{Introduction}\label{Introduction} %

Introduced by Taylor \citep{taylor1715} and Euler \citep{euler1768} in the 
eighteenth century, finite difference methods remain a cornerstone in the 
numerical approximation of partial differences equations. These methods are 
used in a wide array of domains, including the study of immersed boundary 
problems \citep{Fadlun_etal2000,Lai_etal2000,LeVeque_etal1994,PeskinCharles2002} 
and the study of hyperbolic conservation laws \citep{Godlewski_etal1996,
HartenAmi1983, SodGary1978}. Additionally, the local nature of finite difference
schemes make them an adapted tool regarding high performance computing, for 
both CPU \citep{Abide_etal2017, Szustak_etal2016} and GPU \citep{Giles_etal2014,
Hamilton_etal2013, Ye2022}. Traditional schemes, most notably the centred 
and upwind schemes, are well known for their effectiveness and their simplicity.
However, these schemes are often characterised by tight restrictions regarding 
their stability. These restrictions also tend to worsen as the order of the 
schemes increases. Thus, hybrid difference schemes have been proposed for 
convection-diffusion problems \citep{Spalding1972}. The idea is to use 
locally a centred or a upwind/downwind schemes according to the characteristics 
of the problem, and thus to take advantage of the good properties of both schemes. 
An alternative approach is described in the present paper. 

In this work, we introduce a novel finite difference method, termed 
\textit{zigzag schemes}, designed to relax the aforementioned restrictions 
while retaining the simplicity of implementation enjoyed by the centred and 
upwind schemes. This allows the zigzag schemes to be implementable into already 
well established legacy codes based on finite difference methods such as the 
direct numerical simulation (DNS) solver H-Allegro \citep{HAllegro}.
Furthermore, the performance of this new general method does not rely on
properties inherent to a specific class of equation. Contrary to the upwind 
schemes, the coefficients of the zigzag schemes vanish as the order tends to 
infinity, which allows the use of higher-order schemes.
These zigzag schemes are characterised by stencils that are neither fully 
centred nor fully one-sided, but instead represent a hybrid of both. We 
present an explicit formula for the coefficients of these schemes for 
arbitrary orders, as well as a closed-form expression for the infinite-order 
case.

The paper is organised as follows. Sections \ref{S1} and \ref{S2} present, 
respectively, the classical centred and onde-sided (upwind) schemes. In section 
\ref{S3}, we introduce the new collocated and staggered zigzag schemes. Section 
\ref{S4} is devoted to a linear stability analysis of the zigzag schemes, 
as well as a comparison with the centred and upwind schemes 
showcasing situations where the zigzag schemes are the only stable schemes.
In the section \ref{SectionGhostSol}, we highlight a situation where zigzag 
schemes are able to avoid problems inherent to symmetrical schemes (such as 
the formation of ``ghost solutions'') without being limited to lower
orders.

\section{Classical centred finite differences}\label{S1}

\subsection{Schemes definition} 

Consider a function $f$ of one real variable $x$, discretised with a constant step size $\Delta x$. 
The discrete abscissa are noted as $x_i \eqdef i \Delta x$ and the discrete value of $f$ at $x_i$ is 
noted $f_i = f(x_i)$ for $i \in \lbrace 1,2,...,N \rbrace$.
A $2N$-th order centred finite approximation of $f'$ and $f''$ can be written \citep{Fornberg1990}
\begin{align}
f_i' &= \sum_{j=1}^N C_N^j \frac{f_{i+j} - f_{i-j}}{2j\Delta x}, 
 \label{Centred1} \\
f_i'' &= \sum_{j=1}^N \, C_N^j \, \frac{f_{i+j} - 2 f_i + 
f_{i-j}}{(j \, \Delta x)^2},  \label{Centred2}
\end{align}
where
\begin{equation} \label{EQ_centred_coeffs}
C_N^j \eqdef \frac{2 (-1)^{j+1} (N!)^2}{(N+j)!(N-j)!}.
\end{equation}
Note that as $N \rightarrow\infty$, the centred coefficients become  
\citep{fornberg1998}
\begin{equation}\label{Cjinf}
C_{\infty}^j = 2 (-1)^{j+1},
\end{equation}
an infinite-order scheme can then be considered. If the interval of definition 
of $f$ is discretised using $P\geqslant 2$ equally spaced points, 
then necessarily $1\leqslant N\leqslant\floor{(P-1)/2}$ ($\floor{\cdot}$ the 
rounding towards $-\infty$), and appropriate boundary conditions must be taken 
into account to deduce values at virtual nodes outside the definition domain.
Note that the information contained in (\ref{Centred1}) and (\ref{Centred2}) is sufficient to approximate derivatives of any order, since higher--order derivatives can be expressed as combinations of first- and second-order derivatives. 

\subsection{Staggered variant}

Instead of using nodes at integer positions $x=j\/\Delta x$, one may advantageously 
consider nodes at intermediate positions.
The generic staggered centred scheme of order $2N$ then is \citep{Fornberg1990} :
\begin{equation}
f_i' = \sum_{j=1}^N \, \breve{C}_N^j \, \frac{f_{i+j-1/2} - f_{i-j+1/2}}{(2j - 1) \, \Delta x},
\end{equation}
where the coefficients $\breve{C}_N^j$ are defined by
\begin{equation}
\breve{C}_N^j \eqdef \frac{(-1)^{j+1} \, 4^{1 - 2N} \, (2
N)!^2}{(2j-1) \, N!^2 \, (N+j-1)! \, (N-j)!}, \qquad 
\breve{C}_{\infty}^j = \frac{4 \, (-1)^{j+1}}{(2j-1) \, \pi}. 
\end{equation}
Note that the staggered centred scheme coefficients decay faster than the coefficients 
of the collocated centred scheme. The use of such staggered schemes leads to a 
substantial gain in accuracy over the centred ones with identical $N$ and $\Delta x$ 
\citep{Fornberg1990}. 

\subsection{Relation to Lanczos' \texorpdfstring{$\sigma$}{s}-factors}

For a pure frequency of wavenumber $k$, i.e. 
$f = \e^{\ii k x}$ (with $\ii^2 = -1$ not to be confused 
with an index $i$), we have 
\begin{equation}
\left[ \left( \e^{\ii k x} \right)' \right]_i = \ii k \e^{\ii k x} \sum_{j=1}^N C_N^j 
\frac{\sin(k x_j)}{k x_j} = \sigma_c (k,N) \, \ii k \, \e^{\ii k x},
\end{equation}
where 
\begin{equation} \label{sigmaC}
\sigma_c (k,N) \eqdef \sum_{j=1}^N C_N^j \text{sinc}\!\left( \frac{j \pi k}{\kmax} \right) 
\qquad \text{with} \qquad \begin{cases} 
\kmax &\eqdef \pi / \Delta x \\
\text{sinc} &\eqdef \sin(x)/x \\
x_j &= j \Delta x = j \pi / \kmax 
\end{cases}
\end{equation}
Here, $\sigma_c(k,N)$ is called the \textit{centred $\sigma$--factor of order $2N$} 
and $\sigma_c(k,1)$ is the classical $\sigma$-factor of \citet{lanczos1988applied}. 
Note that $\sigma_c$ satisfies the properties: $\sigma_c(k,N) \in \mathbb{R}$, 
$\sigma_c(k,N) = \sigma_c(-k,N)$, $\sigma_c(0,N) = 1$ and $\sigma_c (\kmax,N) = 0$.

For a pure frequency, higher-order derivatives can concisely be defined as 
\begin{align}
\left[ \left( \e^{\ii k x} \right)^{\!(2n)} \right]_i &= (\ii k)^{2n} \, 
\e^{\ii k x_i} \, \sum_{j=1}^N C_N^j \, \sinc\!\left( \frac{j \, \pi \, k}
{2 \, \kmax} \right)^{\!2n}, \\
\left[ \left( \e^{\ii k x} \right)^{\!(2n+1)} \right]_i &= (\ii k)^{2n+1} 
\, \e^{\ii k x_i}\, \sum_{j=1}^N C_N^j \, \cos\!\left( \frac{j \, \pi \, k}
{2 \, \kmax} \right) 
\sinc\!\left(\frac{j\,\pi\,k}{2\,\kmax}\right)^{\!2n+1},
\end{align}
the elementary proof being left to the reader.

\subsection{Practical computation of high orders}\label{practical}

In practice, the coefficients $C_N^j$ are difficult to handle numerically for 
the (not so) large $N$ because of the factorials that appear within the formulas, 
leading to large round-off errors.  
To avoid both inaccuracies and integer overflows, certain computational 
techniques should be used.

One trick is to rely on the logarithm of the $\Gamma$--function. 
Using this approach, one can compute $C_N^j$ for higher values of $N$ faster 
and more accurately using the optimised $\mathtt{Gammaln}$ function \citep{AS65} 
and exploiting the relation
\begin{equation*}
\frac{(N!)^2}{(N+j)!(N-j)!} = \exp\!\left( 2 \/ \log\!\left[ \Gamma (N + 1) \right] - 
\log\!\left[ \Gamma (N + j + 1) \right] - \log \left[ \Gamma (N - j + 1) \right]  \right).
\end{equation*}
However, this approach merely delays the issue by a few orders. Indeed, even using this technique, 
we found that it was sometimes impossible to go beyond a few hundred in order without overflow.
To overcome this limitation, we found that an accurate and robust way to compute $C_N^j$ 
for any $N \in \mathds{N}$ is to exploit the relation 
\begin{equation}
\frac{(N!)^2}{(N+j)!(N-j)!} = \prod_{\ell = 1}^j \left( 1 - \frac{j}{N + l} \right) = \exp\!
\left[ \sum_{\ell = 1}^j \text{log1p}\!\left( \frac{-j}{N+l} \right) \right],
\label{log1p}
\end{equation}
where log1p$(x) \eqdef \log(1+x)$ is specifically designed \citep{Goldberg1991} 
to be accurate to machine precision as $1+x \approx 1$. To our knowledge, this 
technique never appeared before in the literature (for computed the finite 
difference weights). It allows to compute finite differences coefficients of 
large orders accurately within a reasonable time, as illustrated in Table 
\ref{tabtimestagg} for the computation of the staggered coefficients. Note 
that this concise formula performs similarly to Fornberg's efficient algorithm 
\citep{fornberg1998} (c.f. Table \ref{tabtimestagg}).

\begin{table}
\begin{tabular}{ c | c c c c}
order $N$ & direct computation	& $\Gamma$ approach & Fornberg  & log1p  \\ 
\hline
10		& 0.0079 	& 0.0081 	& 0.0088 & 0.0083 \\
80		& 0.2171	& 0.2145	& 0.2180 & 0.2187 \\
100		& overflow	& 0.2815	& 0.2819 & 0.2938 \\
250		& overflow	& 0.7886	& 0.7923 & 0.8014 \\
300		& overflow	& overflow	& 0.9658 & 0.9624 \\
\end{tabular}
\caption{Computational time (seconds) of the centred staggered $\sigma$-factor using different approaches. 
For direct computation (left column), we compute $N!$ using recursive multiplication $N(N-1)(N-2)...1$: this approach is limited to fairly low values of N. 
For large $N$, the log1p approach (last column, Eq. (\ref{log1p})) compares to Fornberg's approach, while allowing a direct computation. \label{tabtimestagg}}
\end{table}

\section{Classical forward/backward finite differences} \label{S2}

It is well known that, in certain situations, centred schemes can be unstable. 
An alternative is to consider stencils that are not symmetric. 
Fully non-centred (i.e. one-sided) schemes can then be considered.
Usual $N$-th order forward/backward schemes for $f'$ and $f''$ approximations can be written as
\begin{align}
f_i' &= \sum_{j=1}^N \, D_N^j \, \frac{f_{i \pm j} - f_i}{\pm j \, \Delta x}, 
\qquad 
f_i'' = \sum_{j=1}^N \, D_N^j \, \frac{f_{i \pm 2j} - 2 \, f_{i \pm j} + f_i}
{(j \, \Delta x)^2},
\end{align}
where
\begin{equation} \label{EQ_forward_coeffs}
D_N^j \eqdef \frac{(-1)^{j+1} \, N!}{j! \, (N-j)!}.
\end{equation}

Thus, for a pure frequency $f = \e^{\ii k x}$, this yields
\begin{align}
\left[ \left( \e^{\ii k x} \right)' \right]_i &= \ii k \, \e^{\ii k x_i} \, 
\sum_{j=1}^N \, D_N^j \, \frac{\exp (\ii k x_{\pm j})-1}{\ii k x_{\pm j}} \\
&= \ii k \, \e^{\ii k x_i} \, \sum_{j=1}^N \, D_N^j \, \left[ \frac{\sin (k x_{\pm j})}
{k x_{\pm j}} + \ii \frac{1 - \cos (k x_{\pm j})}{k x_{\pm j}} \right].
\end{align}
The latter relation suggests the introduction of a \textit{forward/backward $\sigma$-factor of 
order $N$} which, keeping the notation introduced in \eqref{sigmaC}, is defined by 
\begin{equation}
\sigma_{\mathsc{f/b}} (k,N) \eqdef \sum_{j=1}^N \, D_N^j \, \left[ \sinc \left( 
\frac{j \pi k}{\kmax} \right) \pm \ii \frac{1 - \cos (j \pi k / \kmax)}{j \pi k 
/ \kmax} \right].
\end{equation} 
This $\sigma$-factor verifies $\sigma_{\mathsc{f/b}}(k,N)\in\mathds{C}$, 
$\sigma_{\mathsc{f/b}}(-k,N)=\overline{\sigma_{\mathsc{f/b}}(k,N)}$, 
$\sigma_{\mathsc{f/b}}(0,N)=1$, $\sigma_{\mathsc{f/b}}(2\/\kmax,N)=0$ and 
$\Re \lbrace \sigma_{\mathsc{f/b}}(\kmax,N) \rbrace = 0$, ( $\overline{(.)}$ denoting 
the complex conjugate of $(.)$).  

One can notice that, as $N\rightarrow\infty$, the coefficients $D_N^j$ do not 
converge since 
\begin{equation}
\left| D_N^j \right| \sim N^j \, / \, j!
\qquad \text{as} \qquad N \rightarrow \infty.
\end{equation}
Consequently, this approach is restricted to relatively low orders $N$. To overcome these 
limitations, we propose a novel non-centred scheme.

\begin{figure}
\centering
\includegraphics[width=0.8\linewidth]{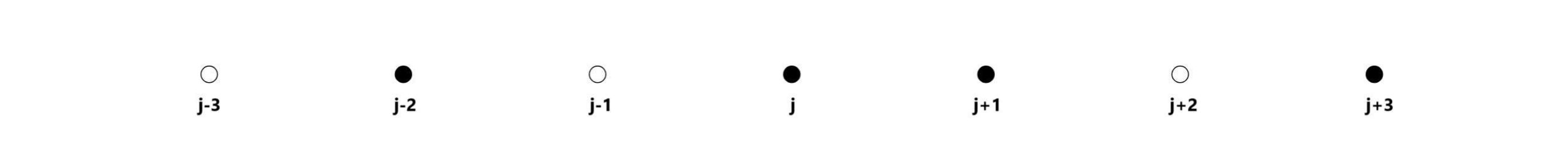}
\caption{Forward-first zigzag scheme stencil used to compute approximation at the $j$-th node. 
Dark points are used, while hollow points are not.}
\label{stencilF}
\end{figure}
\begin{figure}
\centering
\includegraphics[width=0.8\linewidth]{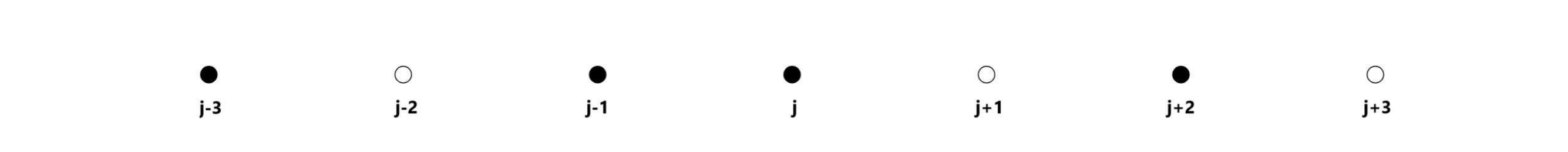}
\caption{Backward-first zigzag scheme stencil used to compute approximation at the $j$-th node. }
\label{stencilB}
\end{figure}

\section{Zigzag finite differences} \label{S3}

Rather than solely relying on forward or backward schemes, we consider here a hybrid approach 
combining both. More specifically, high-order non-centred schemes can be constructed alternating 
forward and backward differences, starting with either a forward or backward difference.
We call these hybrid schemes {\em forward/backward-first zigzag schemes}.

\subsection{Definition of zigzag schemes for regular grids} \label{ZZ}

The figures \ref{stencilF} and \ref{stencilB} illustrate how the stencils of 
such methods look like. 
A $N$th order zigzag scheme approximating $f'$ and $f''$ is defined by
\begin{align}
f_i' &= \sum_{j=1}^N \, Z_N^j \, \frac{f_{i \mp j (-1)^j} - f_i}{\mp (-1)^j 
\, j \, \Delta x}, \\
f_i'' &= \sum_{j=1}^N \, Z_N^j \, \frac{f_{i \mp 2j(-1)^j} - 2 \, f_{i \mp 
j (-1)^j} + f_i}{(j \, \Delta x)^2},
\end{align}
where the coefficients $Z_N^j$ are given analytically by the formulae  
\eqref{defZN1}--\eqref{defZNj} (see Table \ref{TabCoeffsZigzag} for the numerical 
values of the first ones).

\begin{figure}
\begin{NiceTabular}{|l||*{8}{c|}}[cell-space-limits=2pt]\hline
\backslashbox{$N$}{$j$}
&\makebox[2em]{1}&\makebox[2em]{2}&\makebox[2em]{3}
&\makebox[2em]{4}&\makebox[2em]{5}&\makebox[2em]{6}&\makebox[2em]{7}&\makebox[2em]{8}\\\hline\hline
1 & 1 &&&&&&&\\ \hline
2 & $\frac{2}{3}$ & $\frac{1}{3}$&&&&&&\\ \hline
3 & 1 & $\frac{1}{5}$ & $-\frac{1}{5}$ &&&&&\\ \hline
4 & $\frac{4}{5}$ & $\frac{2}{5}$ & $-\frac{4}{35}$ & $-\frac{3}{35}$ &&&& \\\hline
5 & 1 &$\frac{2}{7}$&$-\frac{2}{7}$&$-\frac{1}{21}$&$\frac{1}{21}$&&&\\ \hline
6 & $\frac{6}{7}$ &$\frac{3}{7}$&$-\frac{4}{21}$&$-\frac{1}{7}$&$\frac{2}{77}$&$\frac{5}{231}$&&\\ \hline
7 & 1 &$\frac{1}{3}$&$-\frac{1}{3}$&$-\frac{1}{11}$&$\frac{1}{11}$&$\frac{5}{429}$&$-\frac{5}{429}$&\\ \hline
8 & $\frac{8}{9}$ &$\frac{4}{9}$&$-\frac{8}{33}$&$-\frac{2}{11}$&$\frac{8}{143}$&$\frac{20}{429}$&$-\frac{8}{1287}$&$-\frac{7}{1287}$\\ \hline
\end{NiceTabular}
\caption{First coefficients $Z_N^j$ of the zigzag schemes.} \label{TabCoeffsZigzag}
\end{figure}

The coefficients $Z_N^j$ were obtained using a computer algebra system, low-order 
ones are given in Table \ref{TabCoeffsZigzag}. From this table, the general 
formula \eqref{defZN1}--\eqref{defZNj} for $Z_N^j$ was inferred and validated 
via symbolic computations for orders up to five thousands\footnote{Despite our 
efforts, we did not manage to prove the formula rigorously, nor to simplify it.}
\begin{align}
Z_N^1 &= \frac{2 \, N!}{(N+1)!} \, \floor*{\frac{N + 1}{2}}, \label{defZN1}\\
Z_N^2 &= \frac{2 \, N!}{(N+2)!} \, \floor*{\frac{N}{2}} \, \ceil*{\frac{N+1}{2}}, \\
Z_N^j &= (-1)^{\floor*{j/2}} \frac{2N + j - 2 - (j+2)(-1)^N}
{(j/2 - 1/2)!(N+j)!} \frac{(j-2)! N!}{(j/2 - 3/2)!} \floor*{\frac{N}{2}} 
\ceil*{\frac{N+1}{2}} \\ 
&\times \prod_{\ell=1}^{j/2 - 3/2} \floor*{\frac{N - 2l -1}{2}} \ceil*{\frac{N + 2l + 1}{2}} 
\qquad j\ \text{odd}\geqslant 3, \\
Z_N^j &= (-1)^{1 + j/2} \frac{j! N!}{(j/2)!^2 \, (N+j)!} \floor*{\frac{N}{2}} 
\ceil*{\frac{N+1}{2}} \ceil*{\frac{N-3}{2}} \ceil*{\frac{N+3}{2}} \\
&\times \prod_{\ell = 1}^{j/2 -2} \floor*{\frac{N - 2l -2}{2}} \ceil*{\frac{N + 2l + 4}{2}}  
\qquad j\ \text{even} \geqslant 4, \label{defZNj}
\end{align}
where $\floor*{\cdot}$ and $\ceil*{\cdot}$ are, respectively, the rounding towards 
$-\infty$ and $+\infty$. Note that these schemes are not centred and, unlike the 
classical forward/backward schemes, the coefficients $Z_N^j$ converge as $N \rightarrow 
\infty$, i.e.,  
\begin{equation} \label{vanishingZZcoeffsEq}
Z_{\infty}^{2\ell}\,=\,-Z_{\infty}^{2\ell+1}\,=\,
\frac{(-1)^{\ell+1}\,(2 \ell)!}{4^\ell\,(\ell !)^2}.
\end{equation}
Furthermore, unlike their classical centred counterparts \eqref{Cjinf}, the 
$\infty$-order zigzag coefficients $Z_{\infty}^j$ decay towards zero as $j\to\infty$ since 
$(2 \ell)!/4^\ell(\ell !)^2 \sim 1/ \sqrt{\pi \/\ell}$ as $\ell\to\infty$. A visual
summary of the evolution of the coefficients' magnitude from the forward, centred and 
zigzag schemes is available in Appendix \ref{CoeffsMagnitudeAnnex}.

\subsection{Zigzag \texorpdfstring{$\sigma$}{s}-factors}

Considering a pure frequency $f = \ue^{\ui\/k\/x}$, we have 
\begin{align}
\left( \e^{\ii k x} \right)_i' &= \ii k \, \e^{\ii k x} \, 
\sum_{j=1}^N \, Z_N^j \, \frac{\exp \left( \ii k \, x_{\mp j(-1)^j} \right) -1}
{\ii k \, x_{\mp j(-1)^j}} \nonumber\\
&= \ii k \, \e^{\ii k x} \, \sum_{j=1}^N \, Z_N^j \, \left[ \frac{\sin \left( k 
\, x_{\mp j(-1)^j} \right)}{k \, x_{\mp j(-1)^j}} + \ii \, \frac{1 - \cos \left(  k \, x_{\mp j(-1)^j} \right)}{k \, x_{\mp j(-1)^j} } \right],
\end{align}
which, once again using the notations introduced in \eqref{sigmaC}, suggests the zigzag $\sigma$-factor  
\begin{equation}
\sigma_{\mathsc{z}}(k,N) \eqdef \sum_{j=1}^N \, Z_N^j \, \left[ \sinc \left( \frac{j \, \pi \, k}
{\kmax} \right) \mp (-1)^j \, \ii \, \frac{1 - \cos \left( j  \pi  k \, / \, \kmax \right)}
{j \pi  k \, / \, \kmax} \right],
\end{equation}
where the upper (resp. lower) sign corresponds to the forward-first (resp. backward-first) scheme. 

A closed-form of the $\infty$-order zigzag $\sigma$-factor can be obtained. For the 
forward-first (the backward-first being the complex conjugate of the forward-first), 
we have 
\begin{align}
\sigma_{\mathsc{z}}(k, \infty) &=  \sum_{\ell = 1}^{\infty} \frac{4 (-1)^{\ell +1} 
(2 \ell -2)!}{4^\ell (\ell -1)!^2} 
\frac{\exp ((2 \ell -1) \ii \pi k \, / \, \kmax ) -1}{(2 \ell -1) \ii \pi k \, / \, \kmax} \nonumber\\
&\quad- \sum_{\ell = 1}^{\infty} \frac{ (-1)^{\ell +1} (2 \ell)!}{4^\ell (\ell !)^2} \frac{\exp (-2 \ell \ii 
\pi k \, / \, \kmax ) -1}{2 \ell \ii \pi k \, / \, \kmax} \nonumber\\
&= \frac{\kmax}{\ii \pi k} \log \left( \frac{\exp (\ii \pi k \, / \, \kmax) + \sqrt{1 + \exp (2 \ii \pi k \, 
/ \, \kmax)}}{1 + \sqrt{1 + \exp (-2 \ii \pi k \, / \, \kmax)}} \right),
\end{align}
where we have exploited the two relations 
\begin{gather}
\sum_{\ell = 1}^{\infty} \, \frac{4 (-1)^{\ell + 1} (2 \ell - 2)!}{4^{\ell} (\ell !)^2} = \text{arcsinh}(X) 
= \log \left( X + \sqrt{1 + X^2} \right), \\
\sum_{\ell = 1}^{\infty} \frac{(-1)^{\ell +1} (2 \ell)!}{4^\ell (\ell !)^2} \frac{X^{2 \ell}}{2 \ell} 
= \log \left( \frac{1 + \sqrt{1 + X^2}}{2} \right).
\end{gather}

\subsection{Zigzag schemes for staggered grids}

Zigzag schemes can be also considered for staggered grids, 
the staggered variant of these zigzag schemes being obviously
\begin{equation}
f_i' = \sum_{j=1}^N \, \breve{Z}_N^j \, \frac{f_{i \mp (j-1/2)(-1)^j} - f_i}
{\mp (-1)^j \, (j-1/2) \Delta x},
\end{equation}
where the coefficients $\breve{Z}_N^j$ are defined by  
\begin{equation}
\breve{Z}_N^j \eqdef \frac{(-1)^{\ceil*{1 + j/2}} \, (2N-1)!}{(2j-1) \, 8^{N-1} \, 
(N-1)!} \, \left( \floor*{\frac{N+j-1}{2}}! \, \floor*{\frac{N-j}{2}}! \right)^{-1}.
\end{equation}
Again, this formula was inferred from numerically computed 
coefficients, then verified using a computer algebra system 
up to order one thousand.

\section{Stability analysis on linear advection equation} \label{S4}
\label{vonNeum}
In order to compare the performances of the zigzag schemes to the classical centred and 
forward/backward ones, we conduct a stability analysis\footnote{A comparison of the numerical 
diffusion introduced by these schemes is also available in the Annex \ref{AnnexNumDiff}.} 
for the reference transport equation
\begin{equation} \label{T}
u_t\/+\/c\,u_x\,=\,0, \qquad u(0,x)\,=\,u_0(x),
\qquad (t,x) \in \mathds{R}_+ \times \mathds{T},
\end{equation}
where $\mathds{T}$ is a periodic, one-dimensional domain and $c \in \mathds{R}$ is a 
wave celerity. The exact solution of this equation is $u(t,x) = u_0(x - ct)$. 
Looking for solutions of the form $u(t,x) = U(t)\,\ue^{\ui\/k\/x}$ (with $k \in \mathds{R}$) 
and using a finite-difference scheme for the $x$-derivative, approximation of Eq. \eqref{T} 
yields the semi-discrete equation 
\begin{equation} \label{sdT}
\frac{\ud\/U}{\ud\/t}\/+\/\ui\,k\,c\,\sigma(k)\,U\,=\,0,
\end{equation}
where $\sigma$ denotes the $\sigma$-factor of the scheme used to approach the spatial derivative. 
Note that the scheme is centred,  $\sigma(k)\in\mathds{R}$. The temporal variable $t$ remains to 
be discretised and we investigate below several possibilities. 

\begin{figure}
\centering
\includegraphics[scale=0.6]{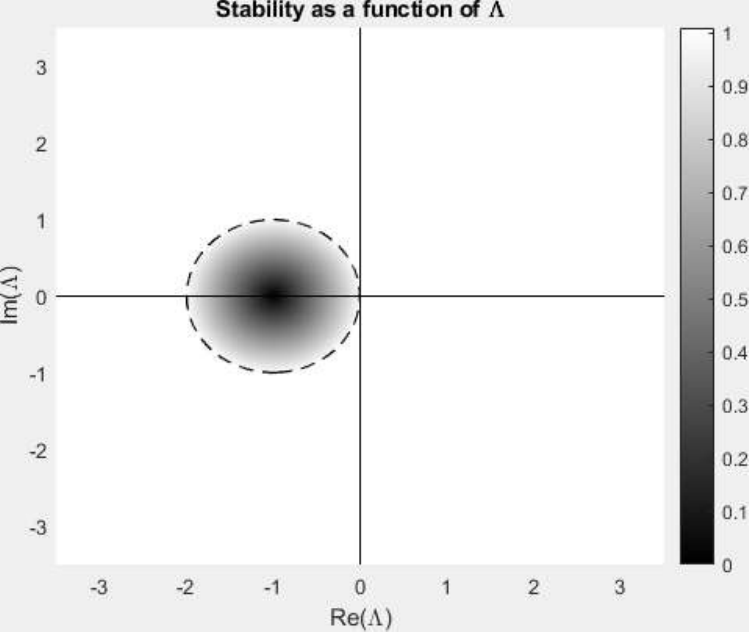}
\caption{\centering Stability region as a function of $\Lambda = -\ii \kappa \pi \lambda \sigma(\kappa)$ in the complex plane. The darker area is where the scheme is stable, and the dotted line represent the critical values where $|G| = 1$.}\label{figeul1}
\end{figure}
\subsection{Euler explicit scheme} 

We start with the first-order explicit Euler scheme $\ud\/U(t)/\ud\/t \approx [U(t+\Delta t) 
- U(t)]/\Delta t$ which, after substitution in \eqref{sdT}, leads to the fully discrete equation 
\begin{equation} \label{E1}
\frac{U(t+\Delta t) - U(t)}{\Delta t}\/ +\/ \ii\,k \, c \, \sigma(k) \, U(t)\, = \, 0.	
\end{equation}

With $U = a \exp{(\ii \omega t})$ ($a \in \mathds{R}$), one obtains
$\ue^{\ui\/\omega\/\Delta t} = 1 - \ui k \, \Delta t \, c \sigma(k)$. Introducing the \hypertarget{amp}{amplification factor} $G\eqdef U(t+ \Delta t)/U(t) = \text{e}^{\ii \omega 
\Delta t}$, the von Neumann stability condition \citep{IsaacsonKeller1994} $||G||_{\infty} \leqslant 1$ shows 
that the solution is bounded if and only if (for any wavenumber $k \in \mathds{R}$) 
\begin{equation} \label{VNE1}
| 1 - \ii \Delta t \, c k \sigma(k) | \leqslant 1.
\end{equation}
\hypertarget{notations}{Let be} the Nyquist wavenumber $\kmax \eqdef \pi / 
\Delta x$, the scaled wavenumber $\kappa = k / \kmax$ and the ``CFL-like''
stability-number $\lambda \eqdef c \, \Delta t / 
\Delta x$. 
Then \eqref{VNE1} is rewritten 
\begin{equation}\label{EulerVonNeumCond}
|1 - \ii \pi \lambda \kappa \sigma(\kappa)| \leq 1, \qquad \forall |\kappa| 
\leq 1,
\end{equation}
this condition allowing the computation of the stability region of the general 
explicit Euler scheme of order $1$ (Figure \ref{figeul1}).  When centred schemes are
considered, $\sigma$ is real and the condition \eqref{EulerVonNeumCond} yields
\begin{equation}
(\pi \lambda \kappa \sigma)^2 \leq 0,
\end{equation}
\noindent
which is never satisfied. 
In other words, this means that centred schemes are unconditionally 
unstable, a very well-known result \citep{leveque2002}. 

\subsubsection{Forward finite differences scheme of order 1}

With the spatial derivative approximated by a first-order forward 
finite difference scheme, the corresponding $\sigma$-factor is 
\begin{equation} \label{FD1}
\sigma(\kappa) = \frac{\text{e}^{\ii \pi \kappa}-1}{\ii \pi \kappa}.
\end{equation}
Substituting \eqref{FD1} into the condition \eqref{EulerVonNeumCond}, 
we obtain the stability condition
\begin{equation} 
\forall \kappa \in [-1,1], \quad |1 + \lambda - \lambda 
\text{e}^{\ii \pi \kappa}| \leq 1 \quad \implies 
\quad \lambda (1 + \lambda) \leq 0 \quad \implies
\quad -1 \leq \lambda \leq 0. \label{cond1}
\end{equation}
In order to satisfy the above condition, the wave celerity $c$ needs to 
be negative (when the time increases). This means that the scheme must be 
``upwind" \citep{leveque2002}. Condition \eqref{cond1} also implies that $|\Delta 
x / \Delta t| \geq c$, meaning that the ``numerical speed'' must be greater 
than the physical speed $c$. (If the celerity of the wave is positive, then 
a backward scheme must be considered instead.) The corresponding stability 
region is given in appendix \ref{EulerAnnex}.

\subsubsection{Forward finite differences scheme of order 2} 

For this scheme, we have
\begin{equation} \label{FD2}
\sigma(\kappa) = \ii \frac{\text{e}^{2\ii \pi \kappa}-4\text{e}^{\ii \pi \kappa}+3}
{2\pi \kappa},
\end{equation}
and the substitution of \eqref{FD2} into \eqref{EulerVonNeumCond} yields 
\begin{equation} \label{cond2}
\left\lvert 1 + \frac{\lambda}{2}\left( \text{e}^{2\ii \pi \kappa}-4\text{e}^{\ii 
\pi \kappa}+3 \right) \right\rvert \leq 1.
\end{equation}
After doing computations similar to those done in \eqref{cond1}, we show that \eqref{cond2} implies $\lambda = 0$. 
Hence, we can conclude that the second-order purely forward scheme (with 
first-order explicit Euler temporal scheme) is stable only for some wavenumbers, 
but never for all wavenumbers. Therefore,  it is unconditionally unstable. 

\subsubsection{Zigzag finite differences scheme of order 2} 

Consider now the zigzag forward-first finite difference scheme of order $2$ 
(recall that the first-order zigzag scheme coincides with the first-order 
purely forward scheme); this scheme can be defined via its $\sigma$-factor
\begin{equation}
\sigma(\kappa) = \ii \frac{\text{e}^{-2\ii \pi \kappa}-4\text{e}^{\ii \pi 
\kappa}+3}{6 \pi \kappa}. 
\end{equation}
Substituting $\sigma$ into \eqref{EulerVonNeumCond} gives the condition 
\begin{equation}
\left\lvert 1 + \frac{\lambda}{6} \left( \text{e}^{-2\ii \pi \kappa}-
4\text{e}^{\ii \pi \kappa}+3 \right) \right\rvert \leq 1.
\end{equation}
After some computations we can show that this also implies $\lambda = 0$.
Hence, as for the pure forward scheme, the zigzag scheme of order two 
(with explicit Euler) is unconditionally unstable (c.f. Appendix \ref{EulerAnnex}) 
due to the lack of stability for low wavenumbers.

\subsection{Second-order Runge--Kutta scheme} 

\begin{figure}
\centering
\includegraphics[scale=0.6]{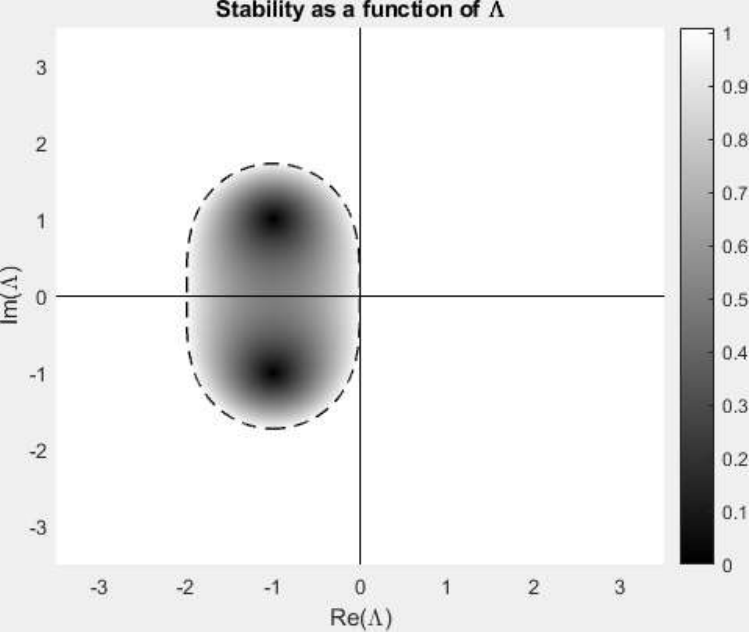}
\caption{Stability region in function of $\Lambda = -\ii \kappa \pi 
\lambda \, \sigma(\kappa)$ in the complex plane. The darker area contains the values of $\Lambda$ for which the scheme is stable. The dotted line represent the critical values where $|G| = 1$.}\label{figRK2stab}
\end{figure}

To increase the order of the time scheme, we consider the second-order 
Runge--Kutta method given by the Butcher's tableau in Table \ref{ButcherTab}.
\renewcommand{\arraystretch}{1.25}
\begin{table}
\begin{tabular}{c | c  c}
$0$ 		& 			&  		\\ 
$\frac{1}{2}$ & $\frac{1}{2}$ 	& 		\\ \hline
		& $0$		& $1$
\end{tabular}
\caption{Butcher's tableau of the second-order Runge--Kutta's method used.} \label{ButcherTab}
\end{table}
\renewcommand{\arraystretch}{1}
Doing so, the semi-discrete equation \eqref{sdT} becomes the fully discrete equation 
\begin{equation} \label{RK2_1}
U(t + \Delta t) = \left[ 1 + \Delta t \, A + \frac{1}{2} (\Delta t \, A)^2 \right] \, U(t),
\qquad A \eqdef - \ii k \, c \, \sigma(k),
\end{equation}
with the von Neumann stability condition 
\begin{equation} \label{RK2}
\left\lvert 1 + \Delta t \, A + \frac{1}{2} (\Delta t \, A)^2 \right\rvert \leq 1 \qquad \forall\,k.
\end{equation}
Using the notations introduced in \hyperlink{notations}{\textit{section 1}}, 
the condition \eqref{RK2} becomes 
\begin{equation} \label{VonN2}
\left\rvert 1 - \ii \pi \kappa \lambda \, \sigma - \frac{1}{2} (\pi \kappa \lambda \, \sigma)^2 \right\rvert \leq 1 \qquad \forall\, |\kappa| \leq 1,
\end{equation}
the corresponding stability region being depicted in figure \ref{figRK2stab}.
\noindent
For centred schemes, $\sigma$ is real and this condition yields
\begin{equation}
(\pi \lambda \kappa \sigma)^4 \leq 0.
\end{equation}
\noindent
It cannot be satisfied 
meaning that centred schemes are unconditionally unstable, like the 
explicit Euler scheme. \\

\subsubsection{Forward finite differences scheme of order 1} 

\hypertarget{FD1}{}
The $\sigma$-factor corresponding to this case, i.e.
\begin{equation}
\sigma(\kappa) = \frac{\text{e}^{\ii \pi \kappa}-1}{\ii \pi \kappa},
\end{equation}
substituted into \eqref{VonN2} gives the condition 
\begin{align}
\left\lvert 1 - \lambda \left( \text{e}^{\ii \pi \kappa} - 1 \right) + 
\frac{1}{2} \lambda^2 \left( \text{e}^{\ii \pi \kappa}-1 \right)^2 \right\rvert 
\leq 1 \quad &\implies \quad -1 \leq \lambda^3 + 2 \, \lambda^2 + 2 \, \lambda \leq 0. \\
 &\implies \quad \lambda \in [-1,0].
\end{align} 
Hence, the forward finite differences scheme of order one is conditionally 
stable, requiring the stability-number $\lambda$ to be in $[-1,0]$. We note that, like in 
the case of the explicit Euler scheme, the scheme must be ``upwind'' to be stable.

\subsubsection{Forward finite differences scheme of order 2} 

In this case, the $\sigma$-factor 
\begin{equation}
\sigma(\kappa) = \ii \frac{\text{e}^{2\ii \pi \kappa}-4\text{e}^{\ii \pi 
\kappa}+3}{2\pi \kappa},
\end{equation} 
substituted into the von Neumann stability condition \eqref{VonN2} yields
\begin{equation}\label{condFD2RK2} 
\left\lvert 1 + \frac{\lambda}{2} \left( \text{e}^{2 \ii \pi \kappa} - 4 
\text{e}^{\ii \pi \kappa} + 3  \right)+ \frac{\lambda^2}{8}\left( \text{e}^{2 
\ii \pi \kappa} - 4 \text{e}^{\ii \pi \kappa} + 3  \right)^2 \right\rvert 
\leq 1.
\end{equation} 
A study of the condition \eqref{condFD2RK2} shows that $\lambda \in [-\frac{1}{2}, 0]$.
Thus, the second-order forward finite differences scheme is conditionally 
stable.
Similarly, for the first-order scheme, the strictest constraint on the stability-number 
is imposed by the Nyquist wavenumber 
(see Annex \ref{RK2Annex}).
Note that the stability strip is half as wide as that of the first-order.

\subsubsection{Forward finite differences scheme of order $3$} 

For this scheme,  the $\sigma$-factor is
\begin{equation} \label{2}
\sigma(\kappa) = - \frac{\ii}{6 \pi \kappa} \left( 2\text{e}^{3 \ii \pi 
\kappa} - 9 \text{e}^{2 \ii \pi \kappa} + 18 \text{e}^{\ii \pi \kappa} - 11 
\right),
\end{equation}
which gives
\begin{equation}
\left\lvert 1 - \frac{\lambda}{6}  \left( 2\text{e}^{3 \ii \pi \kappa} - 9 
\text{e}^{2 \ii \pi \kappa} + 18 \text{e}^{\ii \pi \kappa} - 11 \right) + 
\frac{\lambda^2}{72} \left( 2\text{e}^{3 \ii \pi \kappa} - 9 \text{e}^{2 \ii \pi \kappa} + 18 \text{e}^{\ii \pi \kappa} - 11  \right)^2  \right\rvert \leq 1.
\end{equation}
A study of this condition reveals that $\lambda = 0$.
This scheme is thus unstable.
Moreover, numerical observations seem to show that upwind schemes of 
orders $3$ and above are unstable.

\subsubsection{Zigzag finite differences scheme of order 2}  

The $\sigma$-factor corresponding to this scheme is
\begin{equation}
\sigma(\kappa) = \ii \frac{\text{e}^{-2\ii \pi \kappa}-4\text{e}^{\ii \pi \kappa}+3}{6 \pi \kappa},
\end{equation}
which once injected in \eqref{VonN2} yields
\begin{equation} \label{condZZorder2}
\left\lvert 1 + \frac{\lambda}{6} \left( \text{e}^{-2\ii \pi 
\kappa}  - 4 \text{e}^{4\ii \pi \kappa} + 3\right) + \frac{\lambda^2}{72}  
\left( \text{e}^{-2\ii \pi \kappa}  - 4 \text{e}^{4\ii \pi \kappa} + 
3\right)^2 \right\rvert \leq 1.
\end{equation}
This condition implies $\lambda P(\lambda) \leq 0$, with $P$ defined by
\begin{align*}
P(\lambda) \eqdef & \ \ \lambda^3 \left[64 \cos\left(\kappa  \pi \right)^4+80 
\cos\left(\kappa  \pi \right)^3+105 \cos\left(\kappa  \pi \right)^2+50 \cos 
\! \left(\kappa  \pi \right)+25\right] \\
&+ \lambda^2 \left[-96 \cos\left(\kappa  \pi \right)^3+36 \cos\left(\kappa  
\pi \right)^2+60 \right] \\
&+  \lambda \left[ \cos \left(\kappa  \pi \right)^2-144 \cos \! \left(\kappa 
\pi \right)+72 \right] + 216.
\end{align*}
The polynomial $P$ possesses a global extremum equal to $324 \lambda^3 + 216$ 
at $\kappa = 0$ for any $\lambda$. This extremum is null for the critical 
$\lambda_c = -\sqrt[3]{2} \left( \sqrt[3]{3} \right)^2 / 3 \approx -0.8736$. 
It corresponds to a minimum of $P$ if and only if $\lambda \geq \lambda_c$, 
and it is a maximum if and only if $\lambda \leq \lambda_c$ . Hence, this 
scheme is conditionally stable, under the condition $\lambda_c \leq \lambda 
\leq 0$. \\
The condition on the stability-number for this scheme is less restrictive than that 
of the forward-scheme for the same order. Indeed, the stability strip 
$[-1/2, 0]$ is included in $\left[\lambda_c, 0\right]$. Note that this time, 
the most restrictive stability condition is imposed by the null wavenumber. \\

\subsubsection{Zigzag finite differences scheme of order 3} The $\sigma$-factor of this scheme is
\begin{equation}
\sigma(\kappa) = - \frac{2 \text{e}^{3 \ii \pi \kappa} + 3 \text{e}^{-2\ii 
\pi \kappa} - 30 \text{e}^{\ii \pi \kappa} + 25}{30 \ii \pi \kappa}
\end{equation}
Which after substitution in the von Neumann condition \eqref{VonN2} gives
\begin{equation*}
\left\lvert 1+\frac{\lambda}{30}  \left(2 \,{\mathrm e}^{3 \,\ii 
\kappa  \pi}+3 \,{\mathrm e}^{-2 \,\ii \kappa  \pi}-30 \,{\mathrm 
e}^{\ii \kappa  \pi}+25\right)+\frac{\lambda^{2}}{1800} \left(2 
\,{\mathrm e}^{3 \,\ii \kappa  \pi}+3 \,{\mathrm e}^{-2 \,\ii 
\kappa  \pi}-30 \,{\mathrm e}^{\ii \kappa  \pi}+25\right)^{2} 
\right\rvert \leqslant 1.
\end{equation*}
A similar procedure to the one followed for the condition \eqref{condZZorder2}
gives $-\frac{15}{14} \leq \lambda \leq 0$.
This scheme remains conditionally stable and it is the first one allowing 
$|\lambda| > 1$ thus far. Note that for the same order, the zigzag-scheme is 
conditionally stable whilst the equivalent forward-scheme is not. \\

\subsubsection{Zigzag finite differences scheme of order 4} For this scheme,
\begin{equation}
\sigma (\kappa) = \-\frac{336 \,{\mathrm e}^{\ii \kappa  \pi}-245-84 
\,{\mathrm e}^{-2 \,\ii \kappa  \pi}-16 \,{\mathrm e}^{3 \,\ii \kappa  \pi}+9 
\,{\mathrm e}^{-4 \,\ii \kappa  \pi}}{420 \ii \kappa  \pi}.
\end{equation}
Substituting this $\sigma$-factor in \eqref{VonN2} yields $\lambda = 0$.
This scheme is thus unconditionally unstable. Numerically, we observe that 
the third-order seems to be the highest order allowing the stability of the 
zigzag-schemes while using (RK2). \\

\subsection{Summary on the stability of higher order schemes} \label{Table}

For orders larger than three, the computations of the stability 
regions become highly tedious. Hence, for these higher orders, we 
proceeded by computing the stability-number numerically with an accuracy 
of $\pm 10^{-4}$. In the case of staggered zigzag schemes, we 
extrapolated the results for the infinite order from the order 300. 
A summary of all the critical stability-number values can be found in the 
tables \ref{EulerStabTab},\ref{RK2StabTab},\ref{RK3StabTab},\ref{RK4StabTab},
\ref{RK5StabTab},\ref{RK6StabTab} and \ref{RK7StabTab} and the corresponding
stability regions are in appendix \ref{EulerAnnex} to \ref{RK7Annex}.
Note that tables \ref{RK2StabTab}, \ref{RK5StabTab} and \ref{RK6StabTab},
highlight situations where the zigzag schemes are the only ones to be
conditionally stable.
Furthermore, these tables also reveal that zigzag schemes consistently exhibit greater 
stability for even orders. This alternating pattern is expected, as the 
stencils of zigzag schemes are asymmetrical by design. Consequently, one side 
of the stencil always incorporates more points than the other. However, the 
enhanced stability observed at even orders is unexpected. Indeed, at even 
orders, the backward side of the stencil will contain more points and thus, 
should contribute to a greater extent. This is counter-intuitive, given the 
well-known result that backward schemes are unconditionally unstable for the 
transport equation with negative velocity. \\

\begin{figure}
$\begin{array}{c|c c}
$scheme$ 		& N = 1 & N \geq 2 \\ \hline
$centred$		&  & 0 \\ 
$centred staggered$ &  & 0 \\ 
$forward$		& 1	      & 0 \\ 
$zigzag$		& 1		  & 0 \\ 
$zigzag staggered$ & 1 & 0 \\ 
\end{array} $
\caption{Highest stability-number $|\lambda|= \lvert c \, \Delta t / \Delta x \rvert$ allowed for the forward Euler integration scheme in time with spatial schemes of orders $N=1$ and $N=2$.} \label{EulerStabTab}
\end{figure}

\begin{figure}
$\begin{array}{c | c c c c}
$scheme$ 		& N = 1 & 	N = 2 & N = 3 & N \geq 4\\ \hline
$centred$ &  & 0 &  & 0\\ 
$centred staggered$ &  & 0 &  & 0 \\ 
$forward$		& 1	      		& 1/2  & 0 & 0\\ 
$zigzag$		& 1		  		& \sqrt[3]{2} (\sqrt[3]{3})^2 /3 \approx 0.8736 & 15/14 \approx 1.0714 & 0 \\ 
$zigzag staggered$ & 1.5436 & 1.2599 & 1.7099 & 0 \\ 
\end{array} $
\caption{Highest stability-number $|\lambda|= \lvert c \, \Delta t / \Delta x \rvert$ allowed for RK2 time integrator coupled with spatial schemes of orders $N=1,\; 2,\; 3$ and orders $N \geq 4$. For a given order $N$ (i.e. for a given column), the zigzag schemes enjoy the best stability constraints.} 
\label{RK2StabTab}
\end{figure}

\begin{figure}
$\begin{array}{c|c c c c c c c}
$scheme$ 		& N = 1 & N = 2 & N = 3 & N =4 & N =5 & N =6 & N = \infty \\ \hline
$centred$		&    & 1.7320 &  & 1.2622 &  & 1.0920 & 0.5533 \\ 
$centred staggered$ &  & 1.7320 &  & 1.4845 &  & 1.3949 & 1.1282 \\ 
$forward$		& 1.2563 &	0.6280 & 0 & 0 & 0 & 0 & 0\\ 
$zigzag$		& 1.2563 & 1.8845 & 1.3461 & 1.6490 & 1.3741 & 1.5727 & 1.1037\\ 
$zigzag staggered$ & 1.6791 & 2.5187 & 1.9374 & 2.3199 & 2.0112 & 2.2569 & 2.1407 \\ 
\end{array} $
\caption{Highest stability-number $|\lambda|= \lvert c \, \Delta t / \Delta x \rvert$ allowed for RK3 time integrator with spatial schemes of orders $N=1$ to  $6$, and $N = \infty$ (extrapolated from order $300$ in the case of the zigzag staggered).  For a given order $N$ (i.e. for a given column), the zigzag schemes enjoy the best stability constraints.\label{RK3StabTab}}
\end{figure}

\begin{figure}
$\begin{array}{c | c c c c c c c}
$scheme$ 		& N = 1 & N = 2 & N = 3 & N =4 & N =5 & N =6 & N = \infty \\ \hline
$centred$		&  & 2.8284 &  & 2.0612 &  & 1.7834 & 0.9035 \\ 
$centred staggered$ &  & 2.8284 &  & 2.4243 &  & 2.2779 & 1.8425 \\ 
$forward$		& 1.3926 &	0.6963 & 0 & 0 & 0 & 0 & 0\\ 
$zigzag$		& 1.3926 & 2.0889 & 1.4921 & 1.8278 & 1.5232 & 1.7433 & 1.5800\\ 
$zigzag staggered$ & 1.9122 & 2.8683 & 2.2064 & 2.6419 & 2.2904 & 2.5701 & 2.4378 \\ 
\end{array} $
\caption{Highest stability-number $|\lambda|= \lvert c \, \Delta t / \Delta x \rvert$ allowed for a RK4 time integrator with spatial 
schemes of orders $N=1$ to  $6$, and $N = \infty$ (extrapolated from order $300$ in the case of the zigzag staggered).  \label{RK4StabTab}} 
\end{figure}

\begin{figure}
$\begin{array}{c|c c c c c c c}
$scheme$ 		& N = 1 & N = 2 & N = 3 & N =4 & N =5 & N =6 & N = 7 \\ \hline
$centred$		&    & 0 &  & 0 &  & 0 &  \\ 
$centred staggered$ &    & 0 &  & 0 &  & 0 & \\ 
$forward$		& 1.6085 &	0.8042 & 0 & 0 & 0 & 0 & 0\\ 
$zigzag$ & 1.6085 & 2.4127 & 1.7234 & 2.0100 & 1.7592 & 0.4542 & 0\\ 
$zigzag staggered$ & 2.2646 & 3.3969 & 2.6129 & 3.1287 & 2.7123 & 0.7234 & 0.8121\\ 
\end{array} $
\caption{Highest stability-number $|\lambda|= \lvert c \, \Delta t / \Delta x \rvert$ allowed for a RK5 time integrator with spatial schemes of orders $N=1$ to $7$. 
For a given order $N$ (i.e. for a given column), the zigzag schemes enjoy the best stability constraints.\label{RK5StabTab}}
\end{figure}

\begin{figure}
$\begin{array}{c|c c c c c c c}
$scheme$ 	    & N = 1 & N = 2 & N = 3 & N =4 & N =5 & N =6 & N = 7 \\ \hline
$centred$		&    & 0 &  & 0 &  & 0 &  \\ 
$centred staggered$ &   & 0 &  & 0 &  & 0 & \\ 
$forward$		& 1.7767 &	0.88824 & 0 & 0 & 0 & 0 & 0\\ 
$zigzag$ & 1.7767 & 2.6650 & 1.9035 & 2.3319 & 1.9431 & 2.2240 & 0\\ 
$zigzag staggered$ & 2.4915 & 3.7373 & 2.8749 & 3.4423 & 2.9842 & 3.3488 & 3.0350\\ 
\end{array} $
\caption{Highest stability-number $|\lambda| = \lvert c \, \Delta t / \Delta x \rvert$ allowed for RK6 time integrator with spatial schemes of orders $N=1$ to $7$. 
For a given order $N$ (i.e. for a given column), the zigzag schemes enjoy the best stability constraints.\label{RK6StabTab}}
\end{figure}

\begin{figure}
$\begin{array}{c|c c c c c c c}
$scheme$ 		& N = 1 & N = 2 & N = 3 & N =4 & N =5 & N =6 & N = \infty \\ \hline
$centred$		&    & 1.7644 & & 1.2857 &  & 1.1124 & 0.5638 \\ 
$centred staggered$ &  & 1.7644  &  & 1.5122 &  & 1.4209 & 1.1493\\ 
$forward$		& 1.9770 &	0.9884 & 0 & 0 & 0 & 0 & 0\\ 
$zigzag$ & 1.9770 & 2.9656 & 2.1182 & 2.5949 & 2.1624 & 2.3282 & 1.1198\\ 
$zigzag staggered$ & 2.7440 & 4.1160 & 3.1661 & 3.791 & 3.2865 & 3.6881 & 2.5618\\ 
\end{array} $
\caption{Highest stability-number $|\lambda| = \lvert c \, \Delta t / \Delta x \rvert$ allowed for a RK7 
time integrator with spatial schemes orders $N=1$ to $6$, and $\infty$
(extrapolated from order $300$ in the case of the zigzag staggered).  
For a given order $N$ (i.e. for a given column), the zigzag schemes 
enjoy the best stability constraints. \label{RK7StabTab}}
\end{figure}

\section{Prevention of ghost solutions} \label{SectionGhostSol}

\begin{figure}
    \centering
    \renewcommand{\baselinestretch}{0.9} 

    \begin{minipage}{0.25\textwidth}
        \centering
        \resizebox{\linewidth}{!}{\includegraphics{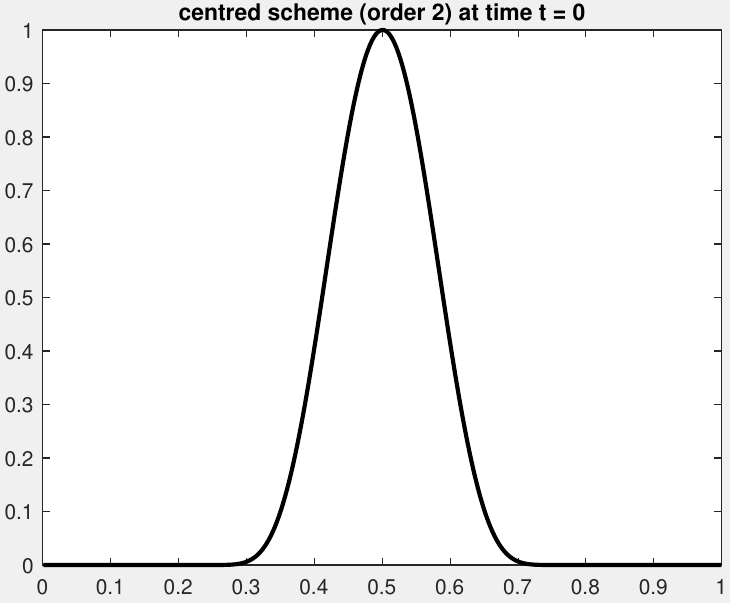}}
    \end{minipage}%
    \hspace{1mm}
    \begin{minipage}{0.25\textwidth}
        \centering
        \resizebox{\linewidth}{!}{\includegraphics{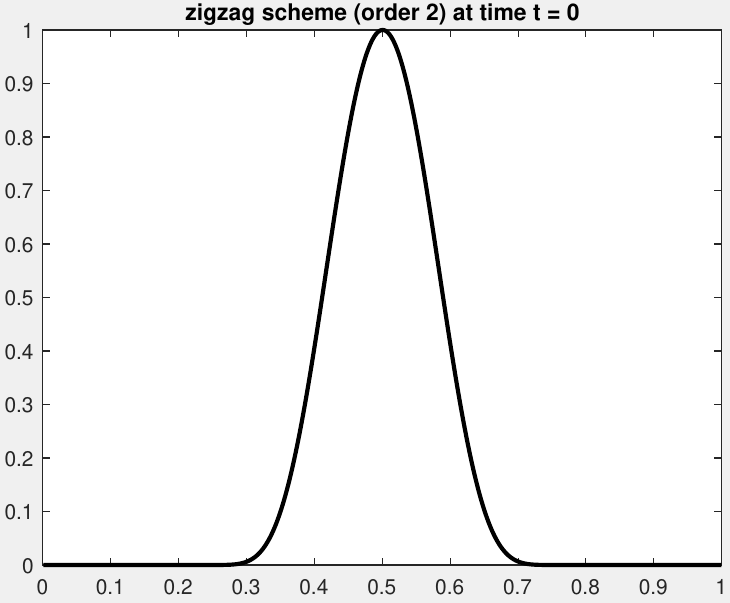}}
    \end{minipage}
    
    \vspace{-2mm} 

    \begin{minipage}{0.25\textwidth}
        \centering
        \resizebox{\linewidth}{!}{\includegraphics{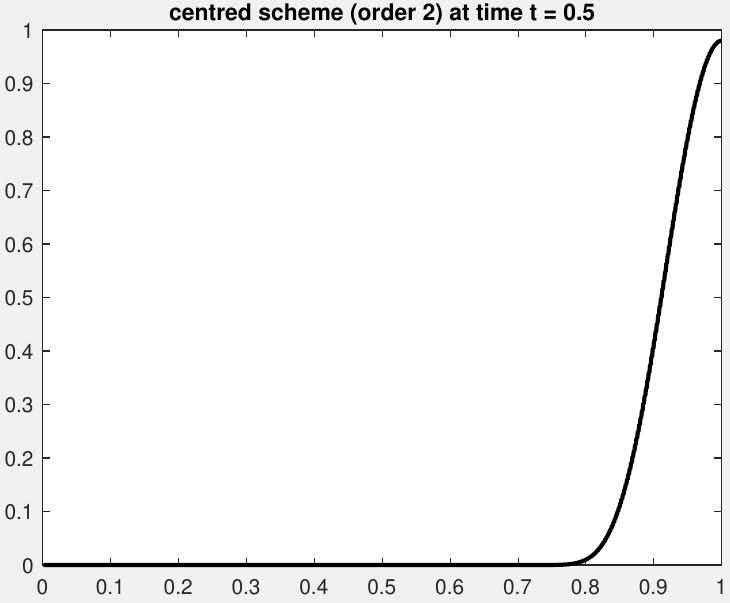}}
    \end{minipage}%
    \hspace{1mm}
    \begin{minipage}{0.25\textwidth}
        \centering
        \resizebox{\linewidth}{!}{\includegraphics{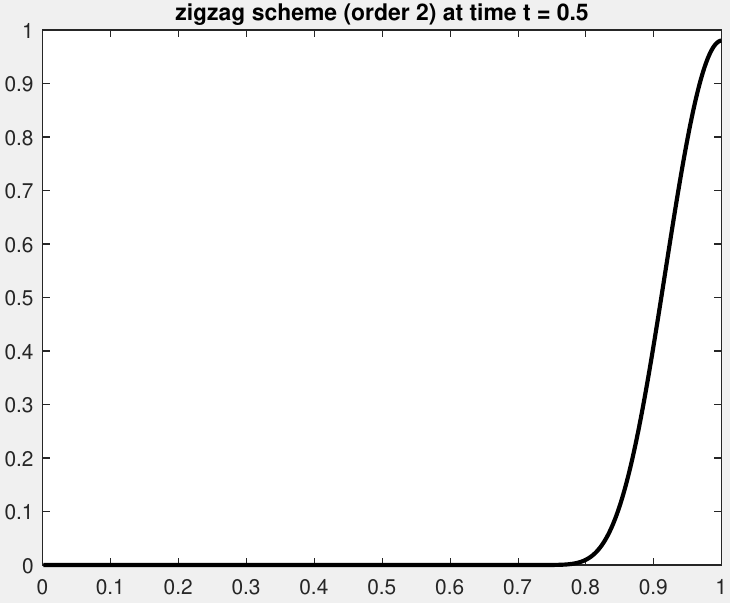}}
    \end{minipage}
    
    \vspace{-2mm}

    \begin{minipage}{0.25\textwidth}
        \centering
        \resizebox{\linewidth}{!}{\includegraphics{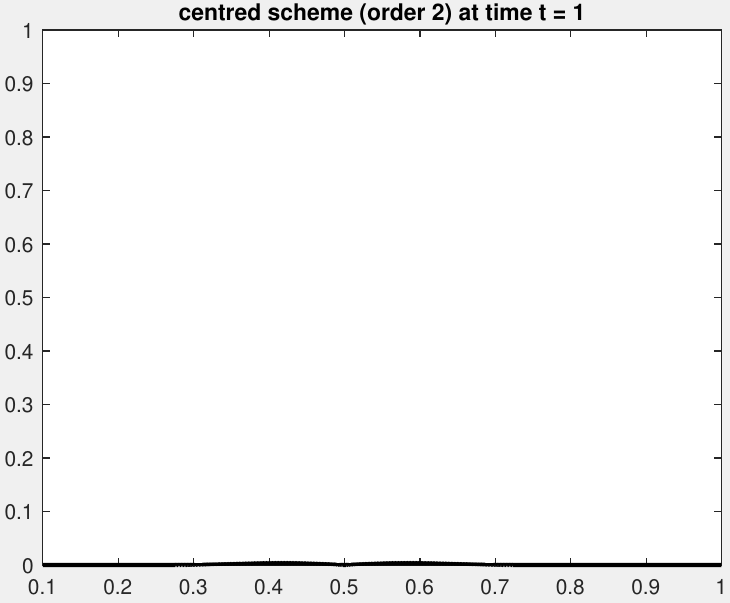}}
    \end{minipage}%
    \hspace{1mm}
    \begin{minipage}{0.25\textwidth}
        \centering
        \resizebox{\linewidth}{!}{\includegraphics{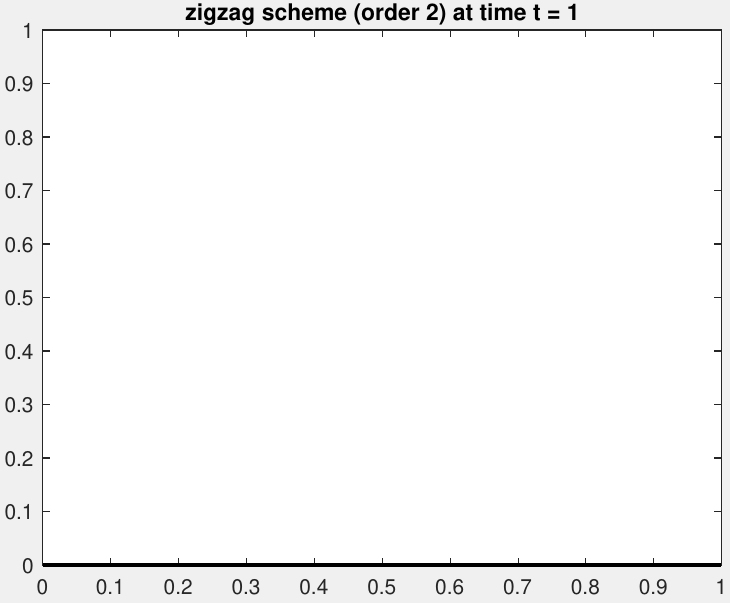}}
    \end{minipage}
    
    \vspace{-2mm}

    \begin{minipage}{0.25\textwidth}
        \centering
        \resizebox{\linewidth}{!}{\includegraphics{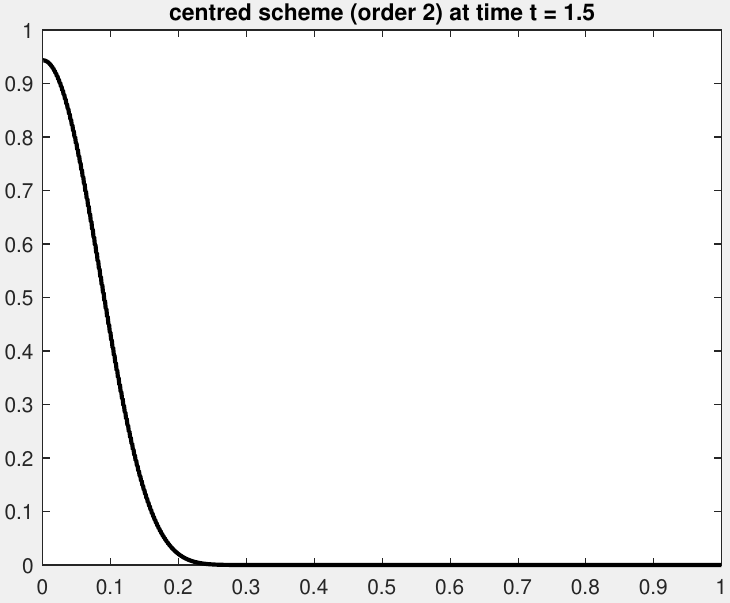}}
    \end{minipage}%
    \hspace{1mm}
    \begin{minipage}{0.25\textwidth}
        \centering
        \resizebox{\linewidth}{!}{\includegraphics{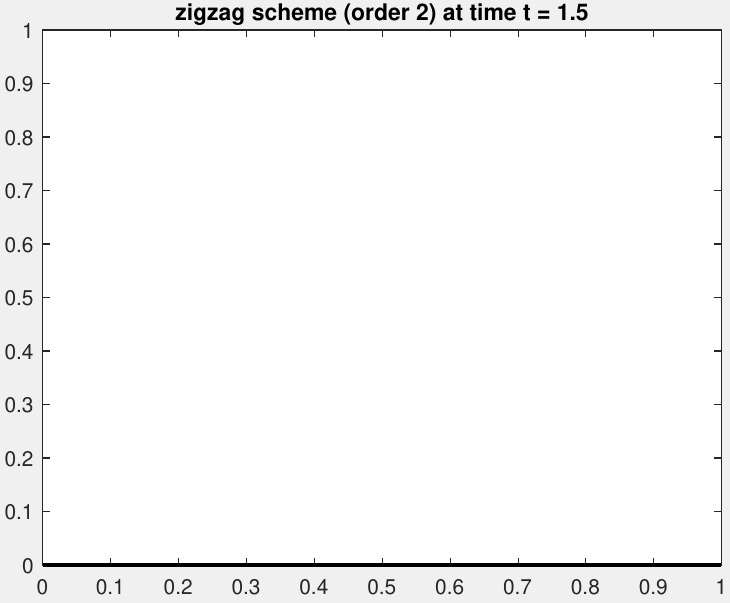}}
    \end{minipage}
    
    \vspace{-2mm}

    \begin{minipage}{0.25\textwidth}
        \centering
        \resizebox{\linewidth}{!}{\includegraphics{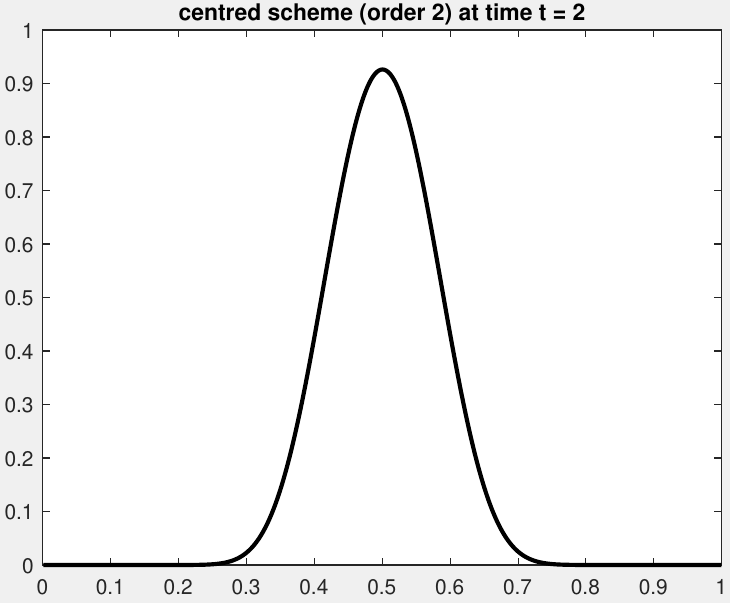}}
    \end{minipage}%
    \hspace{1mm}
    \begin{minipage}{0.25\textwidth}
        \centering
        \resizebox{\linewidth}{!}{\includegraphics{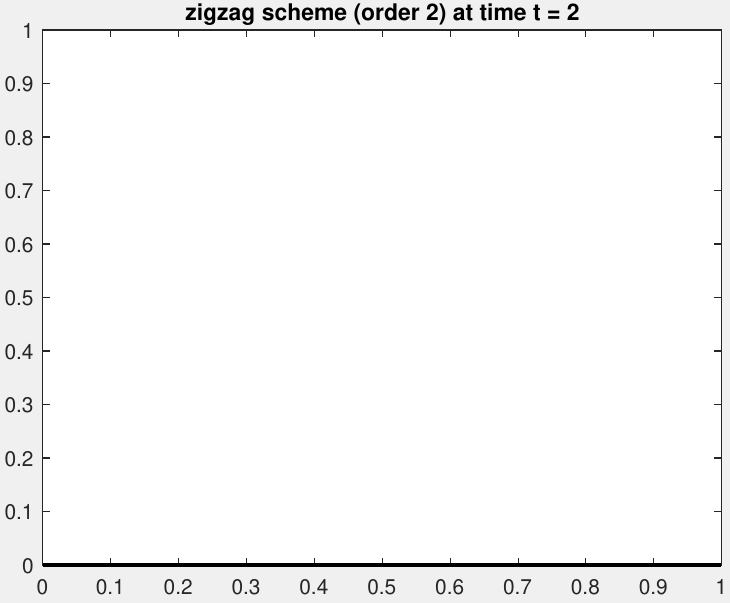}}
    \end{minipage}
    
    \vspace{-2mm}

    \begin{minipage}{0.25\textwidth}
        \centering
        \resizebox{\linewidth}{!}{\includegraphics{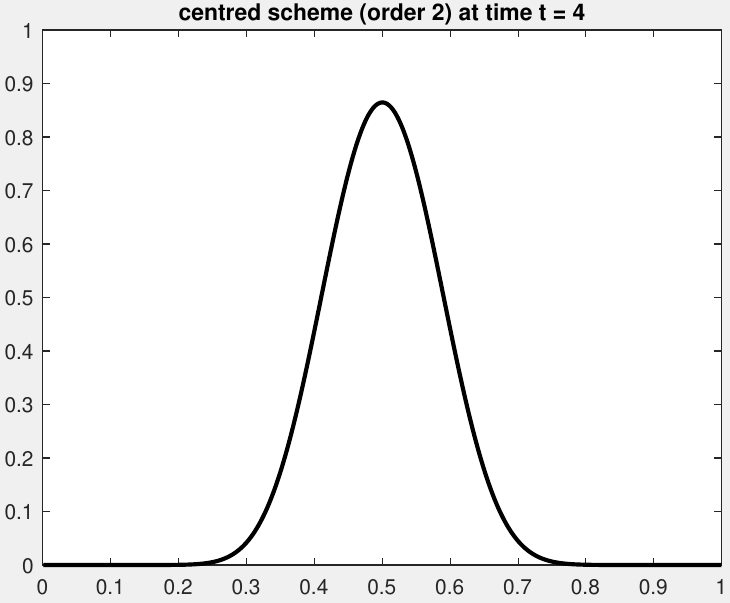}}
    \end{minipage}%
    \hspace{1mm}
    \begin{minipage}{0.25\textwidth}
        \centering
        \resizebox{\linewidth}{!}{\includegraphics{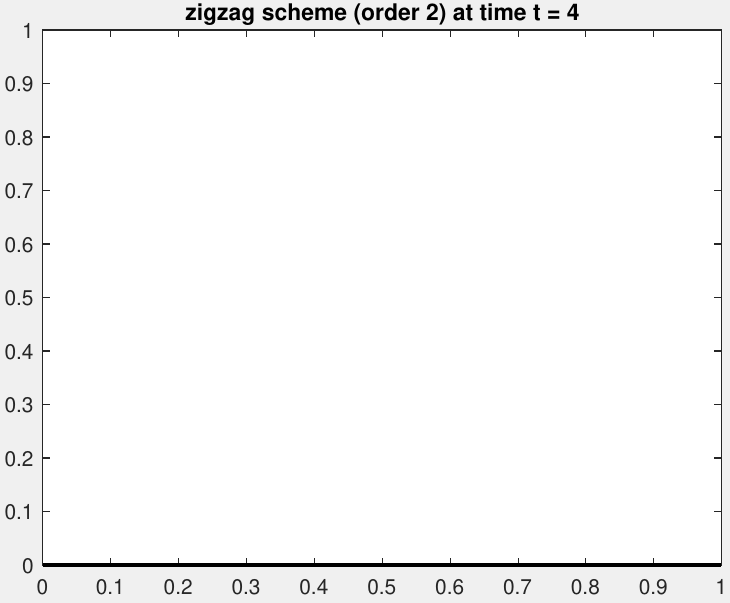}}
    \end{minipage}

    \caption{Numerical simulations of the transport equation \eqref{T} using the second order implicit centred scheme of \cite{Inglard2020} (left column) and the second order implicit zigzag backward-first scheme (right column) at times $t = 0, 0.5, 1, 1.5, 2$ and $4$ using $1000$ points to discretise the space.}
    \label{FIGghost_solution2}
\end{figure}

In the recent work, \citet{Inglard2020} showcased the formation of ghost solutions 
when solving the advection equation \eqref{T} using centred schemes with von Neumann 
boundary conditions. The term ``ghost solution'' designates the anomalous formation 
of a periodic solution, despite imposing von Neumann conditions on the boundaries of 
the domain (see Fig.  \ref{FIGghost_solution2}). 
Under the same conditions formulated in \citep{Inglard2020}, the replacement of the 
centred scheme by a one-sided scheme, or by a zigzag scheme completely prevents the 
formation of ghost solutions.
However, as mentioned in section \ref{S2}, the coefficients of the one-sided scheme diverge
as the order tends to infinity. This property limits the use of one-sided schemes to low orders.
The zigzags schemes do not suffer such a restriction. In other words, the zigzag schemes
allow the use of high-order schemes using von Neumann boundary conditions, without forming ghost solutions.
A comparison of the centred and the zigzag schemes with von Neumann conditions is provided in Fig.
 \ref{FIGghost_solution2}. \\

\section{Conclusion}

In this work, we introduced a novel class of finite differences approximations, 
referred to as {\em zigzag schemes}. These schemes employ hybrid stencils that 
are neither fully centred nor fully one-sided, leading to encouraging results. 
Specifically, as summarised in tables \ref{EulerStabTab} to \ref{RK7StabTab}, 
zigzag schemes exhibit not only a broader range of conditional stability but also, 
in some cases, a significantly higher critical stability-number, sometimes exceeding 
that of standard schemes by a factor of two or more. Moreover, the coefficients
$Z_N^j$ of the zigzag schemes vanish as the order tends to infinity (see Eq. 
\eqref{vanishingZZcoeffsEq}). The zigzag schemes thus allow the use of higher 
order schemes. 
Overall, zigzag schemes offer a compelling alternative for finite difference 
simulations. They enable faster computations at comparable accuracy levels, 
while retaining the simplicity of implementation characteristic of 
conventional schemes. 
Furthermore, a staggered variant of the zigzag schemes was presented in section 
\ref{ZZ}, alongside a closed-form expression for the stencil coefficients at 
infinite order. 

Moreover, the asymmetrical structure of the zigzag schemes prevents the formation
of some issues arising when using symmetrical (centred) schemes. This
was illustrated in section \ref{SectionGhostSol} where we showed that zigzags
schemes are able to prevent the formation of anomalous ``ghost solutions''. 

Lastly, we remark that the zigzag schemes are defined by a signature following 
the pattern $+-+-+-+$ (or $-+-+-+$). This choice was made because it seemed to 
be the most ``natural'' and ``easier'' way to blend centred and one-sided schemes 
together. However, other schemes employing other (potentially random) signatures 
could be investigated.

\bibliographystyle{abbrvnat}
\bibliography{biblio}

\begin{thebibliography}{24}
\providecommand{\natexlab}[1]{#1}
\providecommand{\url}[1]{\texttt{#1}}
\expandafter\ifx\csname urlstyle\endcsname\relax
  \providecommand{\doi}[1]{doi: #1}\else
  \providecommand{\doi}{doi: \begingroup \urlstyle{rm}\Url}\fi

\bibitem[Abide et~al.(2017)Abide, Binous, and Zeghmati]{Abide_etal2017}
S.~Abide, M.~S. Binous, and B.~Zeghmati.
\newblock An efficient parallel high-order compact scheme for the 3{D}
  incompressible {N}avier-{S}tokes equations.
\newblock \emph{Int. J. Comput. Fluid Dyn.}, 31\penalty0 (4-5):\penalty0
  214--229, 2017.
\newblock ISSN 1061-8562,1029-0257.
\newblock \doi{10.1080/10618562.2017.1326592}.
\newblock URL \url{https://doi.org/10.1080/10618562.2017.1326592}.

\bibitem[Abramowitz and Stegun(1965)]{AS65}
M.~Abramowitz and I.~A. Stegun.
\newblock \emph{Handbook of {M}athematical {F}unctions}.
\newblock Dover, 1965.

\bibitem[CORIA-CFD()]{HAllegro}
CORIA-CFD.
\newblock H-allegro. high-order finite difference code that solves the unsteady
  compressible reacting navier-stokes equations system on structured cartesian
  meshes.
\newblock URL \url{https://www.coria-cfd.fr/index.php/H-Allegro}.

\bibitem[Euler(1768)]{euler1768}
L.~Euler.
\newblock \emph{Institutiones calculi integralis}.
\newblock Acad. Imper. scientiarum, 1768.

\bibitem[Fadlun et~al.(2000)Fadlun, Verzicco, Orlandi, and
  Mohd-Yusof]{Fadlun_etal2000}
E.~A. Fadlun, R.~Verzicco, P.~Orlandi, and J.~Mohd-Yusof.
\newblock Combined immersed-boundary finite-difference methods for
  three-dimensional complex flow simulations.
\newblock \emph{J. Comput. Phys.}, 161\penalty0 (1):\penalty0 35--60, 2000.
\newblock ISSN 0021-9991,1090-2716.
\newblock \doi{10.1006/jcph.2000.6484}.
\newblock URL \url{https://doi.org/10.1006/jcph.2000.6484}.

\bibitem[Fornberg(1990)]{Fornberg1990}
B.~Fornberg.
\newblock High-order finite differences and the pseudospectral method on
  staggered grids.
\newblock \emph{SIAM J. Numer. Anal.}, 27\penalty0 (4):\penalty0 904--918,
  1990.

\bibitem[Fornberg(1998)]{fornberg1998}
B.~Fornberg.
\newblock \emph{A Practical Guide to Pseudospectral Methods}.
\newblock Cambridge Monographs on Applied and Computational Mathematics.
  Cambridge University Press, 1998.
\newblock ISBN 9780521645645.

\bibitem[Giles et~al.(2014)Giles, L\'aszl\`o, Reguly, Appleyard, and
  Demouth]{Giles_etal2014}
M.~Giles, E.~L\'aszl\`o, I.~Reguly, J.~Appleyard, and J.~Demouth.
\newblock Gpu implementation of finite difference solvers.
\newblock In \emph{2014 Seventh Workshop on High Performance Computational
  Finance}, pages 1--8, 2014.
\newblock \doi{10.1109/WHPCF.2014.10}.

\bibitem[Godlewski and Raviart(1996)]{Godlewski_etal1996}
E.~Godlewski and P.-A. Raviart.
\newblock \emph{Numerical approximation of hyperbolic systems of conservation
  laws}, volume 118 of \emph{Applied Mathematical Sciences}.
\newblock Springer-Verlag, New York, 1996.
\newblock ISBN 0-387-94529-6.
\newblock \doi{10.1007/978-1-4612-0713-9}.
\newblock URL \url{https://doi.org/10.1007/978-1-4612-0713-9}.

\bibitem[Goldberg(1991)]{Goldberg1991}
D.~Goldberg.
\newblock What every computer scientist should know about floating-point
  arithmetic.
\newblock \emph{ACM Comp. Surv.}, 23\penalty0 (1):\penalty0 5--48, 1991.

\bibitem[Hamilton and Webb(2013)]{Hamilton_etal2013}
B.~Hamilton and C.~Webb.
\newblock Room acoustics modelling using gpu-accelerated finite difference and
  finite volume methods on a face-centered cubic grid.
\newblock 01 2013.

\bibitem[Harten(1983)]{HartenAmi1983}
A.~Harten.
\newblock High resolution schemes for hyperbolic conservation laws.
\newblock \emph{J. Comput. Phys.}, 49\penalty0 (3):\penalty0 357--393, 1983.
\newblock ISSN 0021-9991,1090-2716.
\newblock \doi{10.1016/0021-9991(83)90136-5}.
\newblock URL \url{https://doi.org/10.1016/0021-9991(83)90136-5}.

\bibitem[Inglard et~al.(2020)Inglard, Lagouti\`ere, and Rugh]{Inglard2020}
M.~Inglard, F.~Lagouti\`ere, and H.~H. Rugh.
\newblock Ghost solutions with centered schemes for one-dimensional transport
  equations with {Neumann} boundary conditions.
\newblock \emph{Annales de la Facult\'e des sciences de Toulouse :
  Math\'ematiques}, Ser. 6, 29\penalty0 (4):\penalty0 927--950, 2020.
\newblock \doi{10.5802/afst.1650}.

\bibitem[Isaacson and Keller(1994)]{IsaacsonKeller1994}
E.~Isaacson and H.~B. Keller.
\newblock \emph{Analysis of Numerical Methods}.
\newblock Dover, 1994.

\bibitem[Lai and Peskin(2000)]{Lai_etal2000}
M.-C. Lai and C.~S. Peskin.
\newblock An immersed boundary method with formal second-order accuracy and
  reduced numerical viscosity.
\newblock \emph{J. Comput. Phys.}, 160\penalty0 (2):\penalty0 705--719, 2000.
\newblock ISSN 0021-9991,1090-2716.
\newblock \doi{10.1006/jcph.2000.6483}.
\newblock URL \url{https://doi.org/10.1006/jcph.2000.6483}.

\bibitem[Lanczos(1988)]{lanczos1988applied}
C.~Lanczos.
\newblock \emph{Applied Analysis}.
\newblock Dover Books on Mathematics. Dover Publications, 1988.
\newblock ISBN 9780486656564.

\bibitem[LeVeque(2002)]{leveque2002}
R.~J. LeVeque.
\newblock \emph{Finite Volume Methods for Hyperbolic Problems}.
\newblock Cambridge Texts in Applied Mathematics. Cambridge University Press,
  2002.

\bibitem[LeVeque and Li(1994)]{LeVeque_etal1994}
R.~J. LeVeque and Z.~L. Li.
\newblock The immersed interface method for elliptic equations with
  discontinuous coefficients and singular sources.
\newblock \emph{SIAM J. Numer. Anal.}, 31\penalty0 (4):\penalty0 1019--1044,
  1994.
\newblock ISSN 0036-1429.
\newblock \doi{10.1137/0731054}.
\newblock URL \url{https://doi.org/10.1137/0731054}.

\bibitem[Peskin(2002)]{PeskinCharles2002}
C.~S. Peskin.
\newblock The immersed boundary method.
\newblock \emph{Acta Numer.}, 11:\penalty0 479--517, 2002.
\newblock ISSN 0962-4929,1474-0508.
\newblock \doi{10.1017/S0962492902000077}.
\newblock URL \url{https://doi.org/10.1017/S0962492902000077}.

\bibitem[Sod(1978)]{SodGary1978}
G.~A. Sod.
\newblock A survey of several finite difference methods for systems of
  nonlinear hyperbolic conservation laws.
\newblock \emph{J. Comput. Phys.}, 27\penalty0 (1):\penalty0 1--31, 1978.
\newblock ISSN 0021-9991,1090-2716.
\newblock \doi{10.1016/0021-9991(78)90023-2}.
\newblock URL \url{https://doi.org/10.1016/0021-9991(78)90023-2}.

\bibitem[Spalding(1972)]{Spalding1972}
D.~B. Spalding.
\newblock A novel finite-difference formulation for differential expression
  involving both first and second derivatives.
\newblock \emph{Int. J. Numer. Methods Eng.}, 4:\penalty0 551--559, 1972.

\bibitem[Szustak et~al.(2016)Szustak, Halbiniak, Kulawik, Wrobel, and
  Gepner]{Szustak_etal2016}
L.~Szustak, K.~Halbiniak, A.~Kulawik, J.~Wrobel, and P.~Gepner.
\newblock Toward parallel modeling of solidification based on the generalized
  finite difference method using {I}ntel {X}eon {P}hi.
\newblock In \emph{Parallel processing and applied mathematics. {P}art {I}},
  volume 9573 of \emph{Lecture Notes in Comput. Sci.}, pages 411--422.
  Springer, [Cham], 2016.
\newblock ISBN 978-3-319-32149-3; 978-3-319-32148-6.
\newblock \doi{10.1007/978-3-319-32149-3\_39}.
\newblock URL \url{https://doi.org/10.1007/978-3-319-32149-3_39}.

\bibitem[Taylor(1715)]{taylor1715}
B.~Taylor.
\newblock \emph{{Methodus incrementorum directa inversa. Auctore Brook Taylor,
  LL. D. and Regiae Societatis Secretario}}.
\newblock Typis Pearsonianis Prostant apud Gul. Innys ad Insignia Principis in
  Coemeterio Paulino MDCCXV, Londini, 1715.
\newblock URL \url{https://books.google.com/books?id=iXN1xgEACAAJ}.

\bibitem[Ye et~al.(2022)Ye, Zhang, Wan, Yan, and Sun]{Ye2022}
C.-C. Ye, P.-J.-Y. Zhang, Z.-H. Wan, R.~Yan, and D.-J. Sun.
\newblock Accelerating cfd simulation with high order finite difference method
  on curvilinear coordinates for modern gpu clusters.
\newblock \emph{Advances in Aerodynamics}, 4\penalty0 (1), Feb. 2022.
\newblock ISSN 2524-6992.
\newblock \doi{10.1186/s42774-021-00098-3}.
\newblock URL \url{http://dx.doi.org/10.1186/s42774-021-00098-3}.

\end{thebibliography}

\newpage
\appendix
\section{Magnitude of the finite difference coefficients} \label{CoeffsMagnitudeAnnex}
\begin{figure} \label{FigCoeffsMag}
    \centering
    \begin{minipage}{0.48\textwidth}
        \centering
        \includegraphics[width=\linewidth]{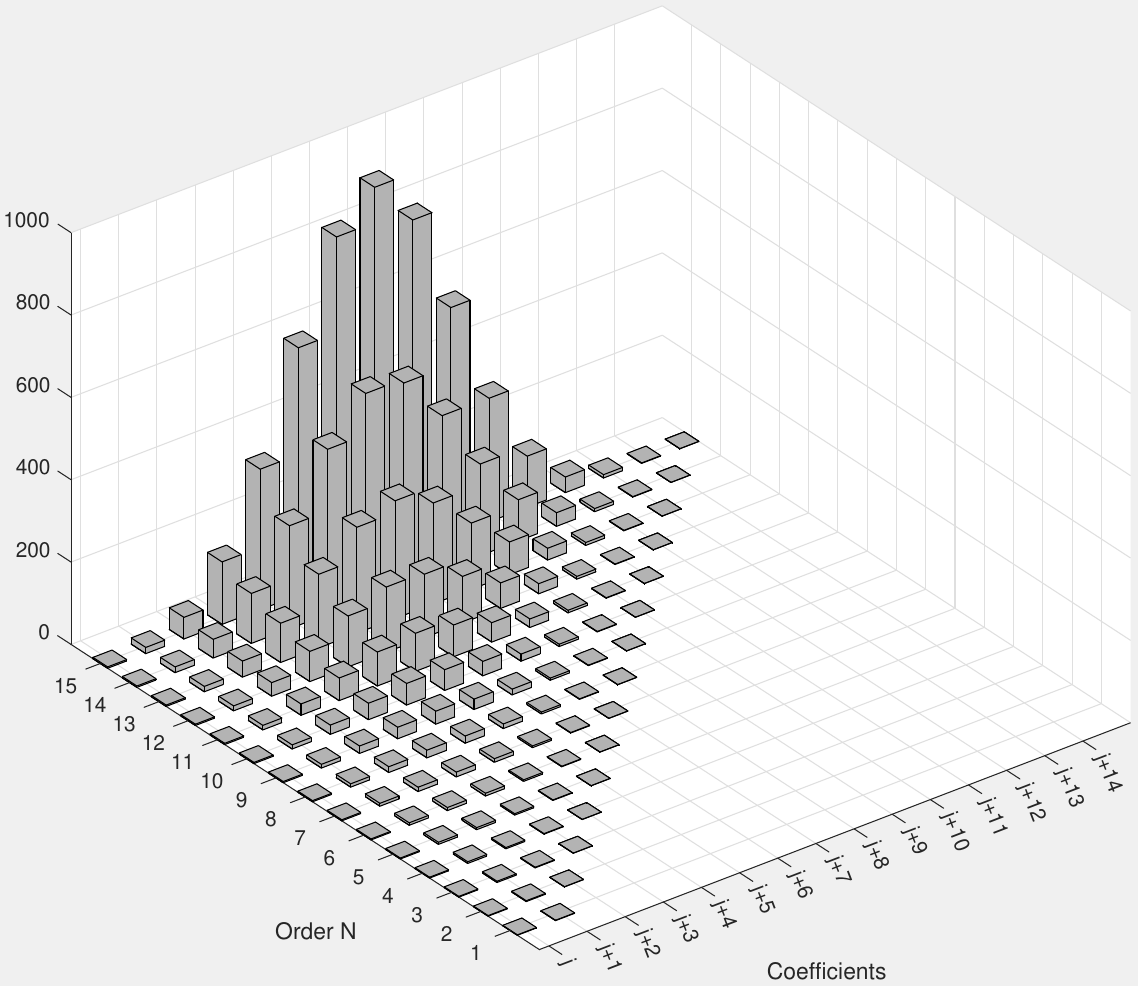}
        \subcaption{Magnitude of the coefficients $D_N^j$ for the forward schemes (c.f. equation \eqref{EQ_forward_coeffs})}
        \label{F_coeff_mag}
    \end{minipage}%
    \hfill
    \begin{minipage}{0.48\textwidth}
        \centering
        \includegraphics[width=\linewidth]{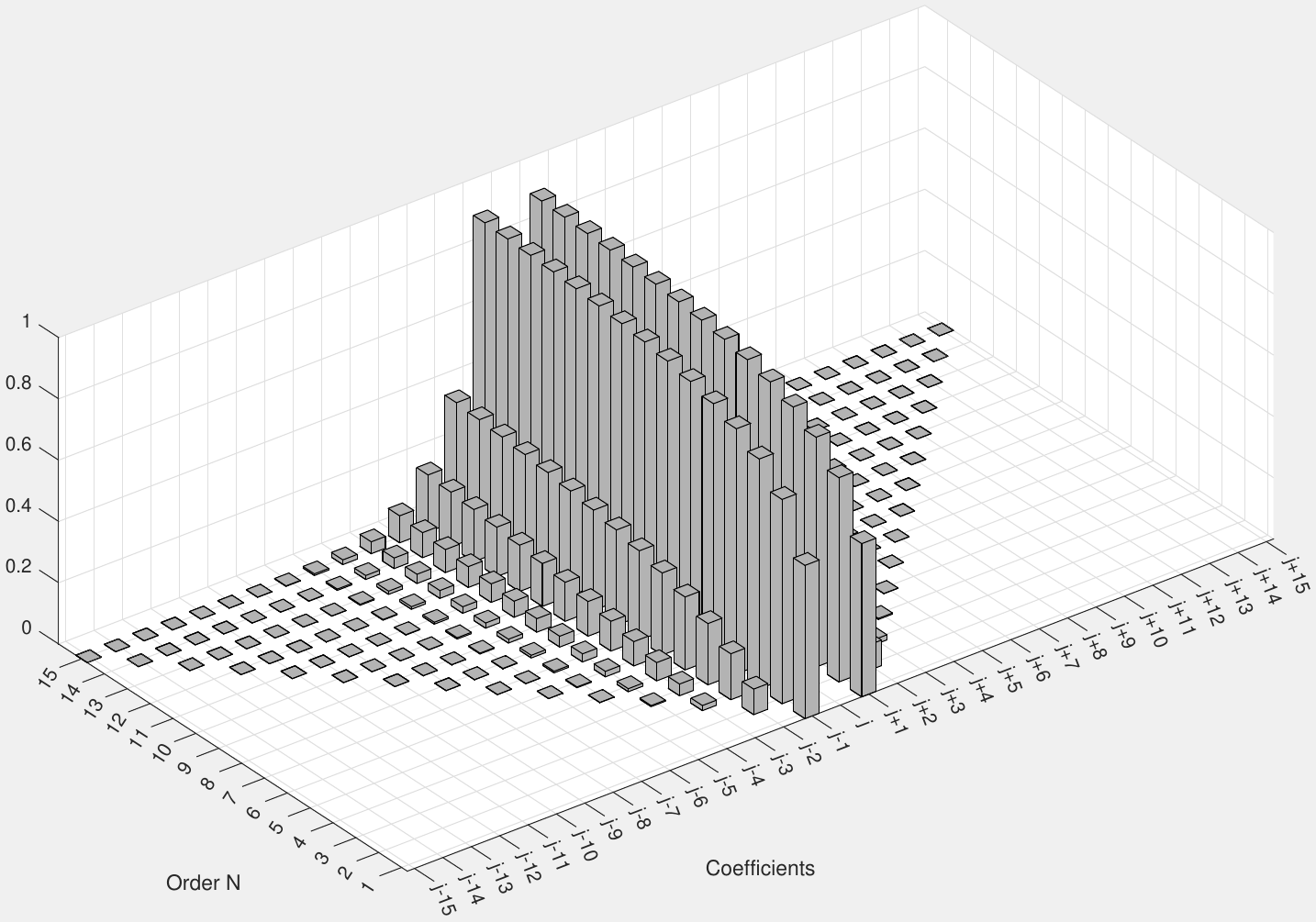}
        \subcaption{Magnitude of the coefficients $C_N^j$ for the centred schemes (c.f. equation \eqref{EQ_centred_coeffs}).}
        \label{C_coeff_mag}
    \end{minipage}
    
    \vskip\baselineskip

    \begin{minipage}{0.48\textwidth}
        \centering
        \includegraphics[width=\linewidth]{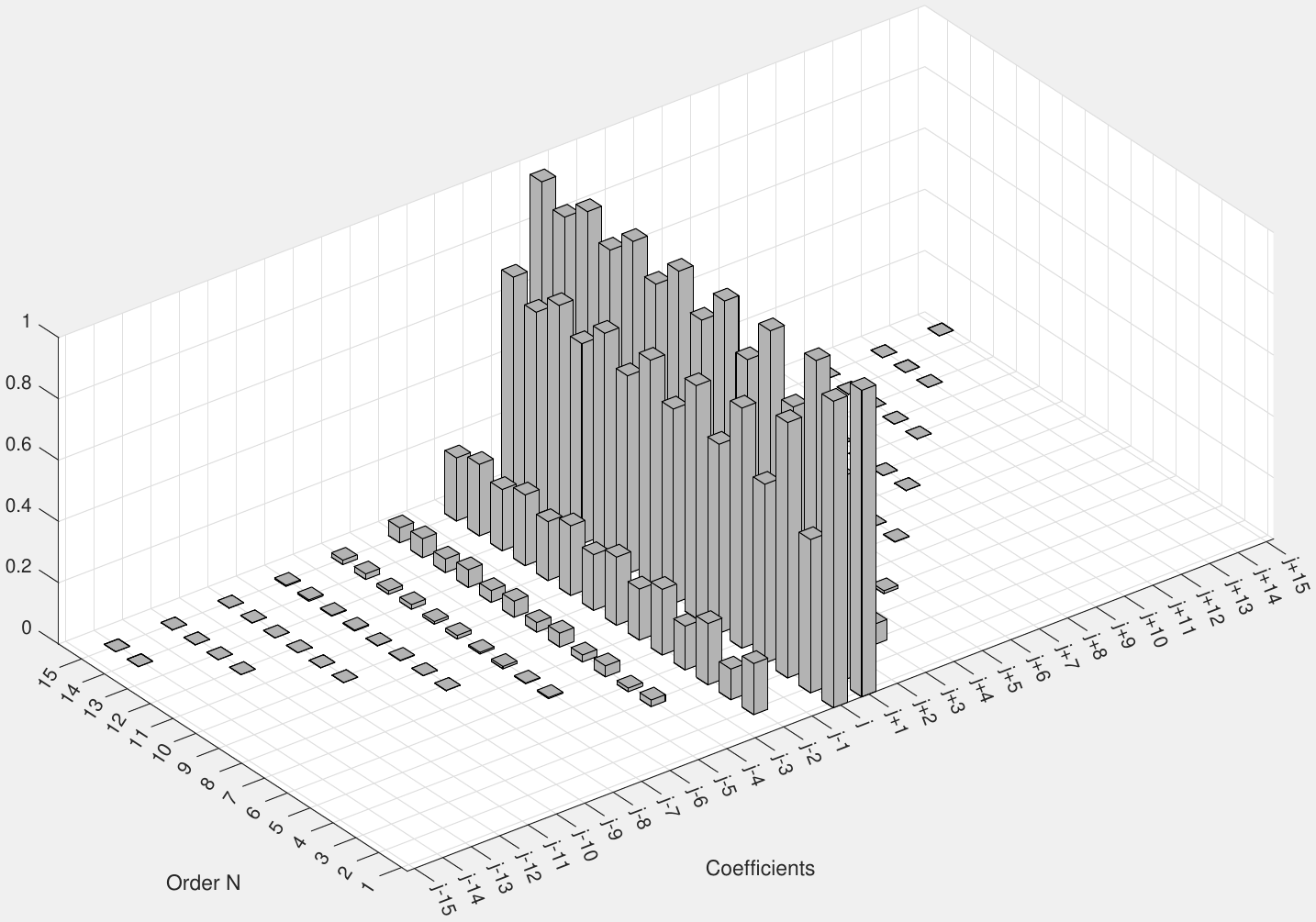}
        \subcaption{Magnitude of the coefficients $Z_N^j$ for the forward-first zigzag schemes (c.f. equations \eqref{defZN1}--\eqref{defZNj}).}
        \label{ZZF_coeff_mag}
    \end{minipage}%
    \hfill
    \begin{minipage}{0.48\textwidth}
        \centering
        \includegraphics[width=\linewidth]{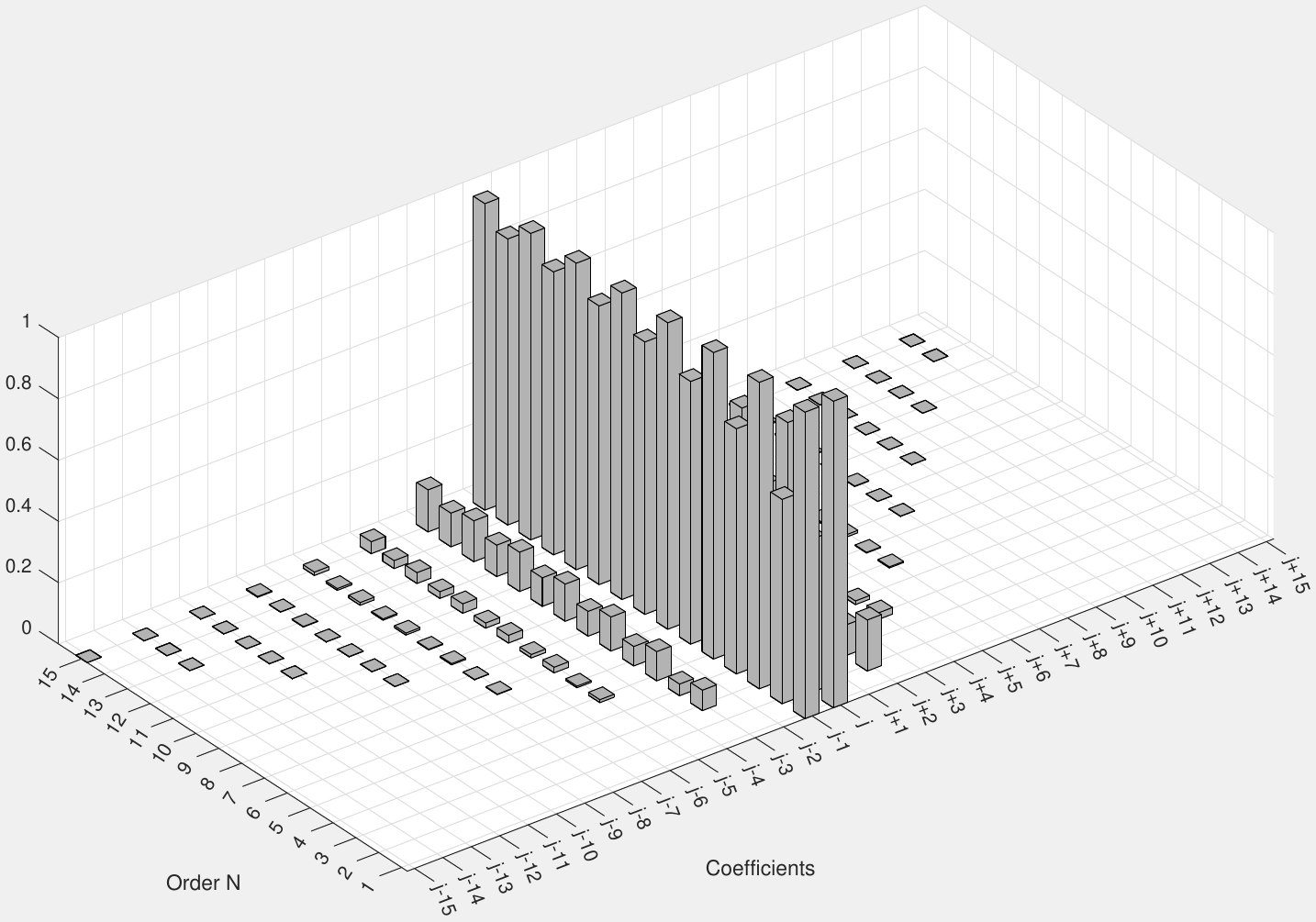}
        \subcaption{Magnitude of the coefficients $Z_N^j$ for the backward-first zigzag scheme (c.f. equations \eqref{defZN1}--\eqref{defZNj}).}
        \label{ZZB_coeff_mag}
    \end{minipage}
\caption{Magnitude of the coefficients used by the centred, one-sided and zigzag schemes.}
\end{figure}

To make it easier to compare the values of the various coefficients involved, we provide four figures (Fig. \ref{F_coeff_mag} to \ref{ZZB_coeff_mag}) illustrating the magnitude of the finite difference coefficients for the forward, centred and zigzag schemes. \\

\newpage
\section{Numerical diffusion of the zigzag schemes}
\label{AnnexNumDiff}

The centred, upwind and zigzag finite differences schemes all 
introduce some
numerical diffusion. If we compare the rate at which these schemes 
dissipate, we notice that the energy loss from the zigzag schemes 
is lower than the one-sided schemes and greater than the centred 
schemes. This middle-ground is not surprising since the zigzag 
schemes can be seen as a mix of one-sided and centred schemes. To 
illustrate this observation with a concrete example, we consider 
the solution of the transport equation \eqref{T} using a 
Runge--Kutta scheme of order $3$ and spatial schemes of order $2$. 
We solve this equation for $(t,x) \in [0, 15] \times [-20,20]$ 
using a periodised error function as the initial condition, with 
the parameters $\Delta t = 0.05$ and $\Delta x = 0.01$. For the 
centred, zigzag and upwind schemes we obtain an energy loss of 
respectively $2.18 \cdot 10^{-9}$, $3.71 \cdot 10^{-8}$ and $1.07 
\cdot 10^{-7}$.

\newpage
\section{Stability regions for the forward Euler method} 
\label{EulerAnnex}

\begin{figure}[H]
\minipage{0.18\textwidth}
  \includegraphics[width=\linewidth]{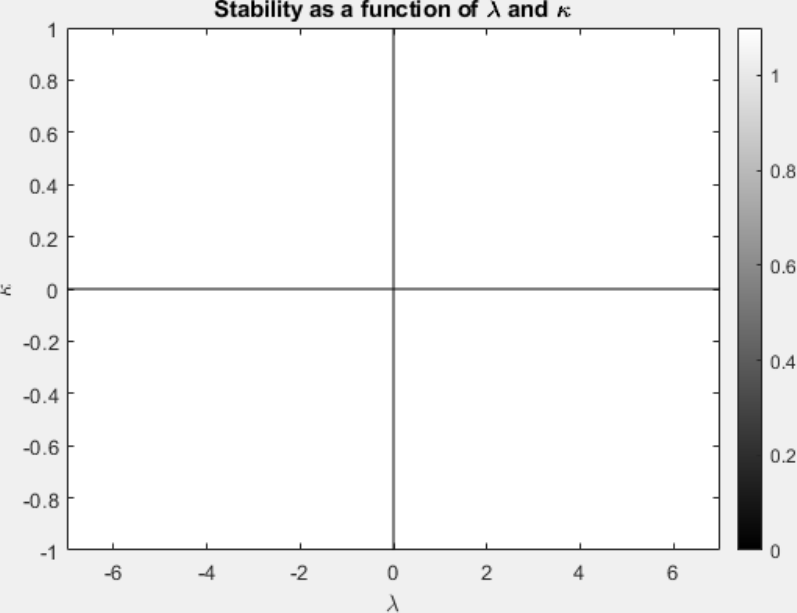}
\endminipage\hfill
\minipage{0.19\textwidth}
  \includegraphics[width=\linewidth]{blank.pdf}
\endminipage\hfill
\minipage{0.19\textwidth}
  \includegraphics[width=\linewidth]{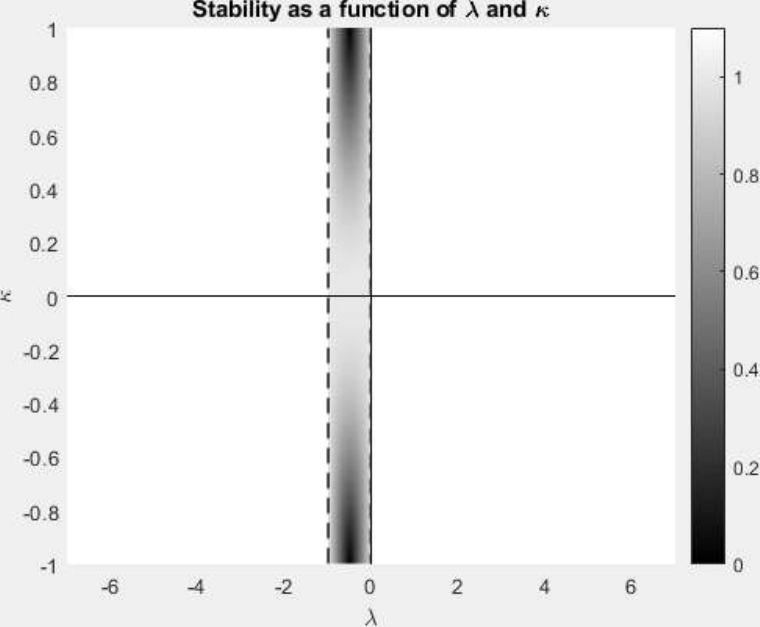}
\endminipage\hfill
\minipage{0.19\textwidth}
  \includegraphics[width=\linewidth]{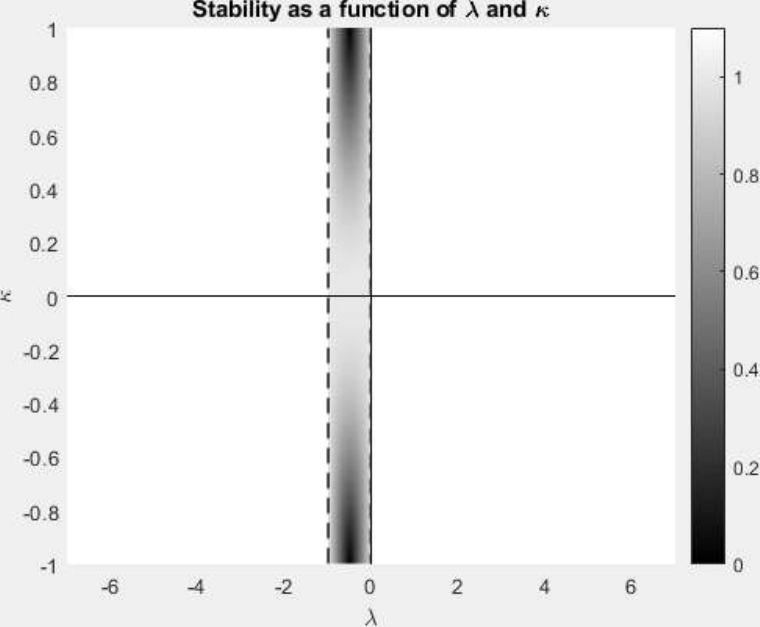}
\endminipage\hfill
\minipage{0.19\textwidth}%
  \includegraphics[width=\linewidth]{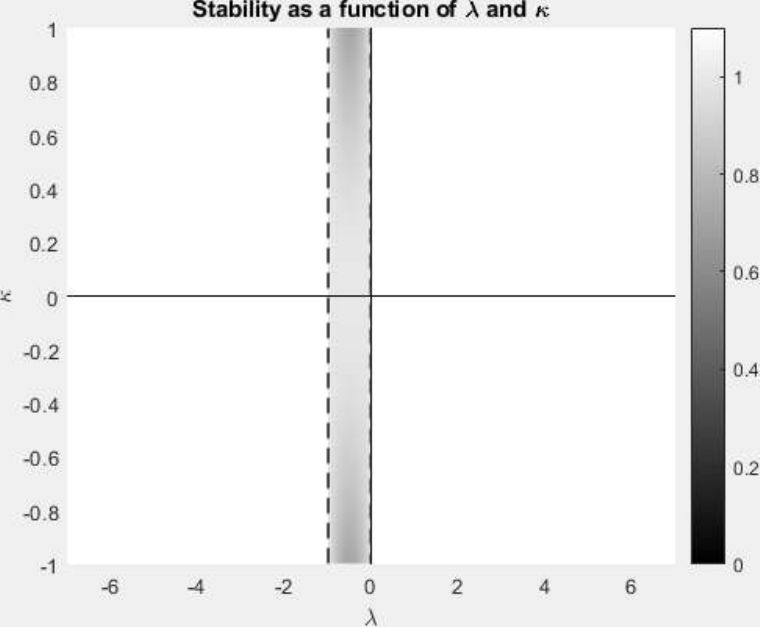}
\endminipage
\caption{FE; left to right : stability regions for \textbf{Order 1} c, cs, f, zff and zffs. \label{fe1} }
\end{figure}

\begin{figure}[H]
\minipage{0.19\textwidth}
  \includegraphics[width=\linewidth]{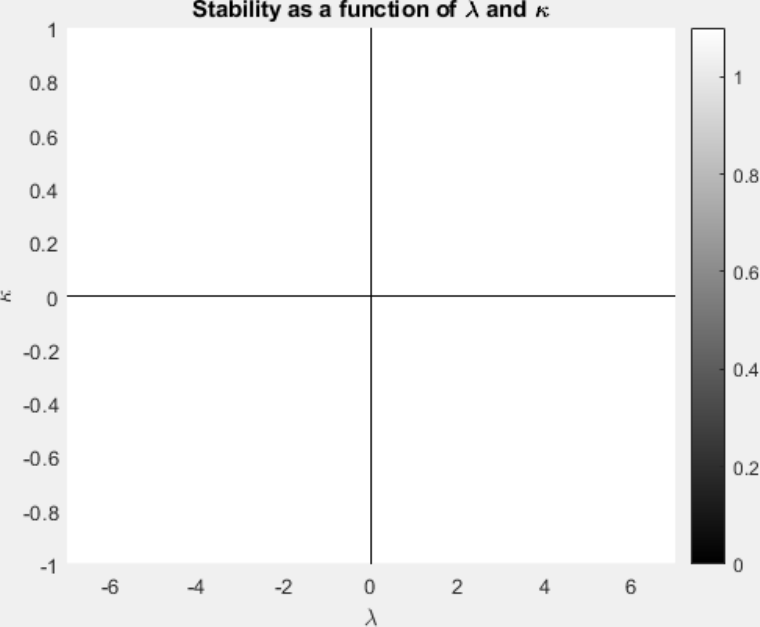}
\endminipage\hfill
\minipage{0.19\textwidth}
  \includegraphics[width=\linewidth]{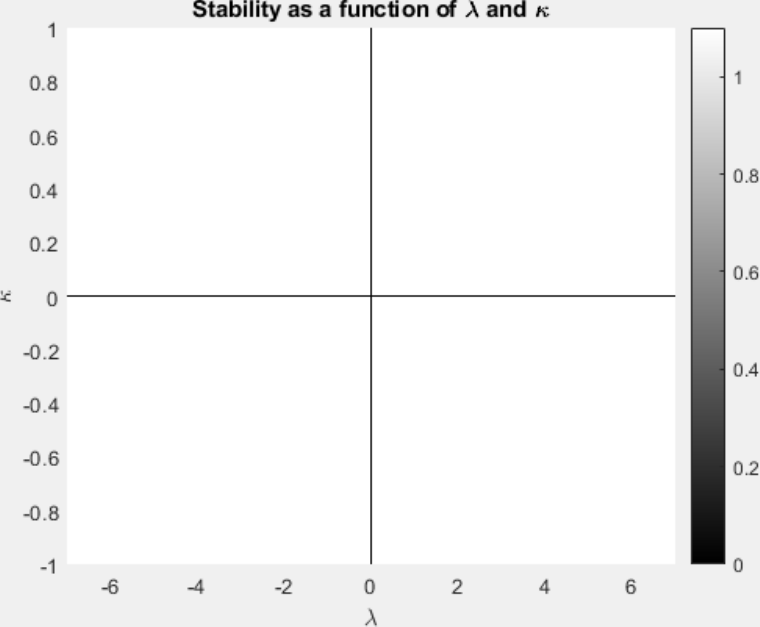}
\endminipage\hfill
\minipage{0.19\textwidth}
  \includegraphics[width=\linewidth]{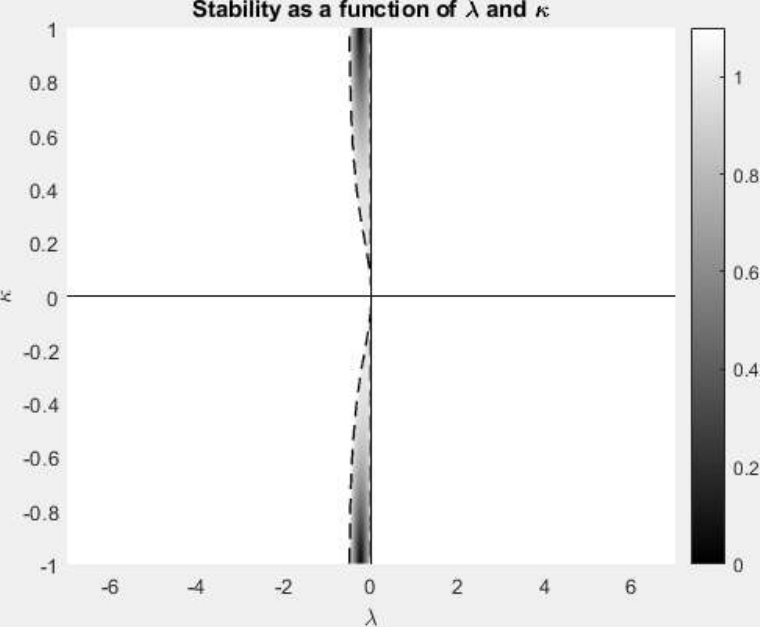}
\endminipage\hfill
\minipage{0.19\textwidth}
  \includegraphics[width=\linewidth]{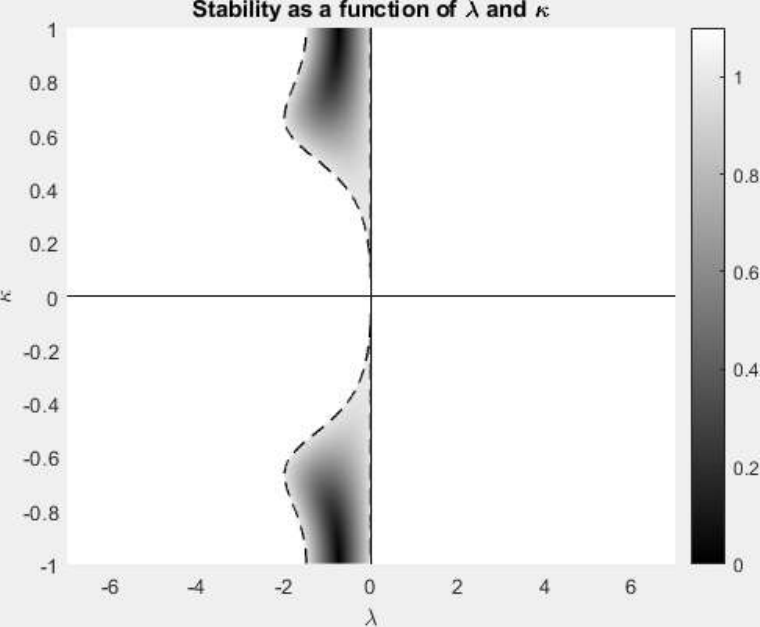}
\endminipage\hfill
\minipage{0.19\textwidth}%
  \includegraphics[width=\linewidth]{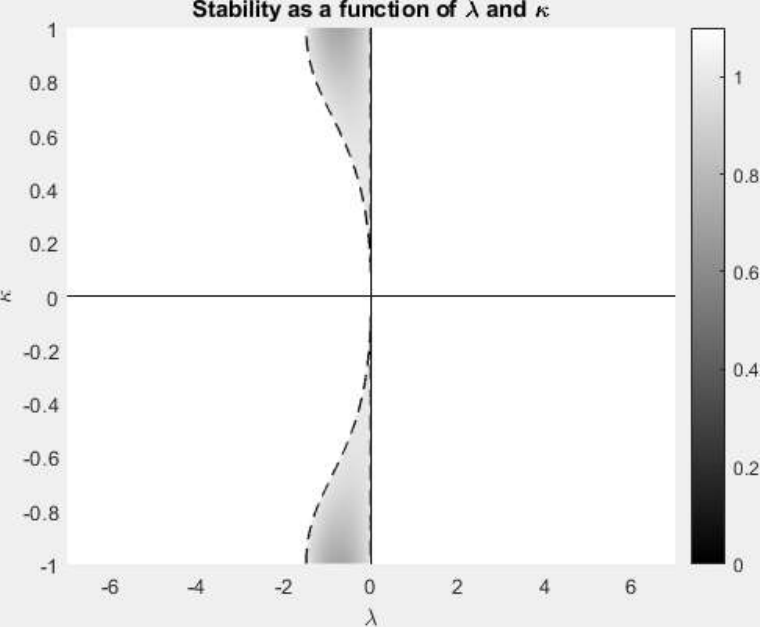}
\endminipage
\caption{FE; left to right : stability regions for \textbf{Order 2} c, cs, f, zff and zffs. \label{fe2} }
\end{figure}

\begin{figure}[H]
\minipage{0.19\textwidth}
  \includegraphics[width=\linewidth]{blank.pdf}
\endminipage\hfill
\minipage{0.19\textwidth}
  \includegraphics[width=\linewidth]{blank.pdf}
\endminipage\hfill
\minipage{0.19\textwidth}
  \includegraphics[width=\linewidth]{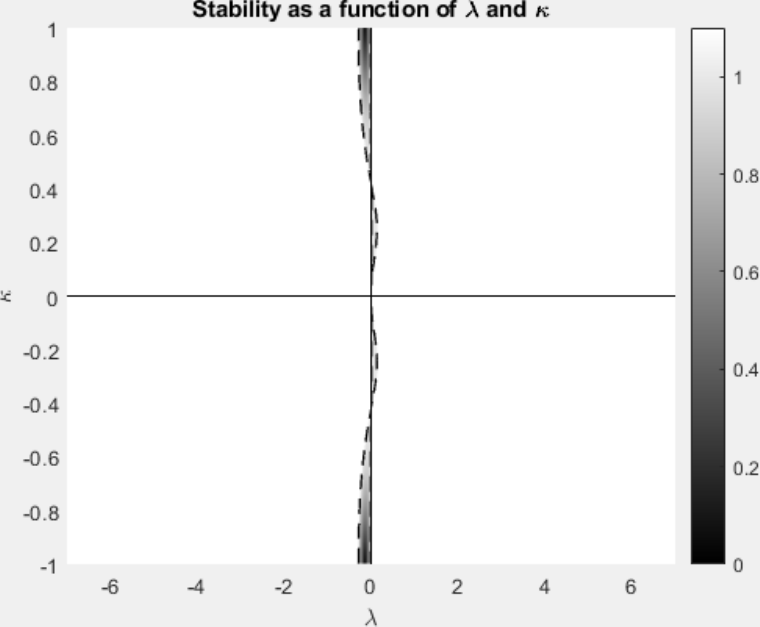}
\endminipage\hfill
\minipage{0.19\textwidth}
  \includegraphics[width=\linewidth]{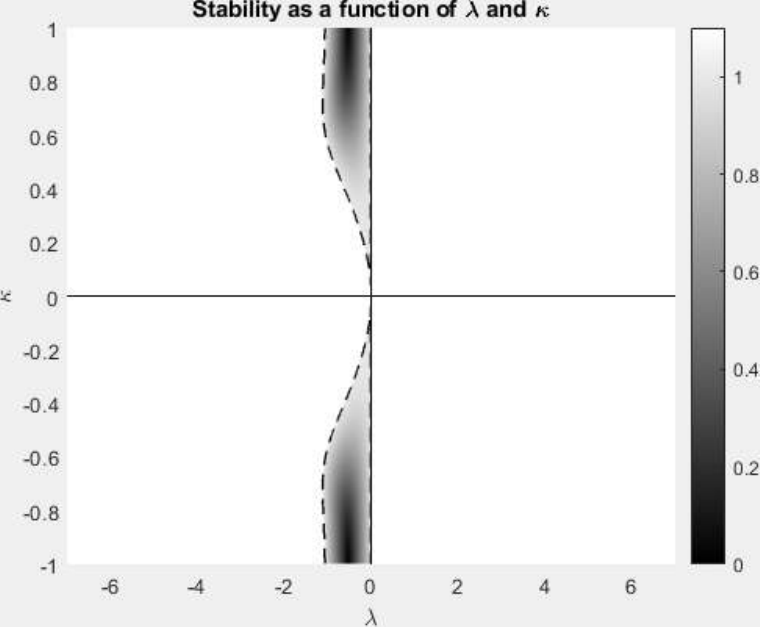}
\endminipage\hfill
\minipage{0.19\textwidth}%
  \includegraphics[width=\linewidth]{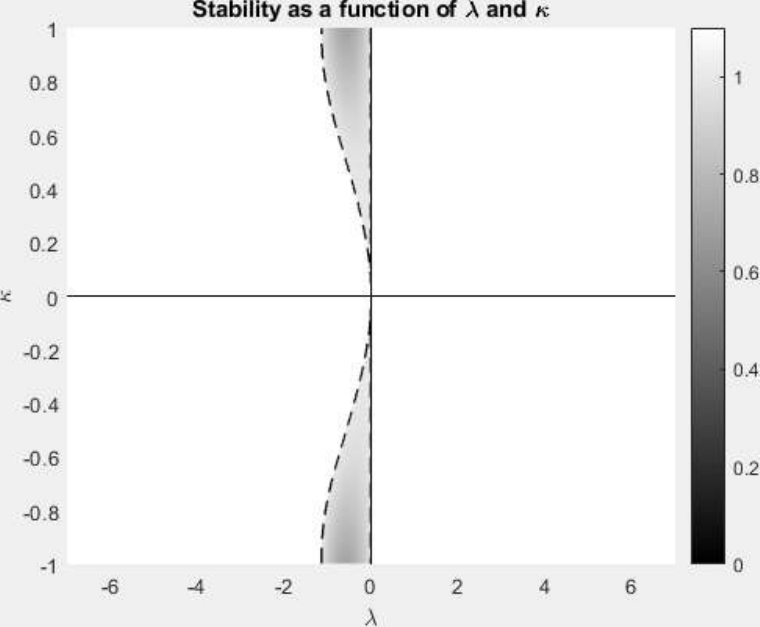}
\endminipage
\caption{FE; left to right : stability regions for \textbf{Order 3} c, cs, f, zff and zffs. \label{fe3} }
\end{figure}

\begin{figure}[H]
\minipage{0.19\textwidth}
  \includegraphics[width=\linewidth]{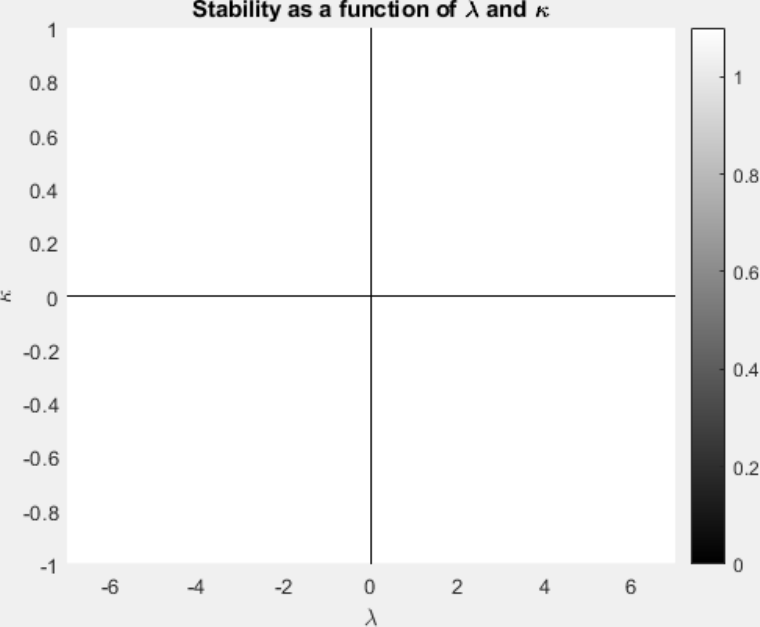}
\endminipage\hfill
\minipage{0.19\textwidth}
  \includegraphics[width=\linewidth]{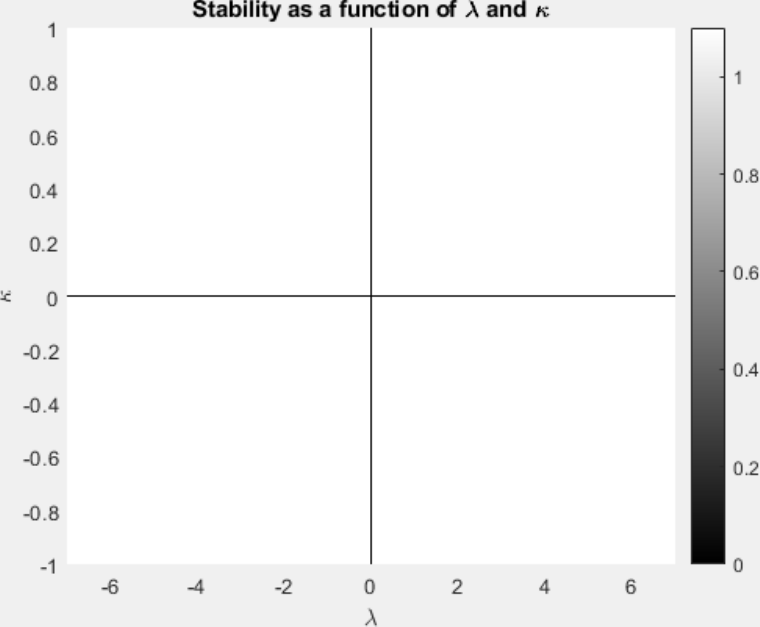}
\endminipage\hfill
\minipage{0.19\textwidth}
  \includegraphics[width=\linewidth]{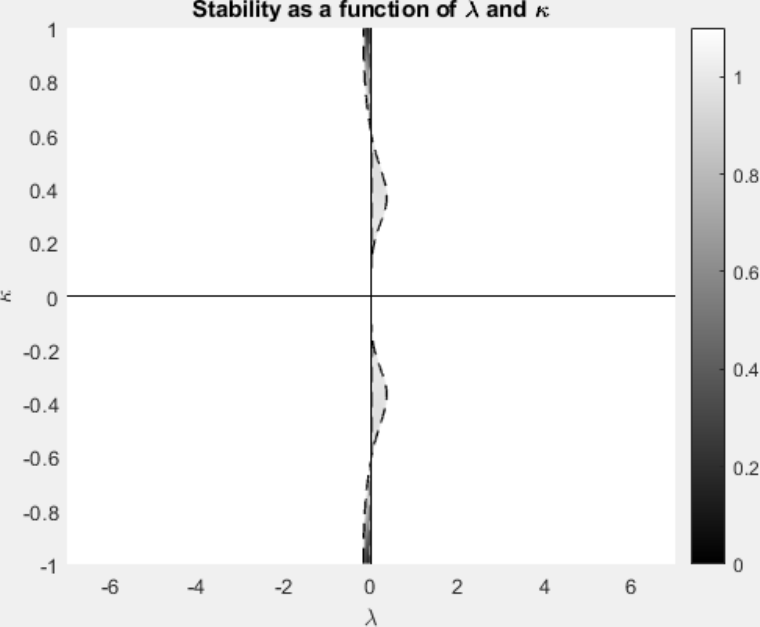}
\endminipage\hfill
\minipage{0.19\textwidth}
  \includegraphics[width=\linewidth]{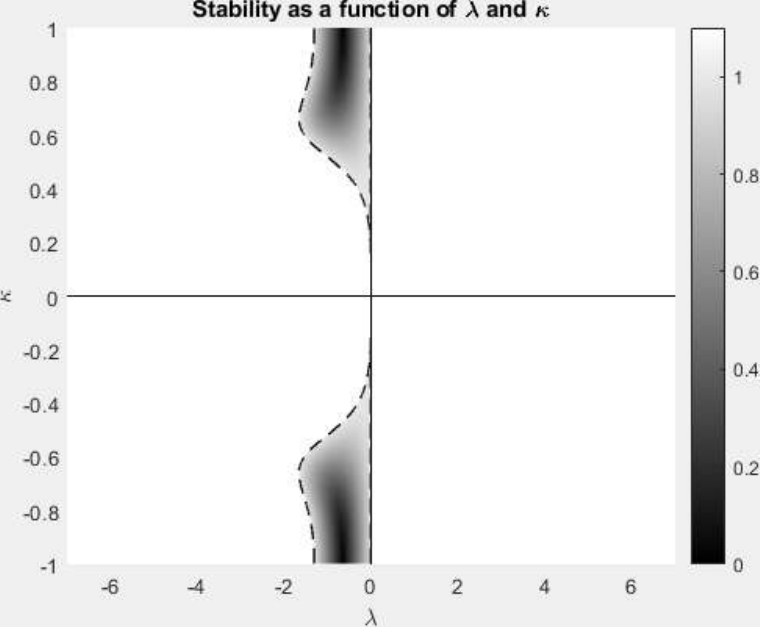}
\endminipage\hfill
\minipage{0.19\textwidth}%
  \includegraphics[width=\linewidth]{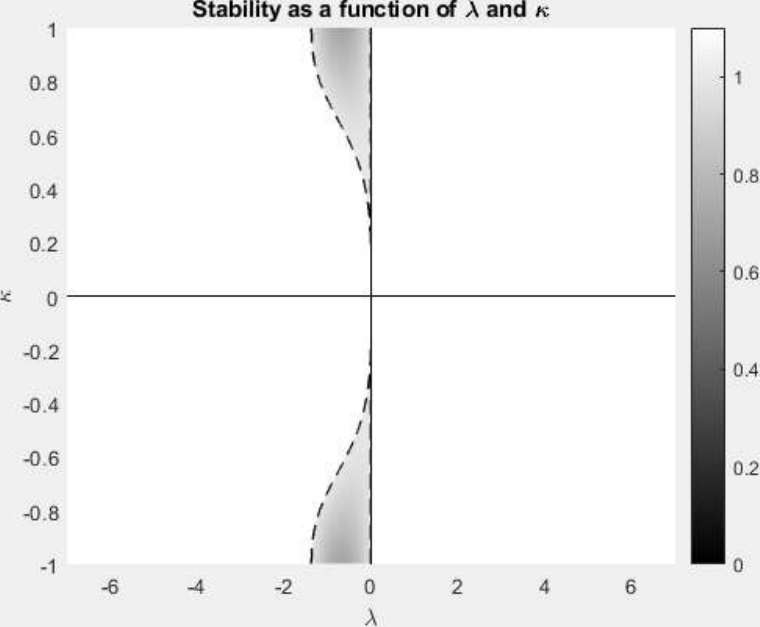}
\endminipage
\caption{FE; left to right : stability regions for \textbf{Order 4} c, cs, f, zff and zffs. \label{fe4} }
\end{figure}

\begin{figure}[H]
\minipage{0.19\textwidth}
  \includegraphics[width=\linewidth]{blank.pdf}
\endminipage\hfill
\minipage{0.19\textwidth}
  \includegraphics[width=\linewidth]{blank.pdf}
\endminipage\hfill
\minipage{0.19\textwidth}
  \includegraphics[width=\linewidth]{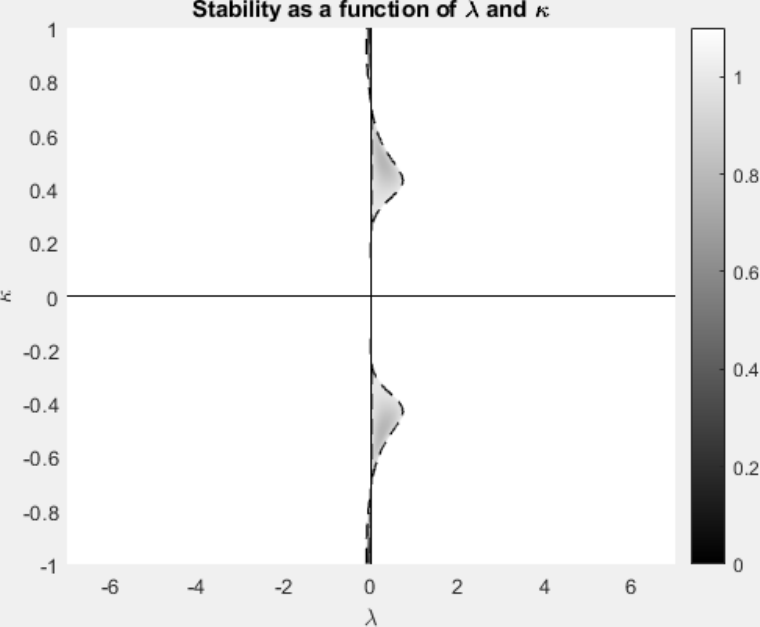}
\endminipage\hfill
\minipage{0.19\textwidth}
  \includegraphics[width=\linewidth]{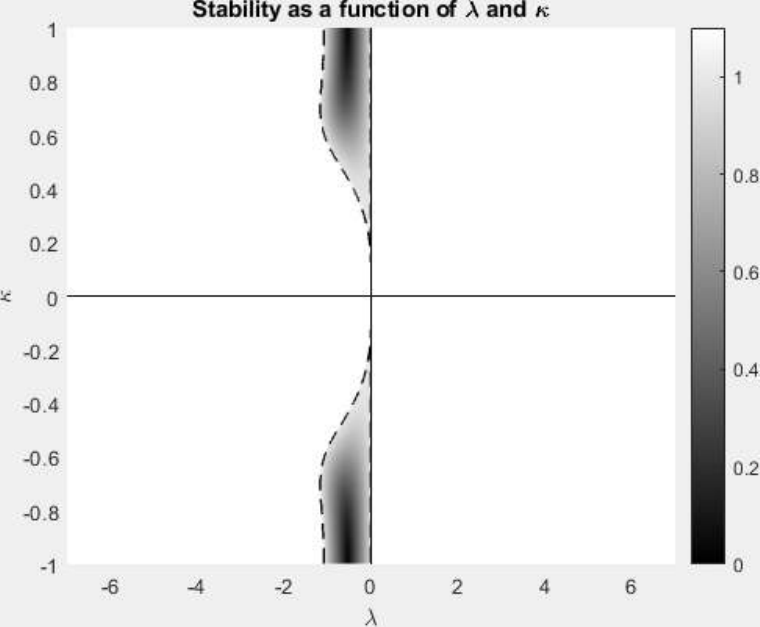}
\endminipage\hfill
\minipage{0.19\textwidth}%
  \includegraphics[width=\linewidth]{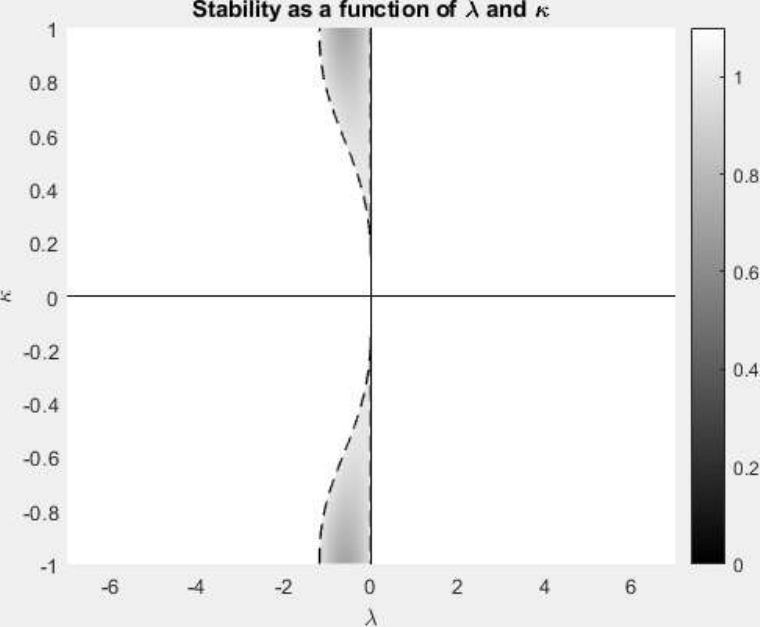}
\endminipage
\caption{FE; left to right : stability regions for \textbf{Order 5} c, cs, f, zff and zffs. \label{fe5} }
\end{figure}

We provide Fig. \ref{fe1} to \ref{fe5} the stability regions in the 
$\lambda-\kappa$ plane ($\lambda$ and $\kappa$ respectively being 
the stability-number and the scaled wavenumber) for the forward Euler 
(FE) method in time,  for spatial orders $N=1$ to $5$, respectively 
(left to right) for centred (c), centred staggered (cs), forward (f), 
zigzag forward-first (zff) and zigzag 
forward-first staggered (zffs) schemes. The greyed out areas 
represent the couples $(\lambda,\kappa)$ for which the scheme is stable 
(i.e. the amplification factor $|G| \leq 1$), while the dotted contours 
correspond to the critical case $|G| = 1$. A given scheme is conditionally
stable if and only if there exist a $\lambda_c \in \mathds{R}$ such that
the line of equation $\lambda = \lambda_c$ is included in the grey area.

\newpage
\section{Stability regions for RK2} \label{RK2Annex}
\begin{figure}[H]
\minipage{0.18\textwidth}
  \includegraphics[width=\linewidth]{blank.pdf}
\endminipage\hfill
\minipage{0.19\textwidth}
  \includegraphics[width=\linewidth]{blank.pdf}
\endminipage\hfill
\minipage{0.19\textwidth}
  \includegraphics[width=\linewidth]{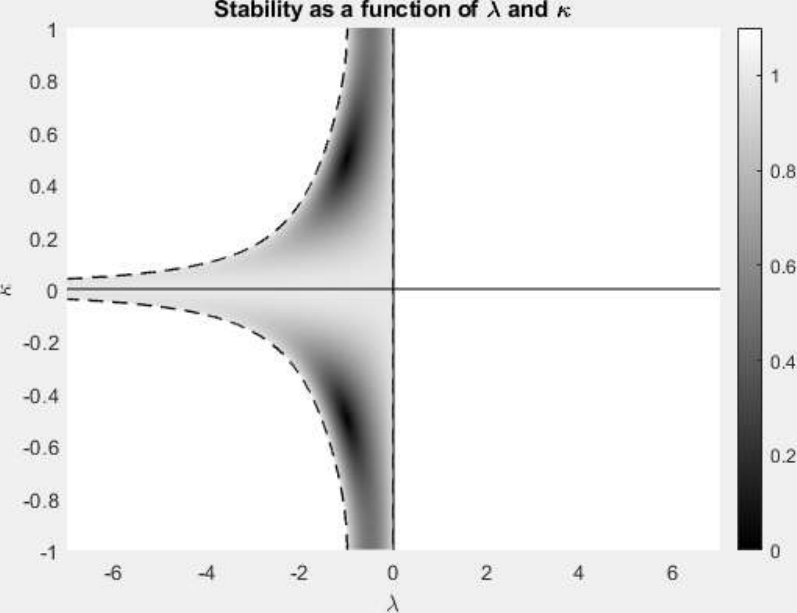}
\endminipage\hfill
\minipage{0.19\textwidth}
  \includegraphics[width=\linewidth]{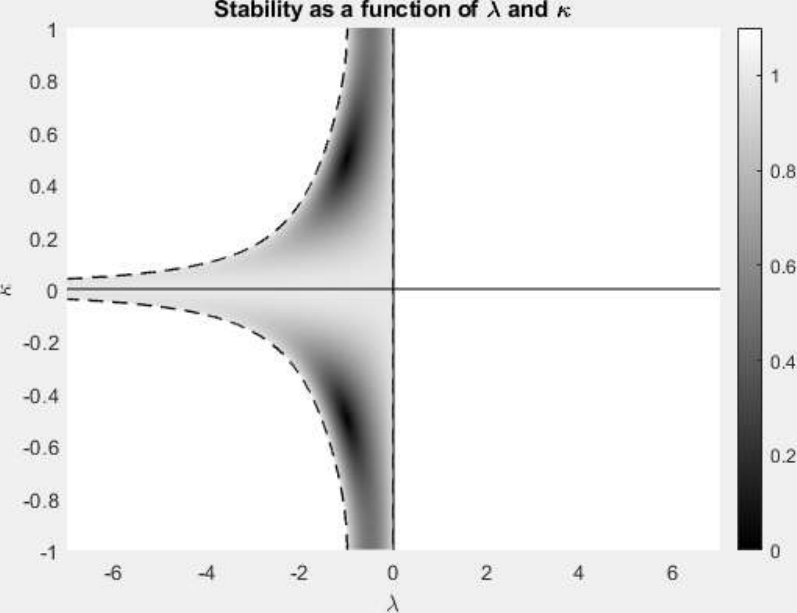}
\endminipage\hfill
\minipage{0.19\textwidth}%
  \includegraphics[width=\linewidth]{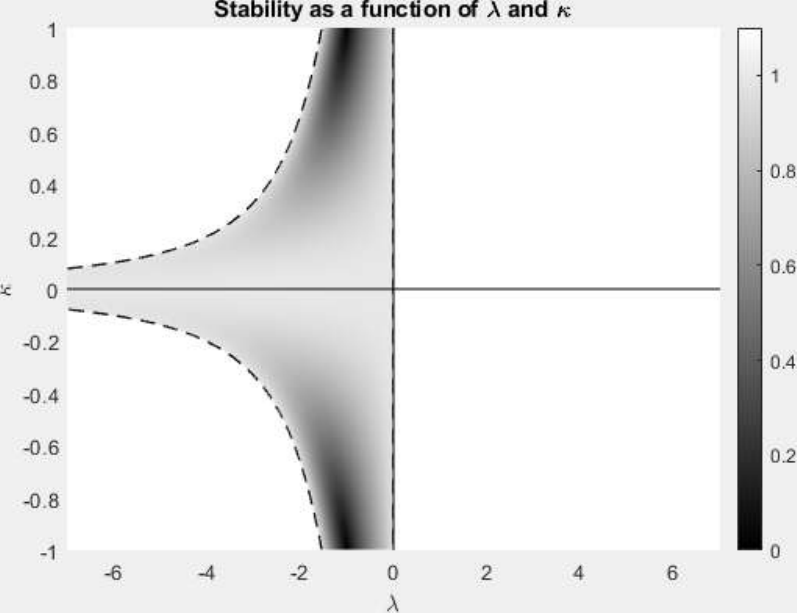}
\endminipage
\caption{RK2; left to right : stability regions for \textbf{Order 1} c, cs, f, zff and zffs. \label{rk21} }
\end{figure}

\begin{figure}[H]
\minipage{0.19\textwidth}
  \includegraphics[width=\linewidth]{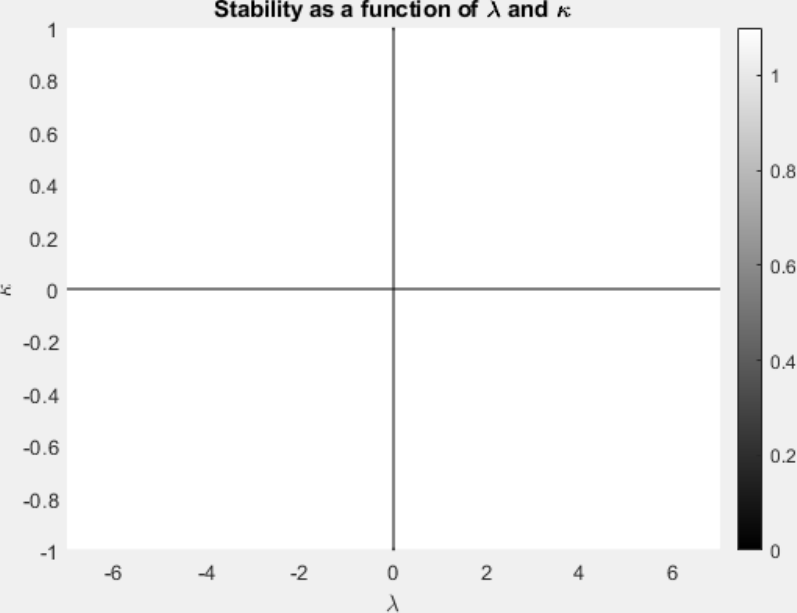}
\endminipage\hfill
\minipage{0.19\textwidth}
\includegraphics[width=\linewidth]{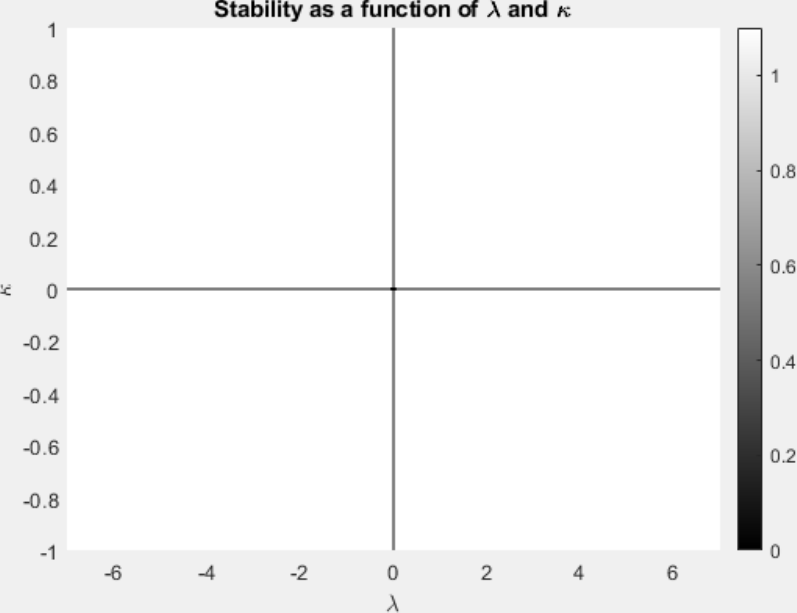}
\endminipage\hfill
\minipage{0.19\textwidth}
  \includegraphics[width=\linewidth]{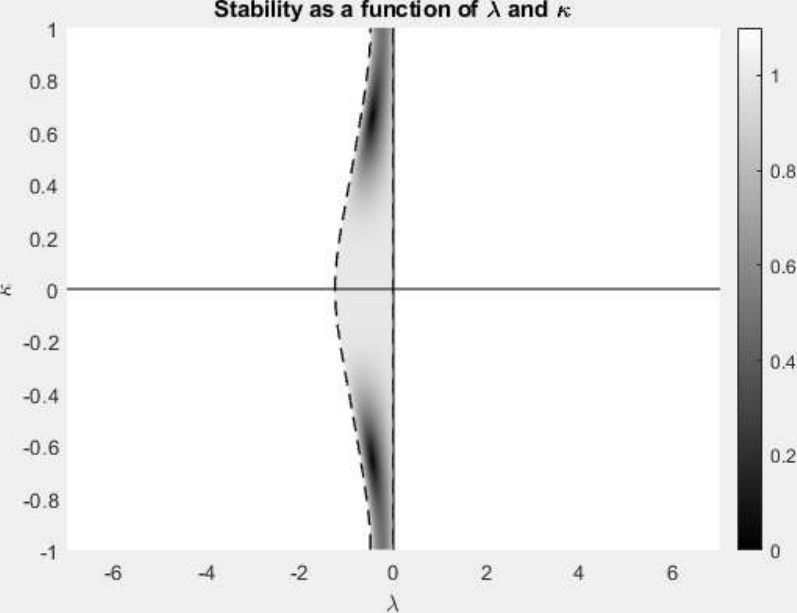}
\endminipage\hfill
\minipage{0.19\textwidth}
  \includegraphics[width=\linewidth]{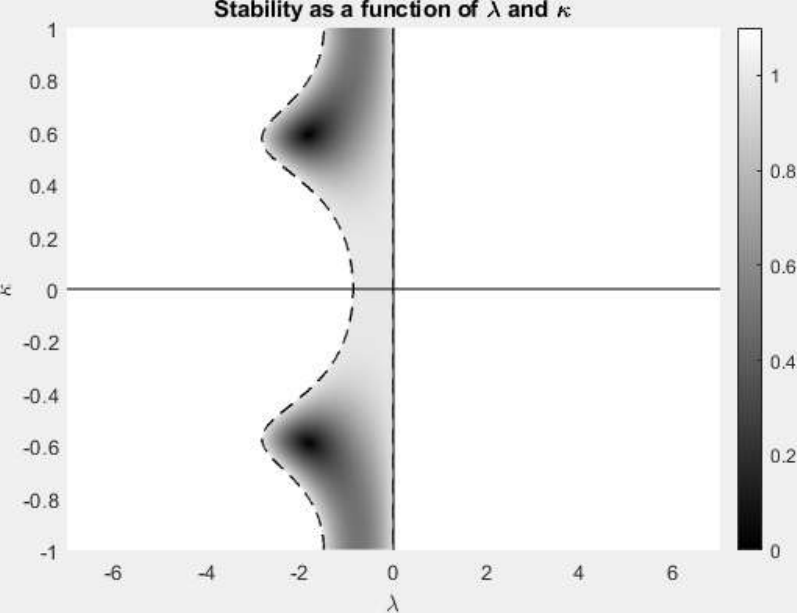}
\endminipage\hfill
\minipage{0.19\textwidth}%
  \includegraphics[width=\linewidth]{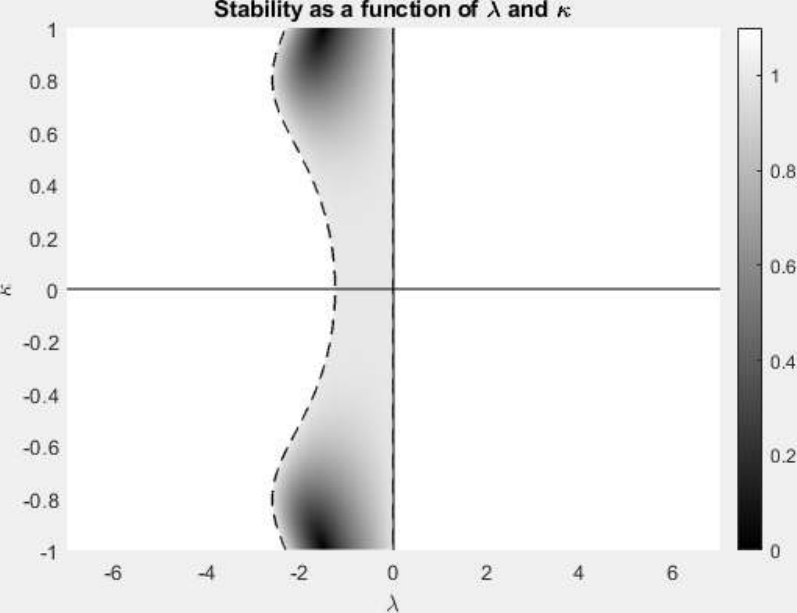}
\endminipage
\caption{RK2; left to right : stability regions for \textbf{Order 2} c, cs, f, zff and zffs. \label{rk22}  }
\end{figure}

\begin{figure}[H]
\minipage{0.19\textwidth}
  \includegraphics[width=\linewidth]{blank.pdf}
\endminipage\hfill
\minipage{0.19\textwidth}
  \includegraphics[width=\linewidth]{blank.pdf}
\endminipage\hfill
\minipage{0.19\textwidth}
  \includegraphics[width=\linewidth]{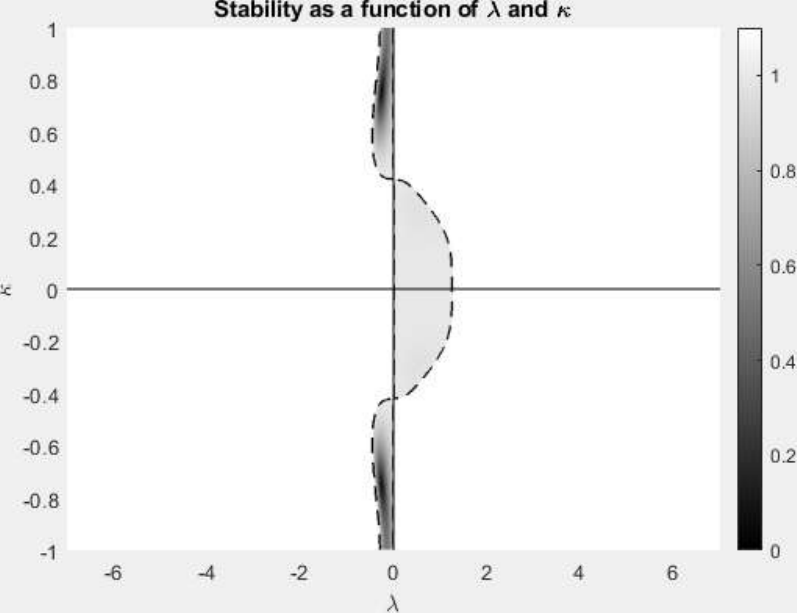}
\endminipage\hfill
\minipage{0.19\textwidth}
  \includegraphics[width=\linewidth]{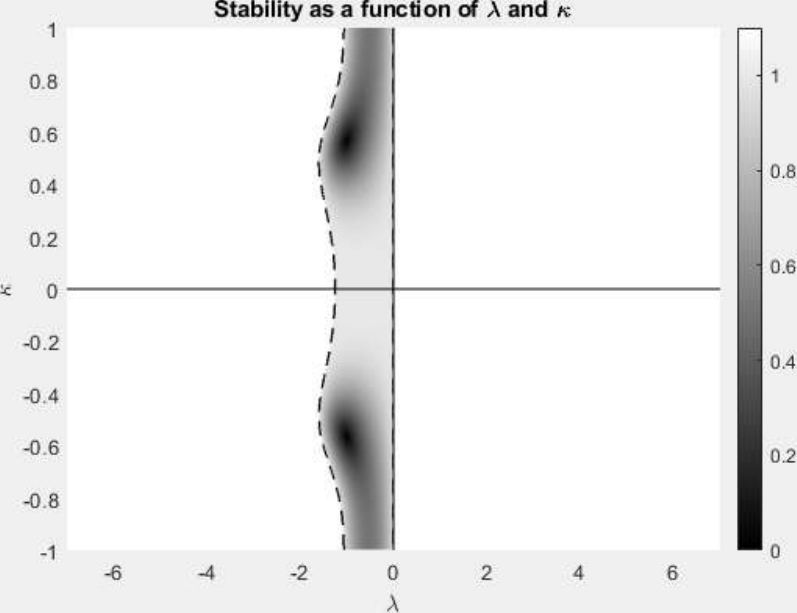}
\endminipage\hfill
\minipage{0.19\textwidth}%
  \includegraphics[width=\linewidth]{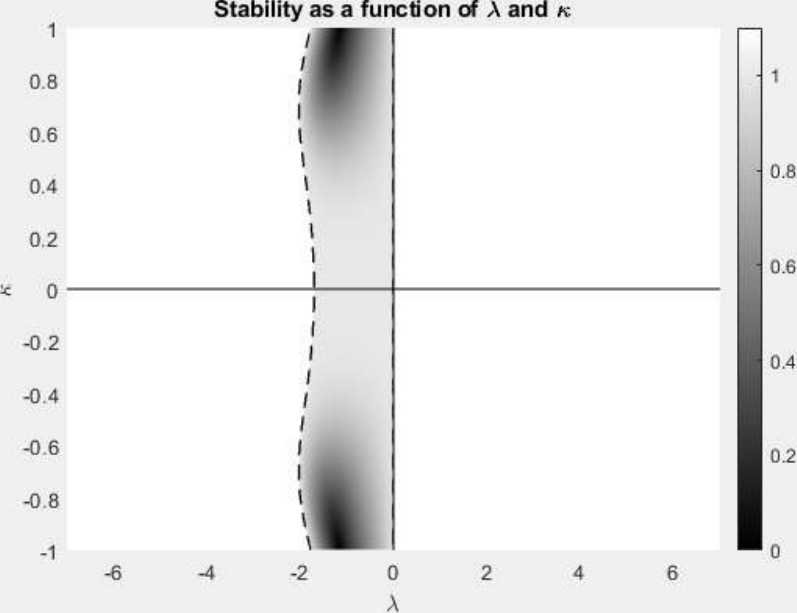}
\endminipage
\caption{RK2; left to right : stability regions for \textbf{Order 3} c, cs, f, zff and zffs. \label{rk23}  }
\end{figure}

\begin{figure}[H]
\minipage{0.19\textwidth}
  \includegraphics[width=\linewidth]{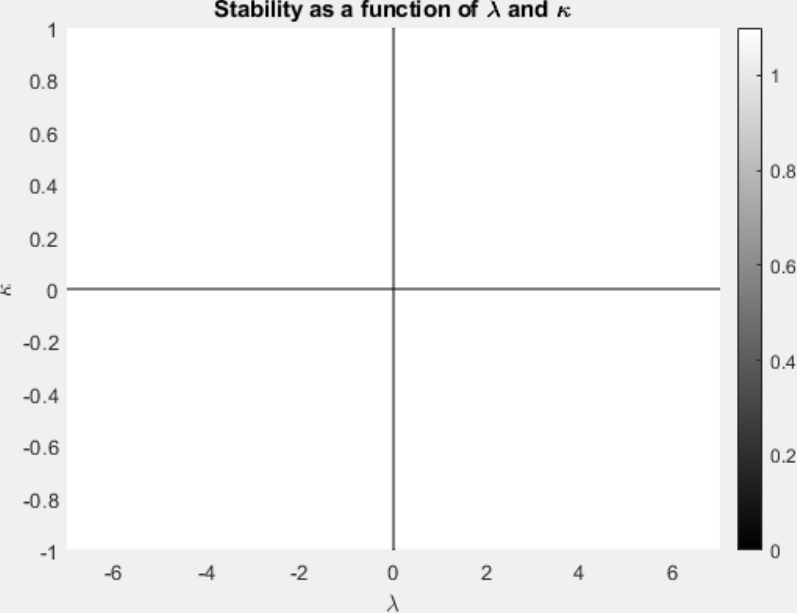}
\endminipage\hfill
\minipage{0.19\textwidth}
  \includegraphics[width=\linewidth]{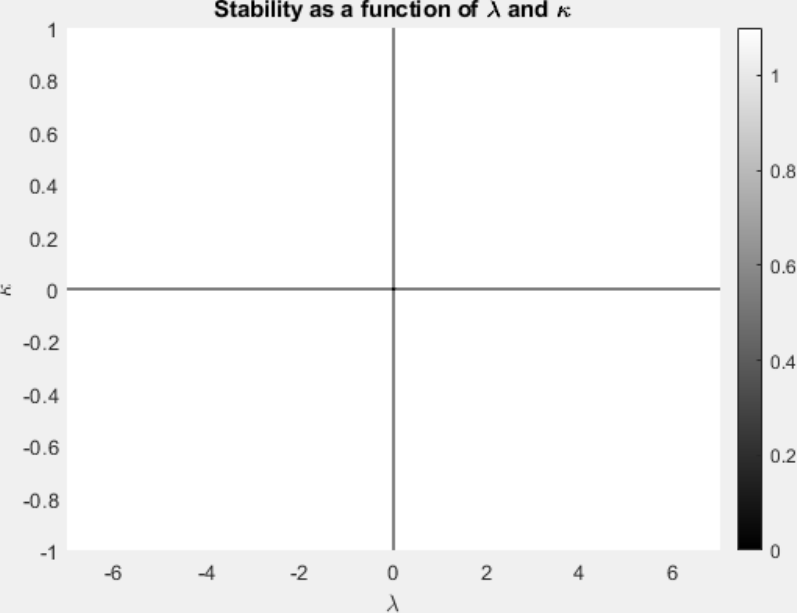}
\endminipage\hfill
\minipage{0.19\textwidth}
  \includegraphics[width=\linewidth]{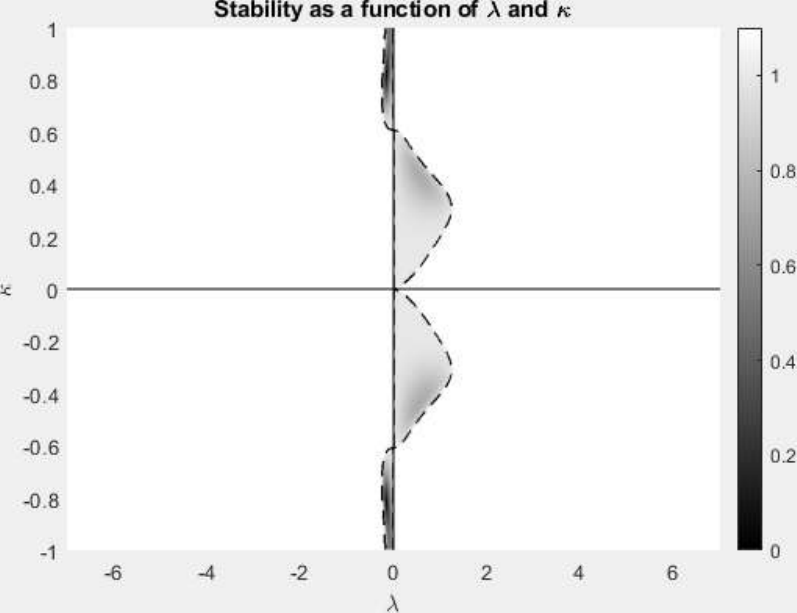}
\endminipage\hfill
\minipage{0.19\textwidth}
  \includegraphics[width=\linewidth]{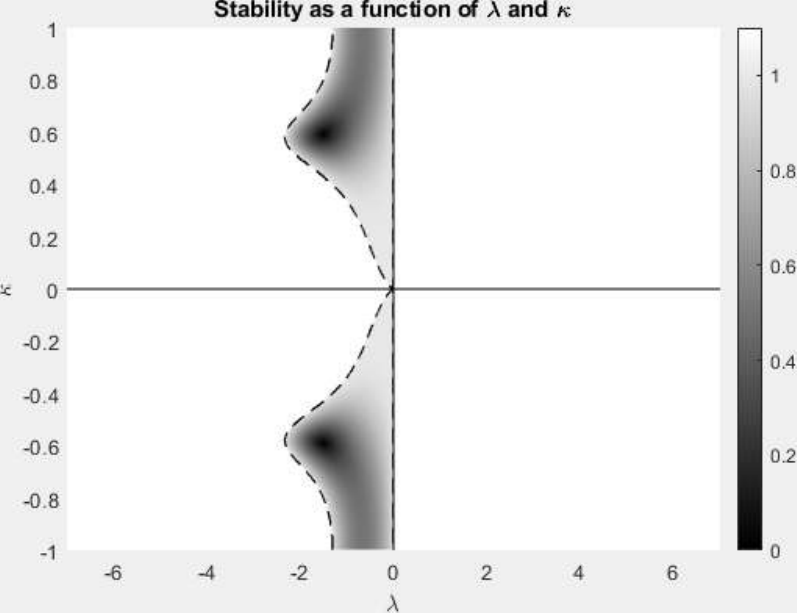}
\endminipage\hfill
\minipage{0.19\textwidth}%
  \includegraphics[width=\linewidth]{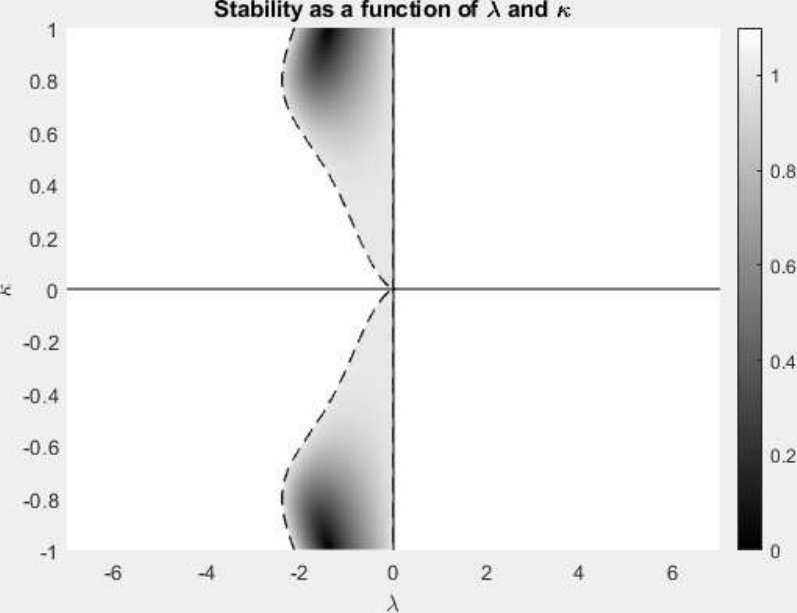}
\endminipage
\caption{RK2; left to right : stability regions for \textbf{Order 4} c, cs, f, zff and zffs. \label{rk24}  }
\end{figure}

\begin{figure}[H]
\minipage{0.19\textwidth}
  \includegraphics[width=\linewidth]{blank.pdf}
\endminipage\hfill
\minipage{0.19\textwidth}
  \includegraphics[width=\linewidth]{blank.pdf}
\endminipage\hfill
\minipage{0.19\textwidth}
  \includegraphics[width=\linewidth]{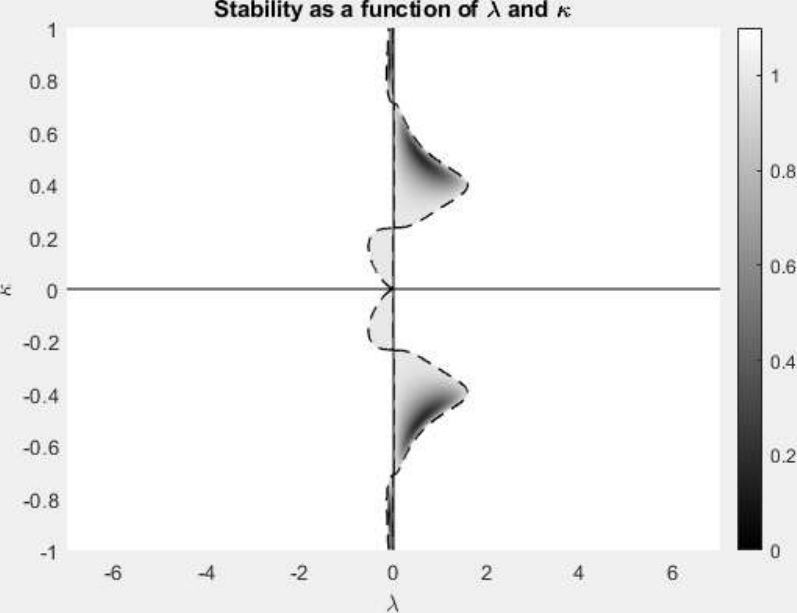}
\endminipage\hfill
\minipage{0.19\textwidth}
  \includegraphics[width=\linewidth]{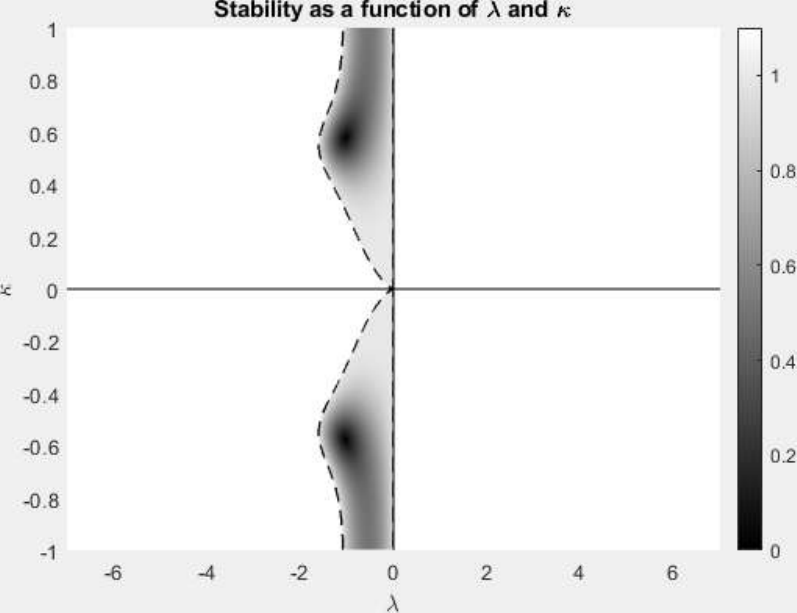}
\endminipage\hfill
\minipage{0.19\textwidth}%
  \includegraphics[width=\linewidth]{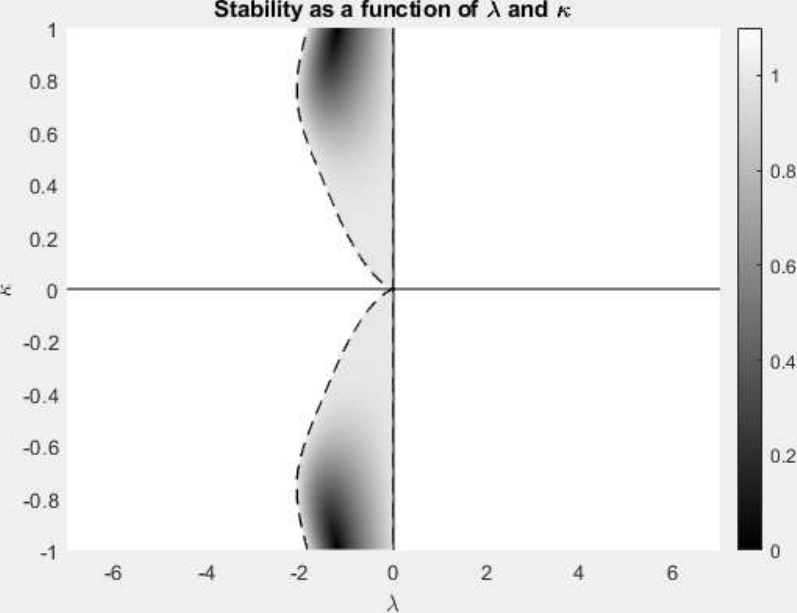}
\endminipage
\caption{RK2; left to right : stability regions for \textbf{Order 5} c, cs, f, zff and zffs. \label{rk25}  }
\end{figure}

\begin{figure}[H]
\minipage{0.19\textwidth}
  \includegraphics[width=\linewidth]{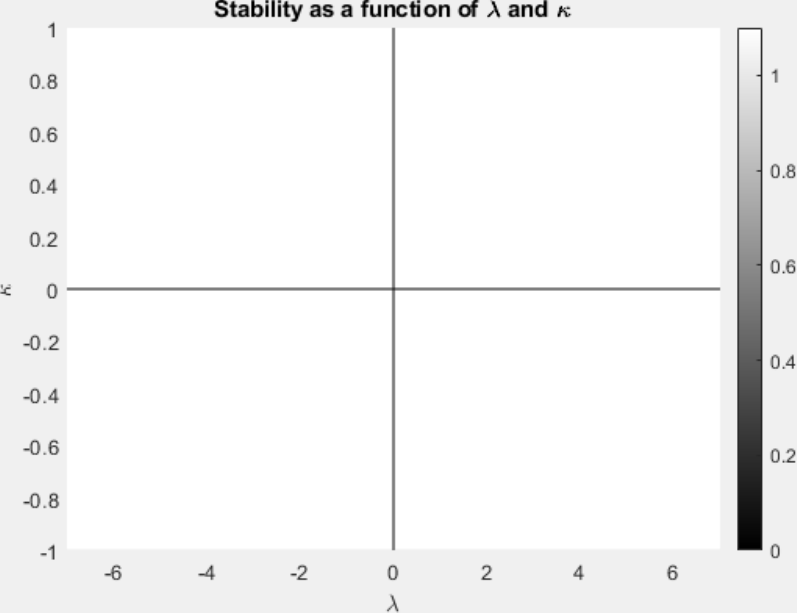}
\endminipage\hfill
\minipage{0.19\textwidth}
  \includegraphics[width=\linewidth]{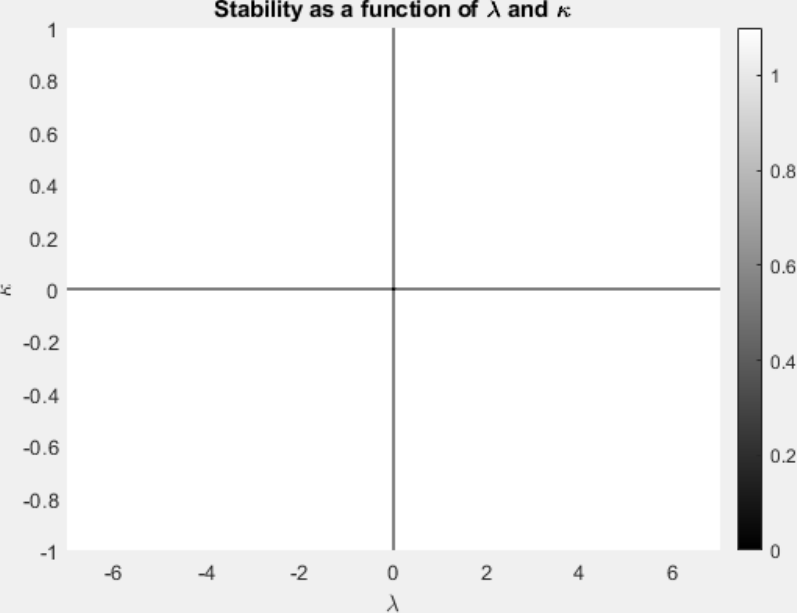}
\endminipage\hfill
\minipage{0.19\textwidth}
  \includegraphics[width=\linewidth]{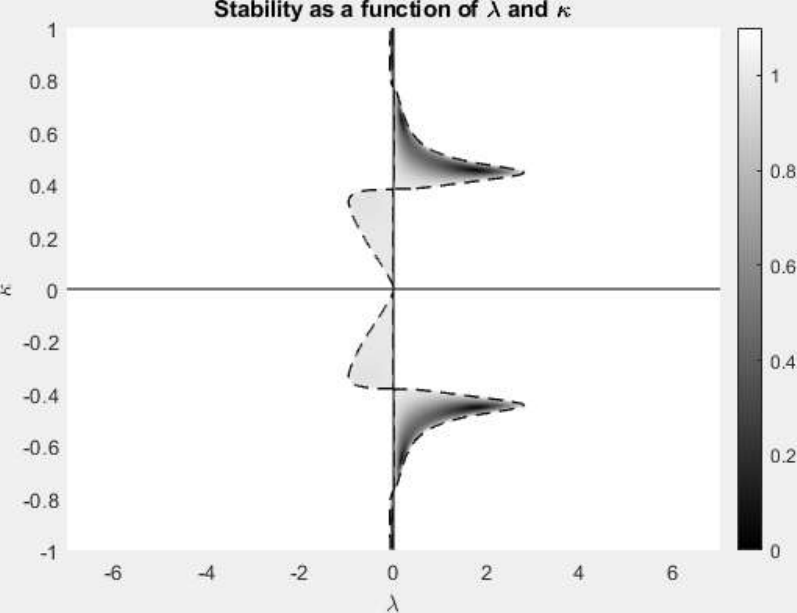}
\endminipage\hfill
\minipage{0.19\textwidth}
  \includegraphics[width=\linewidth]{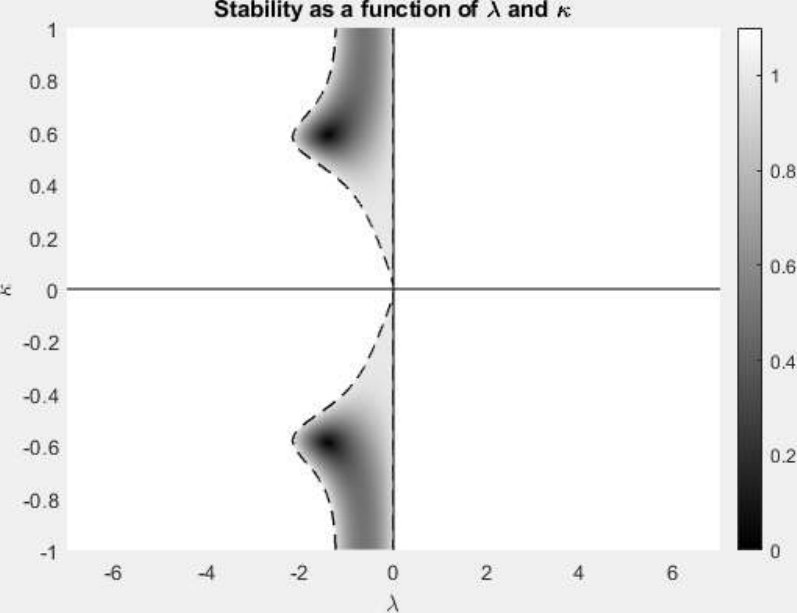}
\endminipage\hfill
\minipage{0.19\textwidth}%
  \includegraphics[width=\linewidth]{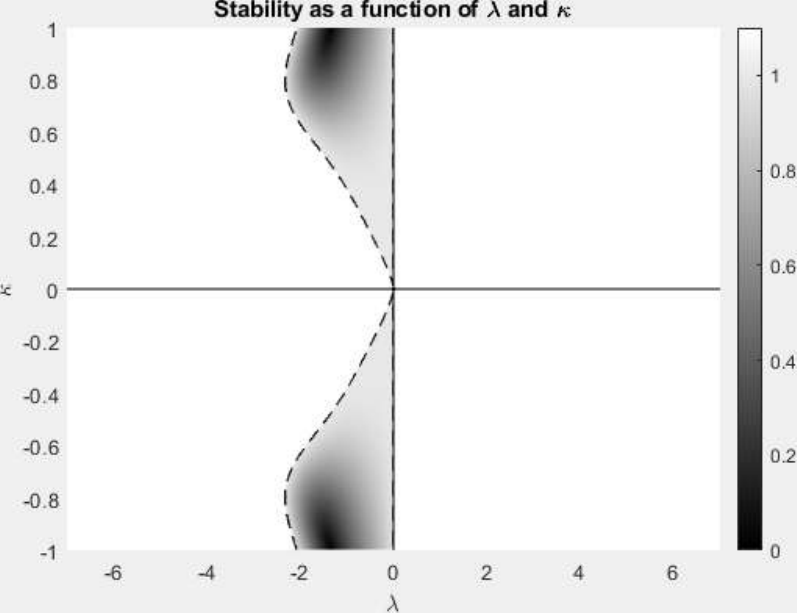}
\endminipage
\caption{RK2; left to right : stability regions for \textbf{Order 6} c, cs, f, zff and zffs. \label{rk26}  }
\end{figure}

Fig. \ref{rk21} to \ref{rk26} show the stability regions in the $\lambda$---$\kappa$ plane 
($\lambda$ and $\kappa$ respectively being  the stability-number and the scaled wavenumber) 
for RK2 time integrator, for spatial orders $N=1$ to $6$, respectively (left to right) for 
centred (c), centred staggered (cs), forward (f), zigzag forward-first (zff) and zigzag 
forward-first staggered (zffs) schemes. The greyed out areas 
represent the couples $(\lambda,\kappa)$ for which the scheme is stable 
(i.e. the amplification factor $|G| \leq 1$), while the dotted contours 
correspond to the critical case $|G| = 1$. A given scheme is conditionally
stable if and only if there exist a $\lambda_c \in \mathds{R}$ such that
the line of equation $\lambda = \lambda_c$ is included in the grey area.

\newpage
\section{Stability regions for RK3}
\begin{figure}[H]
\minipage{0.18\textwidth}
  \includegraphics[width=\linewidth]{blank.pdf}
\endminipage\hfill
\minipage{0.19\textwidth}
  \includegraphics[width=\linewidth]{blank.pdf}
\endminipage\hfill
\minipage{0.19\textwidth}
  \includegraphics[width=\linewidth]{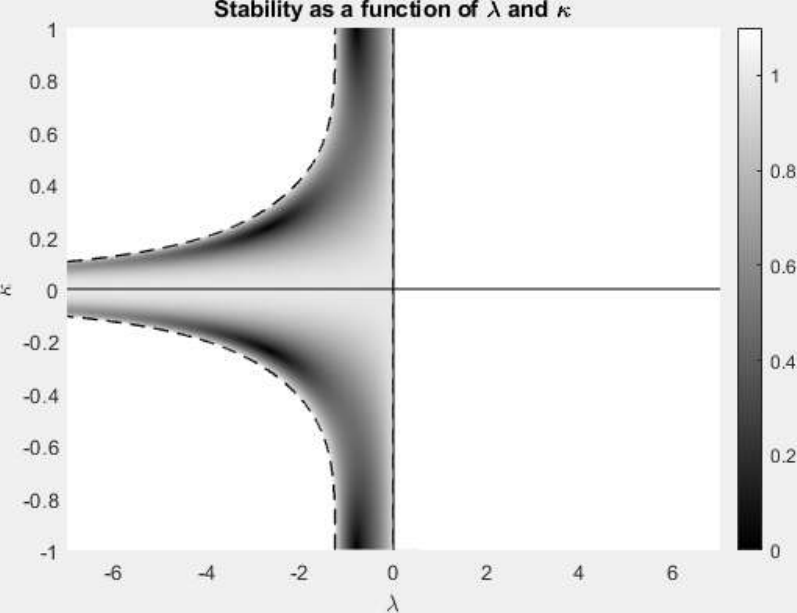}
\endminipage\hfill
\minipage{0.19\textwidth}
  \includegraphics[width=\linewidth]{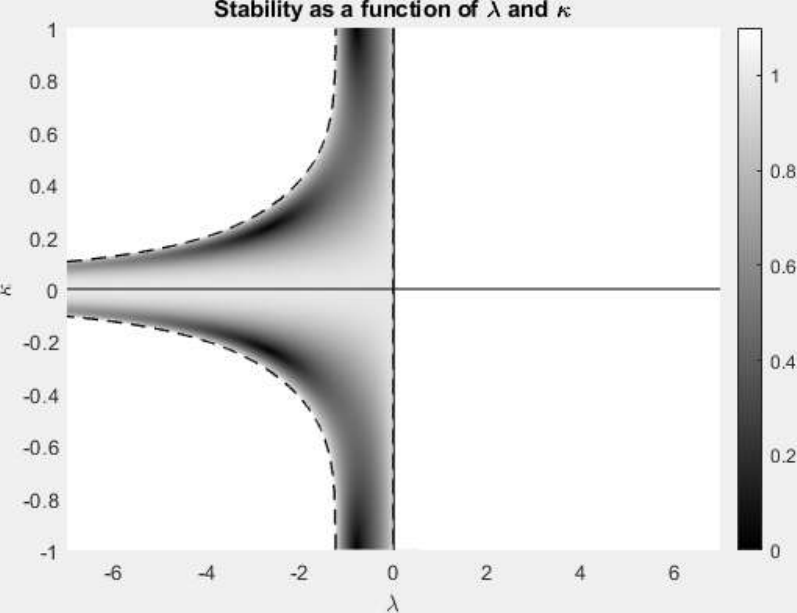}
\endminipage\hfill
\minipage{0.19\textwidth}%
  \includegraphics[width=\linewidth]{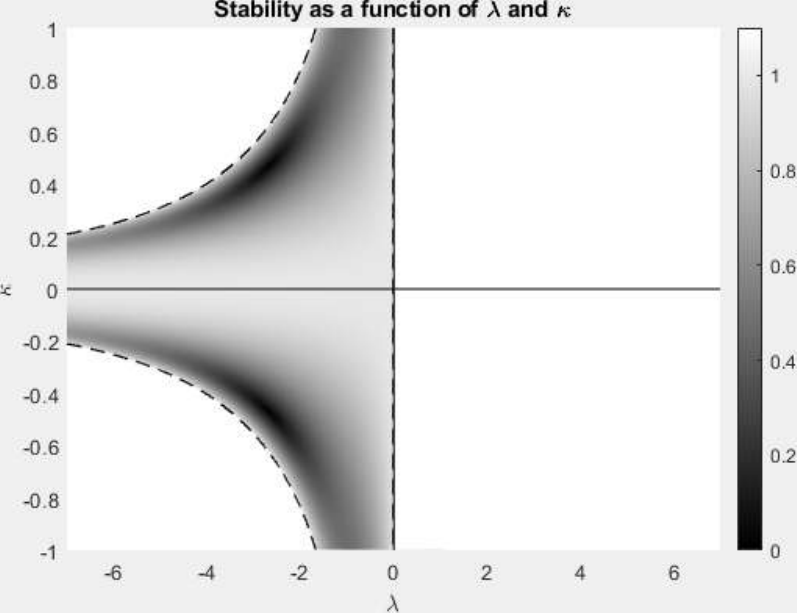}
\endminipage
\caption{RK3; left to right : stability regions for \textbf{Order 1} c, cs, f, zff and zffs.\label{rk31}}
\end{figure}

\begin{figure}[H]
\minipage{0.19\textwidth}
  \includegraphics[width=\linewidth]{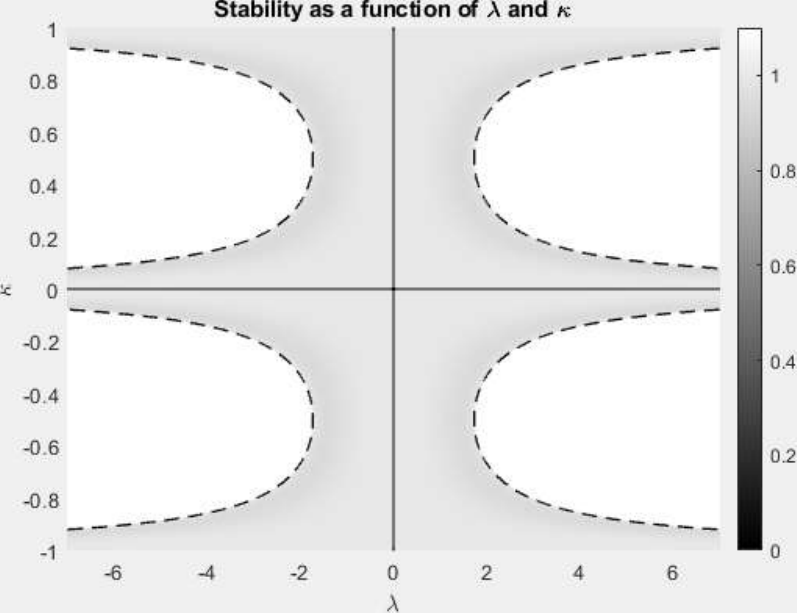}
\endminipage\hfill
\minipage{0.19\textwidth}
  \includegraphics[width=\linewidth]{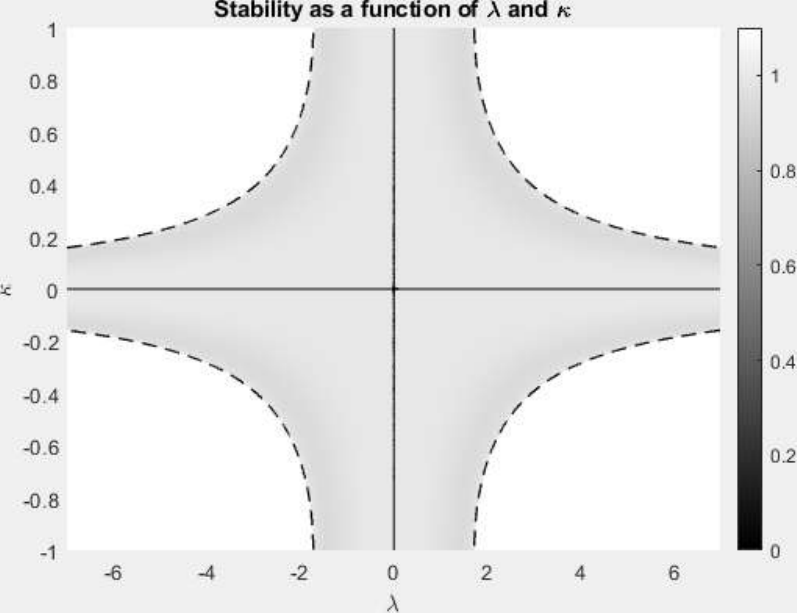}
\endminipage\hfill
\minipage{0.19\textwidth}
  \includegraphics[width=\linewidth]{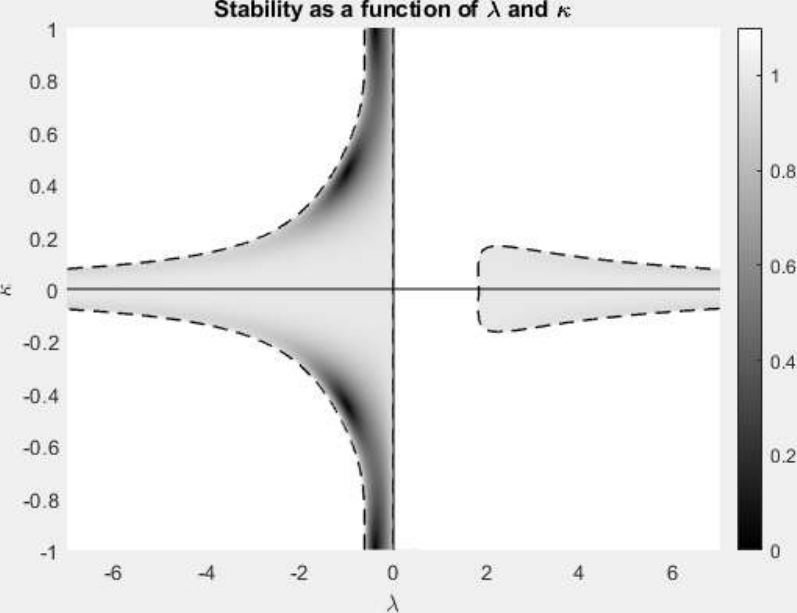}
\endminipage\hfill
\minipage{0.19\textwidth}
  \includegraphics[width=\linewidth]{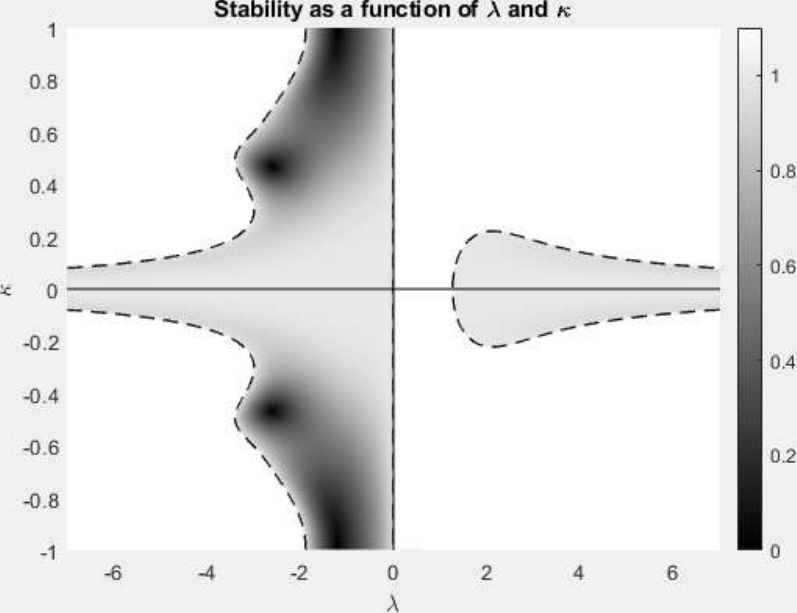}
\endminipage\hfill
\minipage{0.19\textwidth}%
  \includegraphics[width=\linewidth]{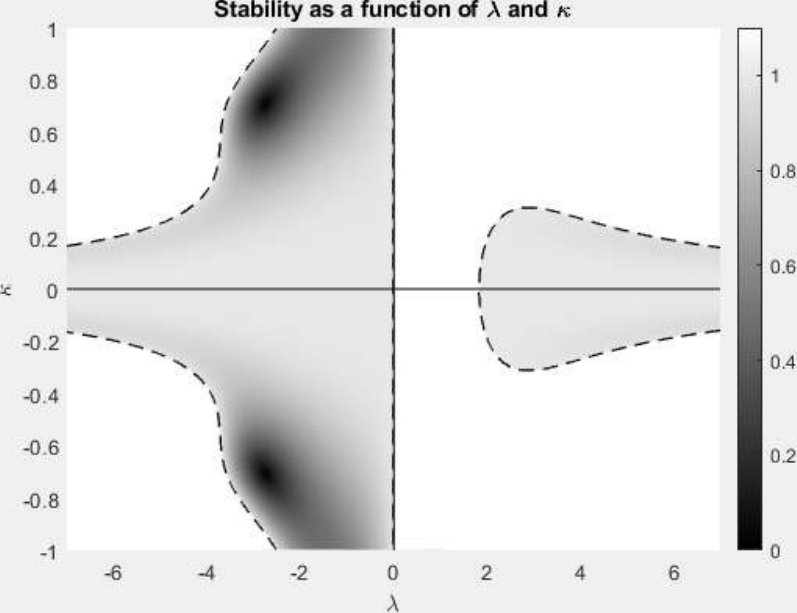}
\endminipage
\caption{RK3; left to right : stability regions for \textbf{Order 2} c, cs, f, zff and zffs.\label{rk32}}
\end{figure}

\begin{figure}[H]
\minipage{0.19\textwidth}
  \includegraphics[width=\linewidth]{blank.pdf}
\endminipage\hfill
\minipage{0.19\textwidth}
  \includegraphics[width=\linewidth]{blank.pdf}
\endminipage\hfill
\minipage{0.19\textwidth}
  \includegraphics[width=\linewidth]{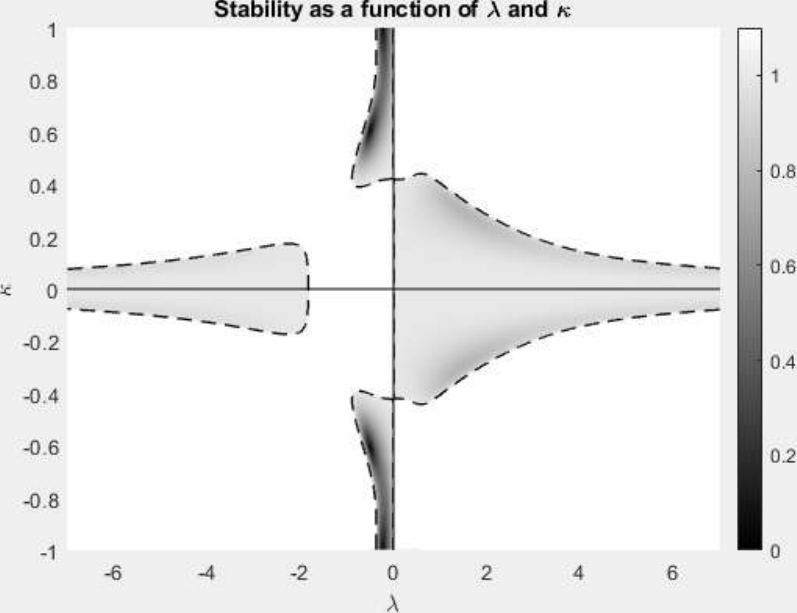}
\endminipage\hfill
\minipage{0.19\textwidth}
  \includegraphics[width=\linewidth]{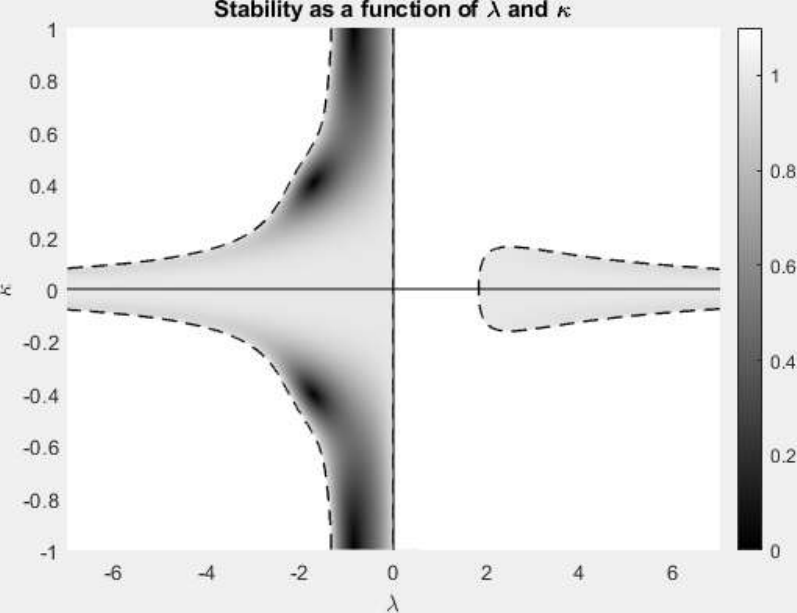}
\endminipage\hfill
\minipage{0.19\textwidth}%
  \includegraphics[width=\linewidth]{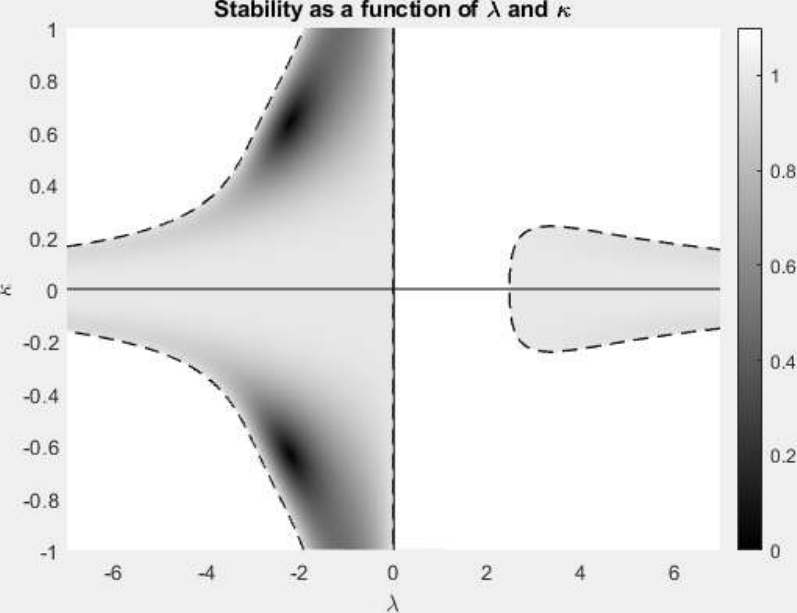}
\endminipage
\caption{RK3; left to right : stability regions for \textbf{Order 3} c, cs, f, zff and zffs.\label{rk33}}
\end{figure}

\begin{figure}[H]
\minipage{0.19\textwidth}
  \includegraphics[width=\linewidth]{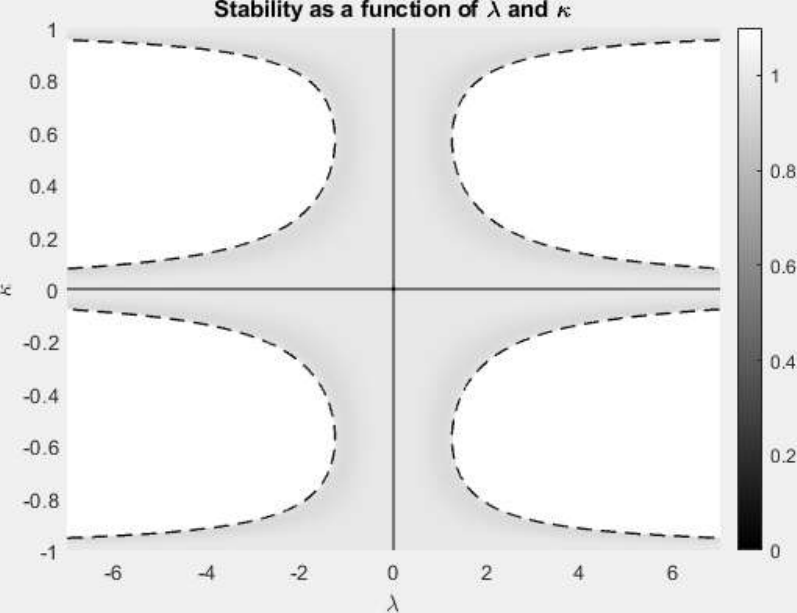}
\endminipage\hfill
\minipage{0.19\textwidth}
  \includegraphics[width=\linewidth]{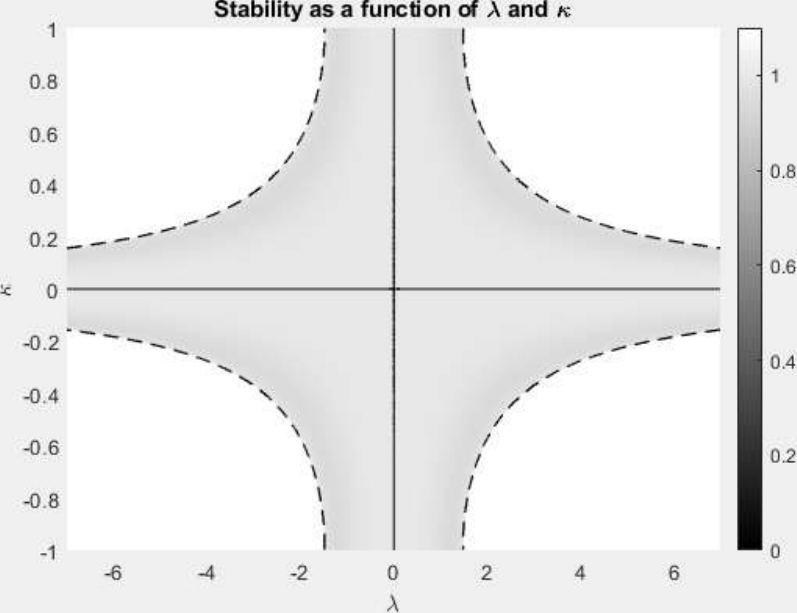}
\endminipage\hfill
\minipage{0.19\textwidth}
  \includegraphics[width=\linewidth]{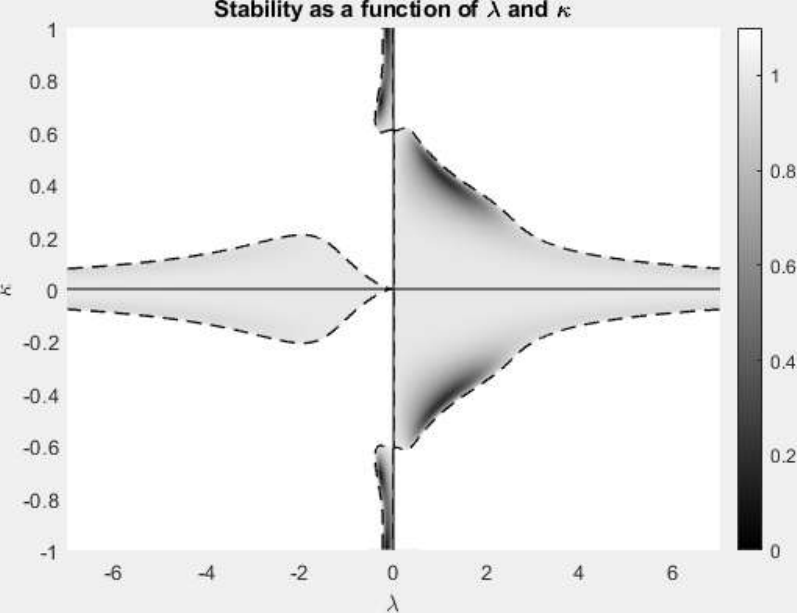}
\endminipage\hfill
\minipage{0.19\textwidth}
  \includegraphics[width=\linewidth]{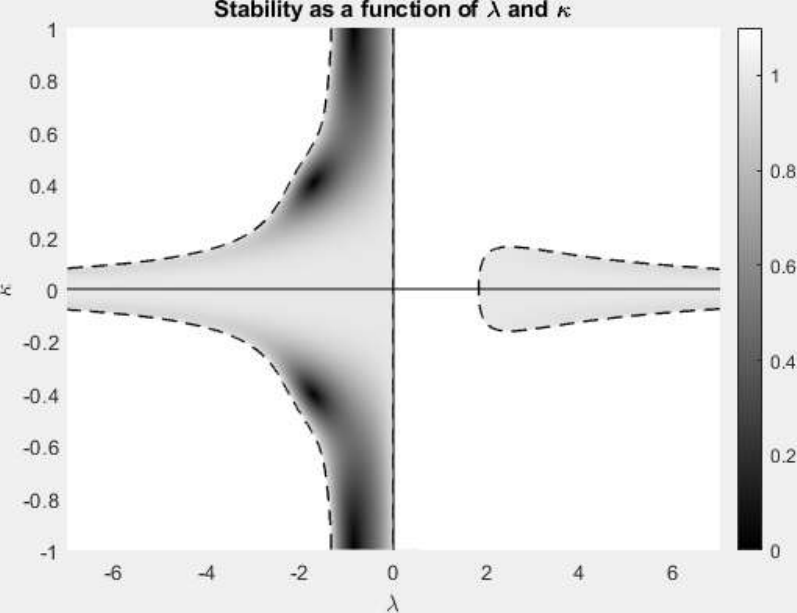}
\endminipage\hfill
\minipage{0.19\textwidth}%
  \includegraphics[width=\linewidth]{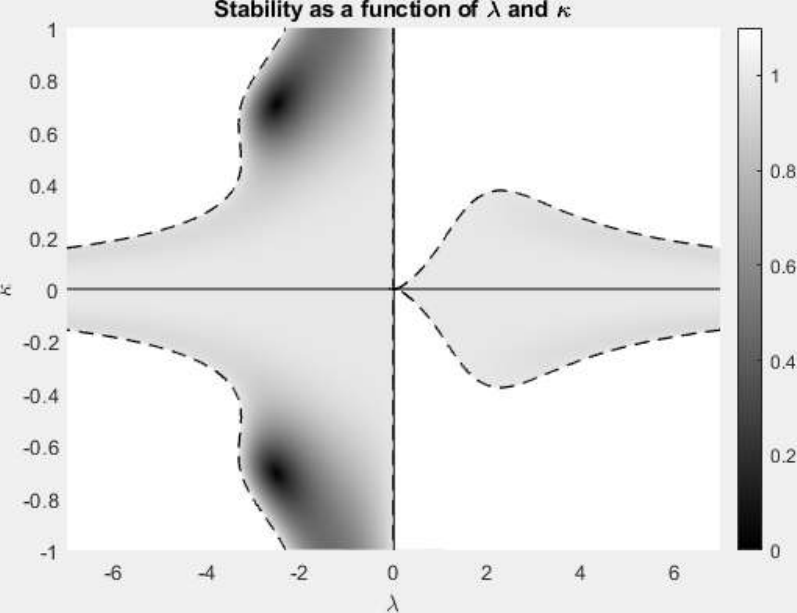}
\endminipage
\caption{RK3; left to right : stability regions for \textbf{Order 4} c, cs, f, zff and zffs.\label{rk34}}
\end{figure}

\begin{figure}[H]
\minipage{0.19\textwidth}
  \includegraphics[width=\linewidth]{blank.pdf}
\endminipage\hfill
\minipage{0.19\textwidth}
  \includegraphics[width=\linewidth]{blank.pdf}
\endminipage\hfill
\minipage{0.19\textwidth}
  \includegraphics[width=\linewidth]{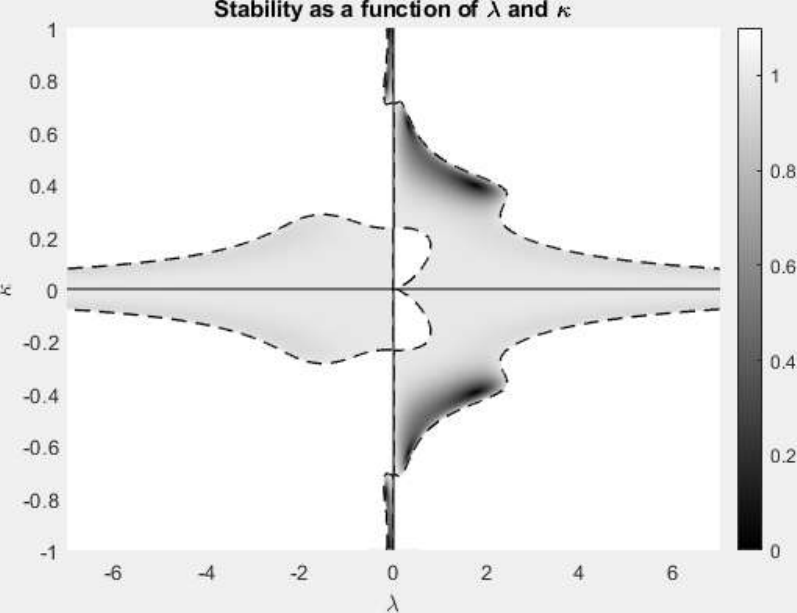}
\endminipage\hfill
\minipage{0.19\textwidth}
  \includegraphics[width=\linewidth]{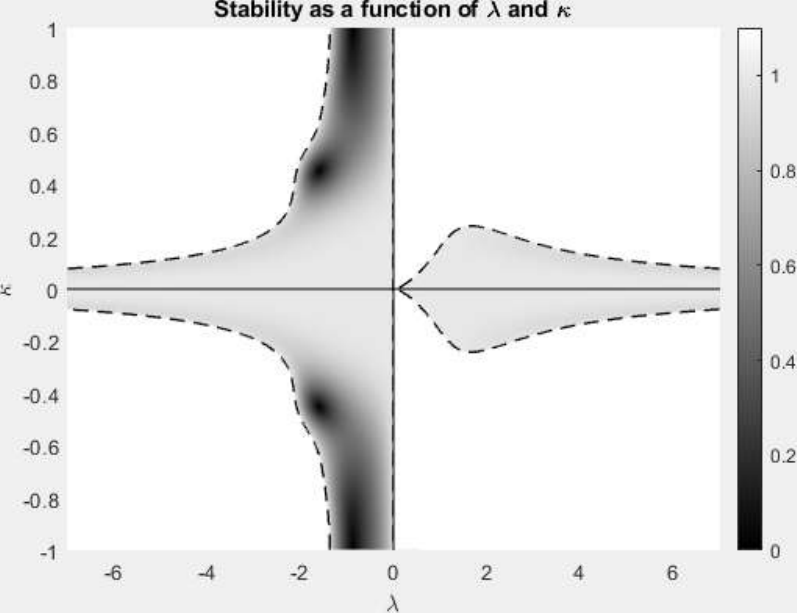}
\endminipage\hfill
\minipage{0.19\textwidth}%
  \includegraphics[width=\linewidth]{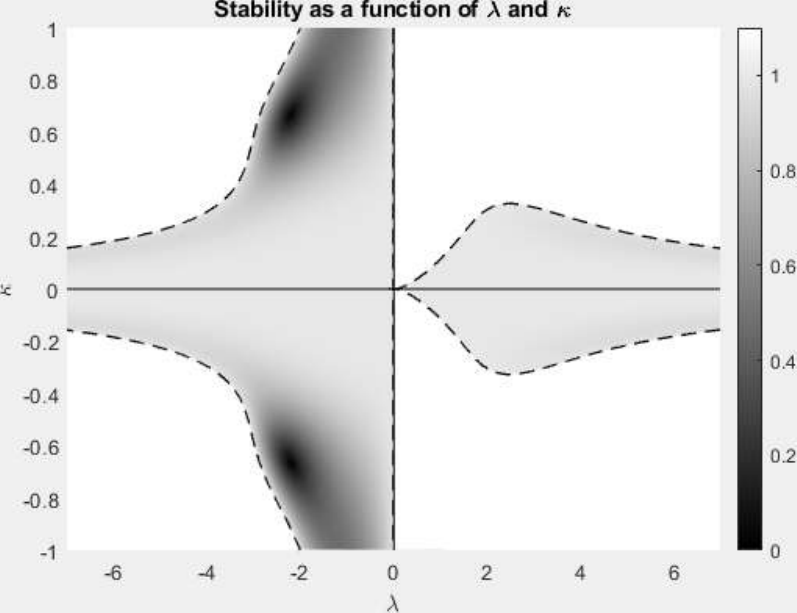}
\endminipage
\caption{RK3; left to right : stability regions for \textbf{Order 5} c, cs, f, zff and zffs.\label{rk35}}
\end{figure}

\begin{figure}[H]
\minipage{0.19\textwidth}
  \includegraphics[width=\linewidth]{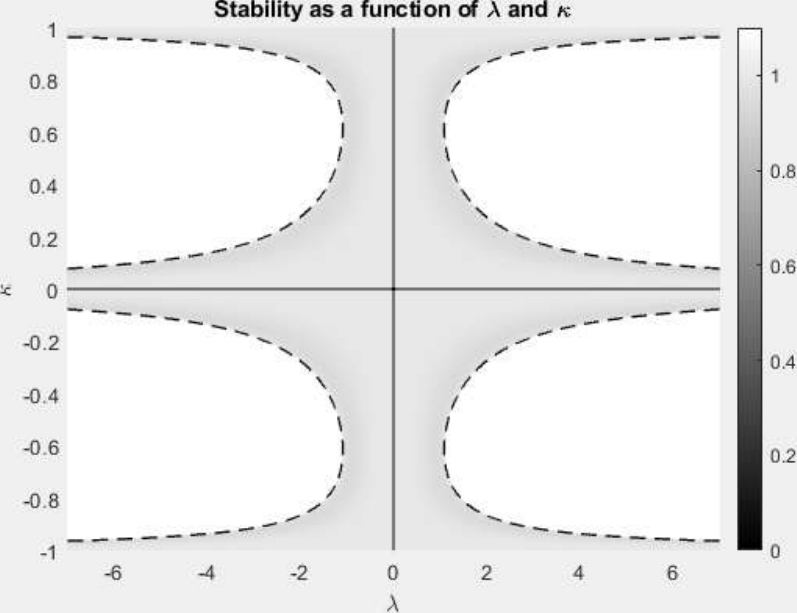}
\endminipage\hfill
\minipage{0.19\textwidth}
  \includegraphics[width=\linewidth]{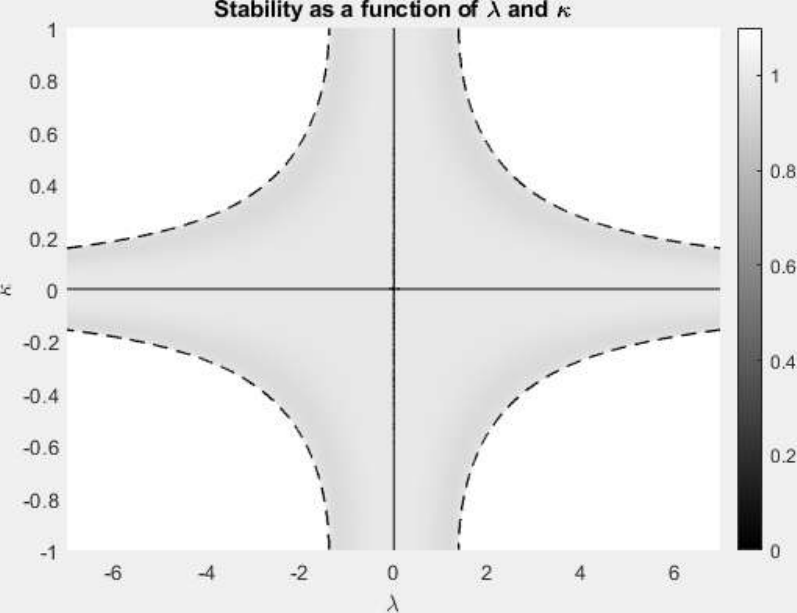}
\endminipage\hfill
\minipage{0.19\textwidth}
  \includegraphics[width=\linewidth]{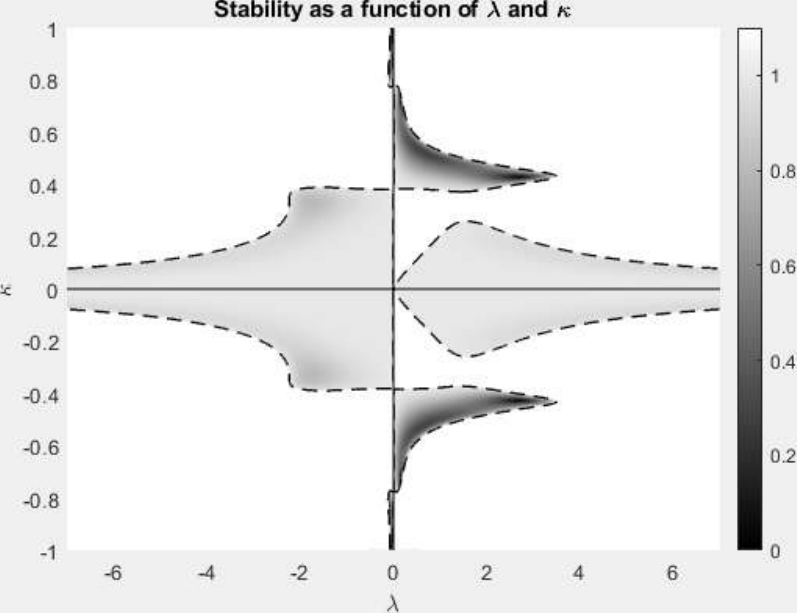}
\endminipage\hfill
\minipage{0.19\textwidth}
  \includegraphics[width=\linewidth]{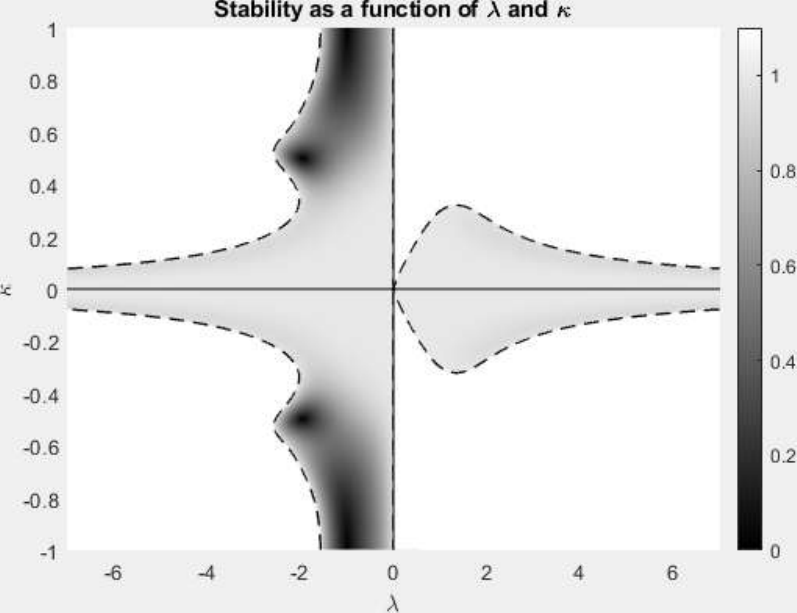}
\endminipage\hfill
\minipage{0.19\textwidth}%
  \includegraphics[width=\linewidth]{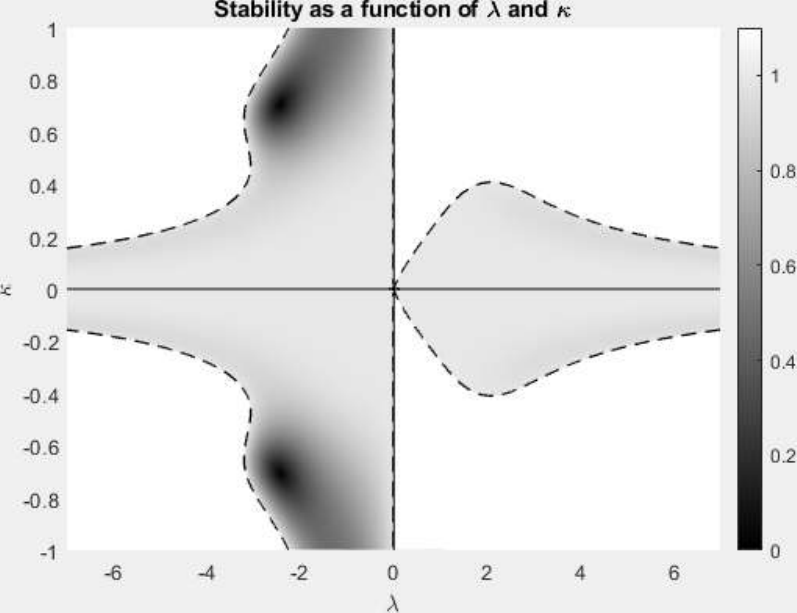}
\endminipage
\caption{RK3; left to right : stability regions for \textbf{Order 6} c, cs, f, zff and zffs.\label{rk36}}
\end{figure}

Fig. \ref{rk31} to \ref{rk36} show the stability regions in the $\lambda$---$\kappa$ plane 
($\lambda$ and $\kappa$ respectively being the stability-number and the scaled wavenumber) for RK3 time integrator, 
for spatial orders $N=1$ to $6$, 
respectively (left to right) for centred (c), centred staggered (cs), forward (f), zigzag forward-first (zff) and zigzag forward-first staggered (zffs) schemes.
The greyed out areas 
represent the couples $(\lambda,\kappa)$ for which the scheme is stable 
(i.e. the amplification factor $|G| \leq 1$), while the dotted contours 
correspond to the critical case $|G| = 1$. A given scheme is conditionally
stable if and only if there exist a $\lambda_c \in \mathds{R}$ such that
the line of equation $\lambda = \lambda_c$ is included in the grey area.

\newpage
\section{Stability regions for RK4}
\begin{figure}[H]
\minipage{0.18\textwidth}
  \includegraphics[width=\linewidth]{blank.pdf}
\endminipage\hfill
\minipage{0.19\textwidth}
  \includegraphics[width=\linewidth]{blank.pdf}
\endminipage\hfill
\minipage{0.19\textwidth}
  \includegraphics[width=\linewidth]{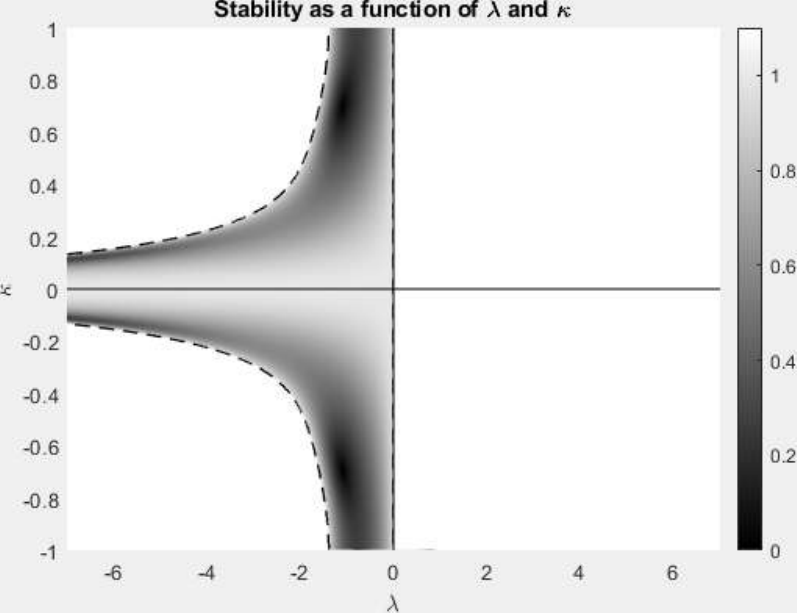}
\endminipage\hfill
\minipage{0.19\textwidth}
  \includegraphics[width=\linewidth]{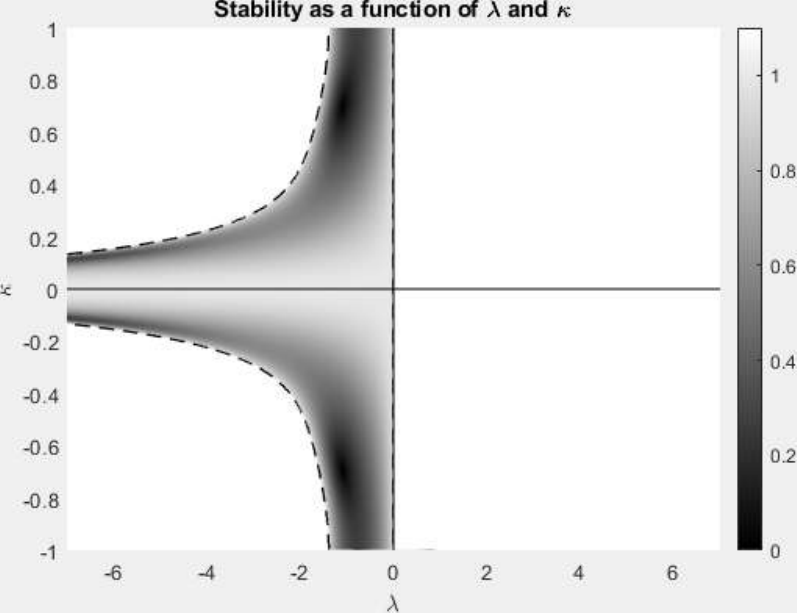}
\endminipage\hfill
\minipage{0.19\textwidth}%
  \includegraphics[width=\linewidth]{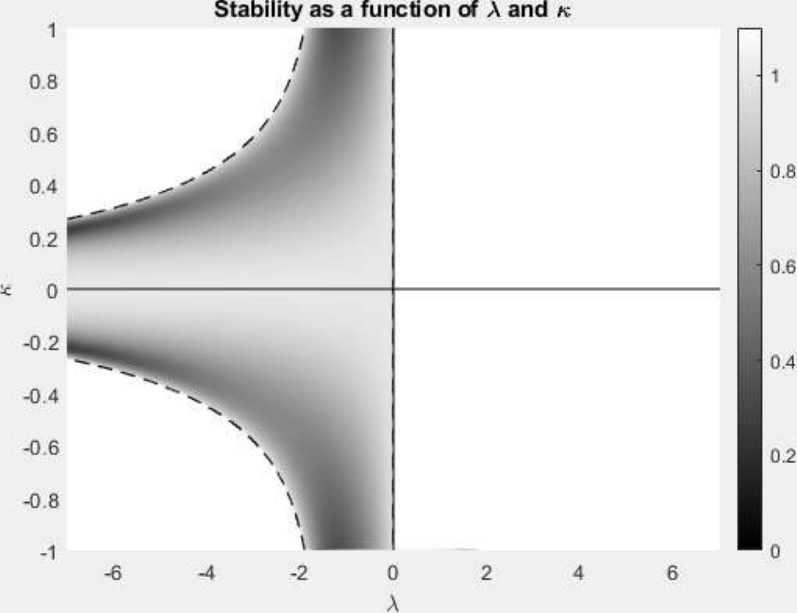}
\endminipage
\caption{RK4; left to right : stability regions for \textbf{Order 1} c, cs, f, zff and zffs.\label{rk41}}
\end{figure}

\begin{figure}[H]
\minipage{0.19\textwidth}
  \includegraphics[width=\linewidth]{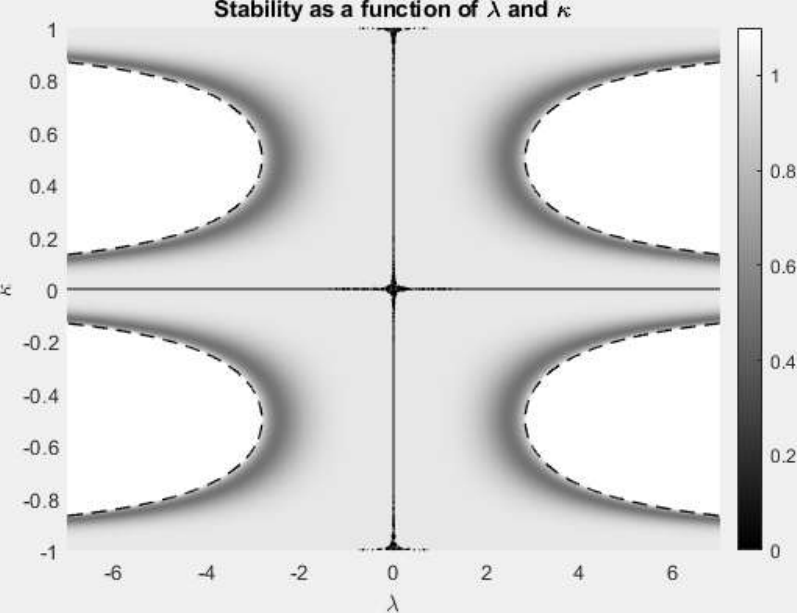}
\endminipage\hfill
\minipage{0.19\textwidth}
  \includegraphics[width=\linewidth]{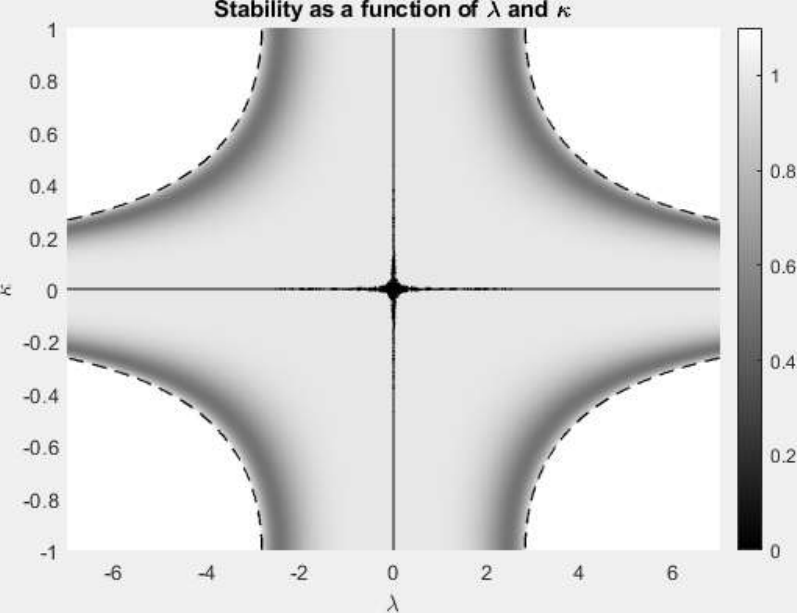}
\endminipage\hfill
\minipage{0.19\textwidth}
  \includegraphics[width=\linewidth]{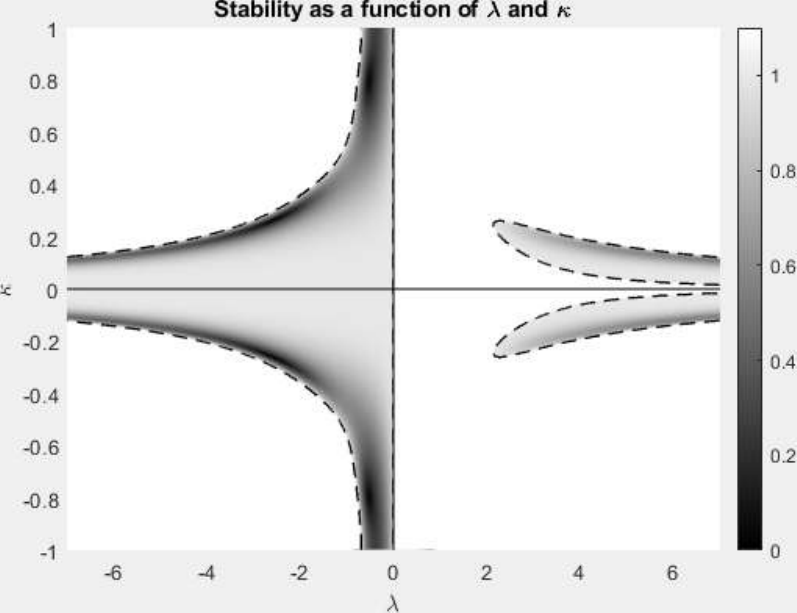}
\endminipage\hfill
\minipage{0.19\textwidth}
  \includegraphics[width=\linewidth]{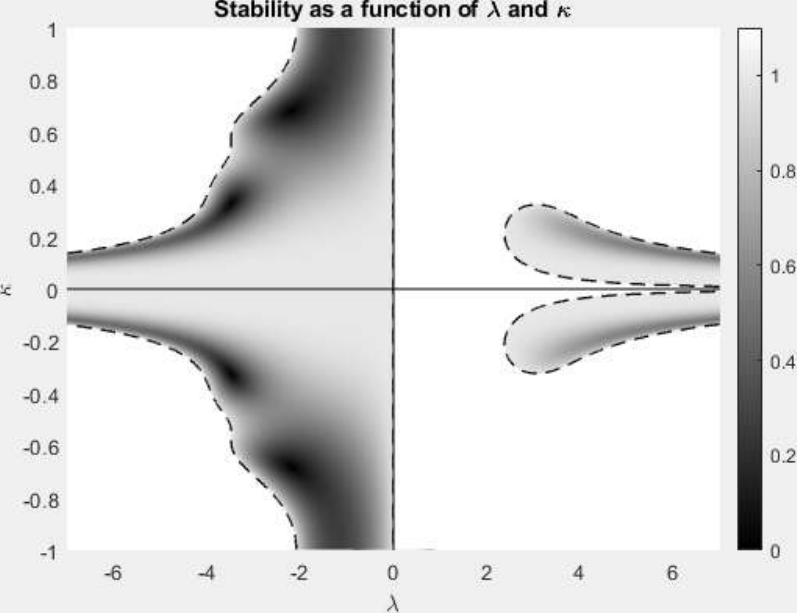}
\endminipage\hfill
\minipage{0.19\textwidth}%
  \includegraphics[width=\linewidth]{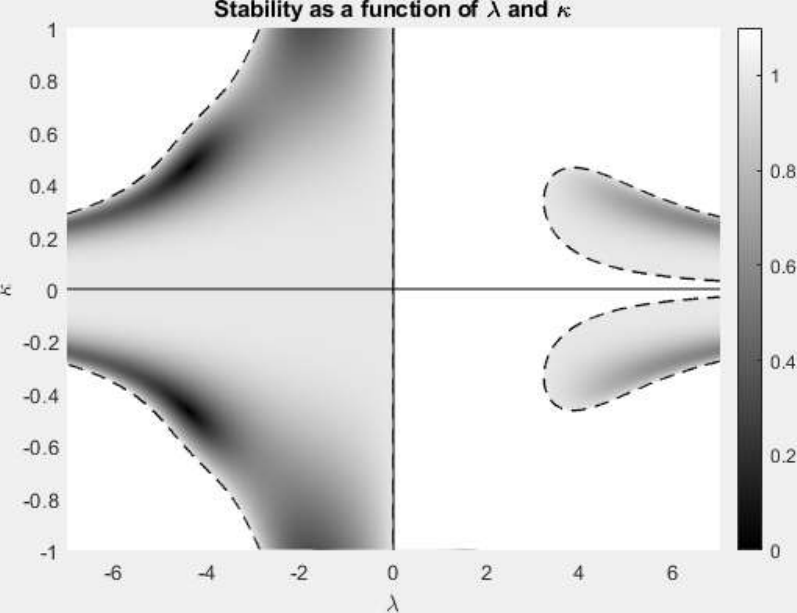}
\endminipage
\caption{RK4; left to right : stability regions for \textbf{Order 2} c, cs, f, zff and zffs.\label{rk42}}
\end{figure}

\begin{figure}[H]
\minipage{0.19\textwidth}
  \includegraphics[width=\linewidth]{blank.pdf}
\endminipage\hfill
\minipage{0.19\textwidth}
  \includegraphics[width=\linewidth]{blank.pdf}
\endminipage\hfill
\minipage{0.19\textwidth}
  \includegraphics[width=\linewidth]{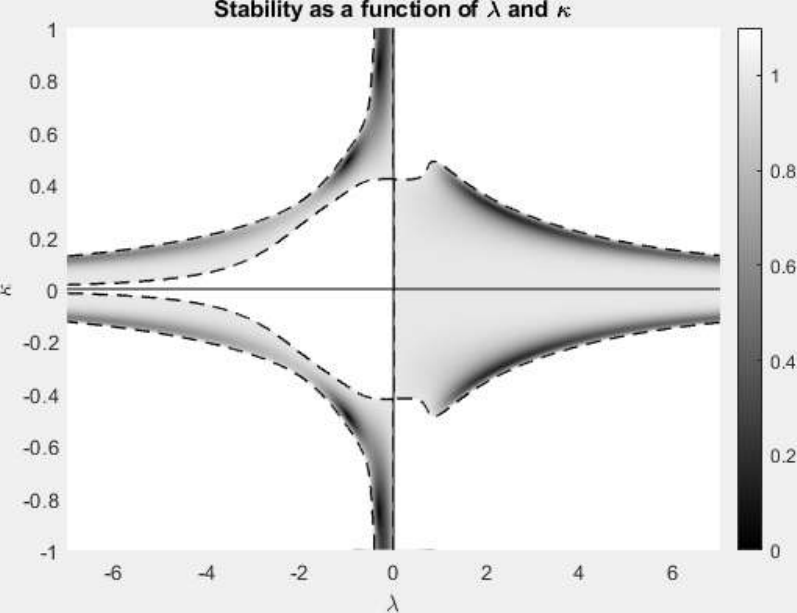}
\endminipage\hfill
\minipage{0.19\textwidth}
  \includegraphics[width=\linewidth]{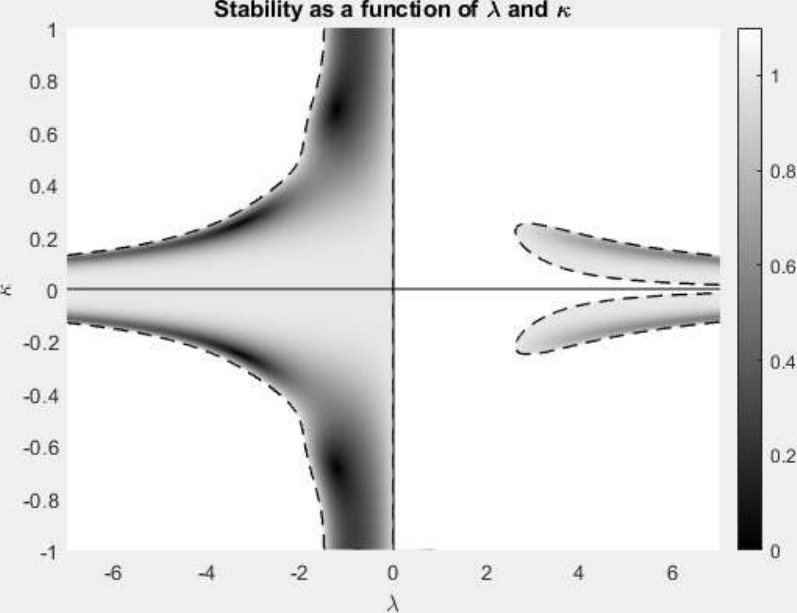}
\endminipage\hfill
\minipage{0.19\textwidth}%
  \includegraphics[width=\linewidth]{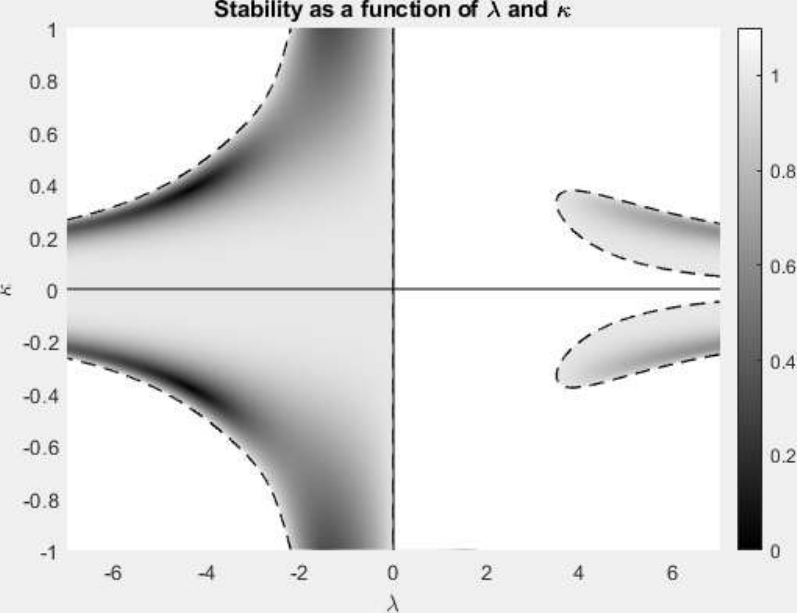}
\endminipage
\caption{RK4; left to right : stability regions for \textbf{Order 3} c, cs, f, zff and zffs.\label{rk43}}
\end{figure}

\begin{figure}[H]
\minipage{0.19\textwidth}
  \includegraphics[width=\linewidth]{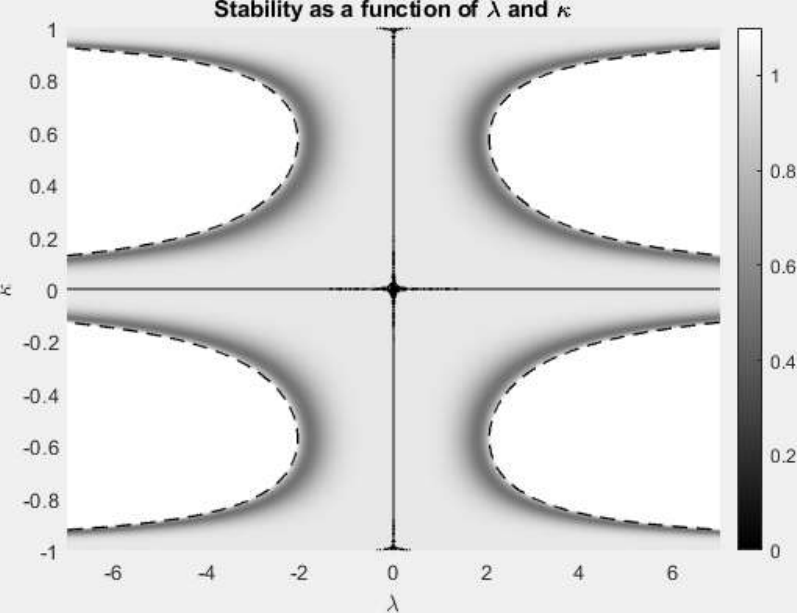}
\endminipage\hfill
\minipage{0.19\textwidth}
  \includegraphics[width=\linewidth]{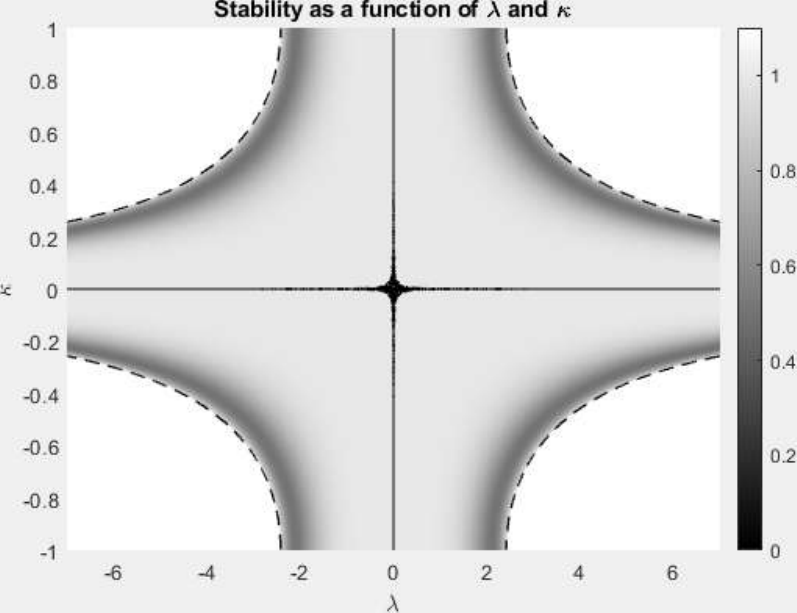}
\endminipage\hfill
\minipage{0.19\textwidth}
  \includegraphics[width=\linewidth]{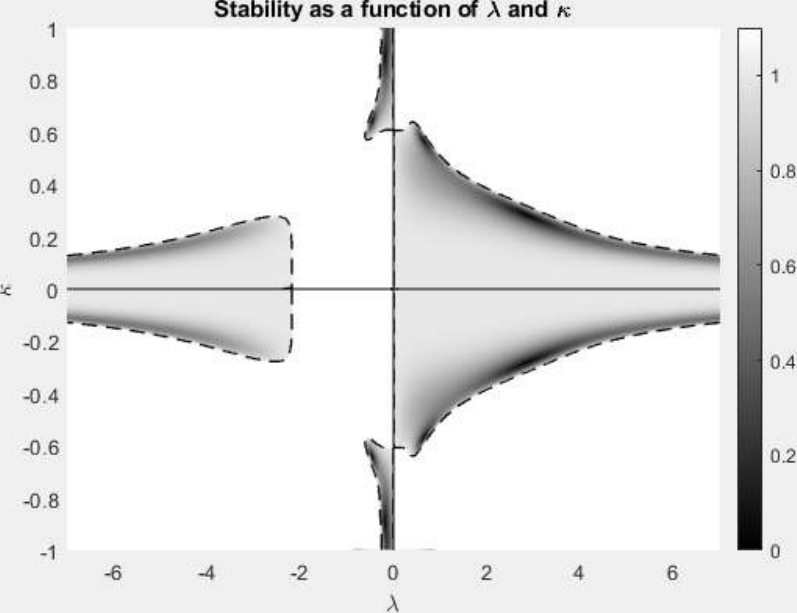}
\endminipage\hfill
\minipage{0.19\textwidth}
  \includegraphics[width=\linewidth]{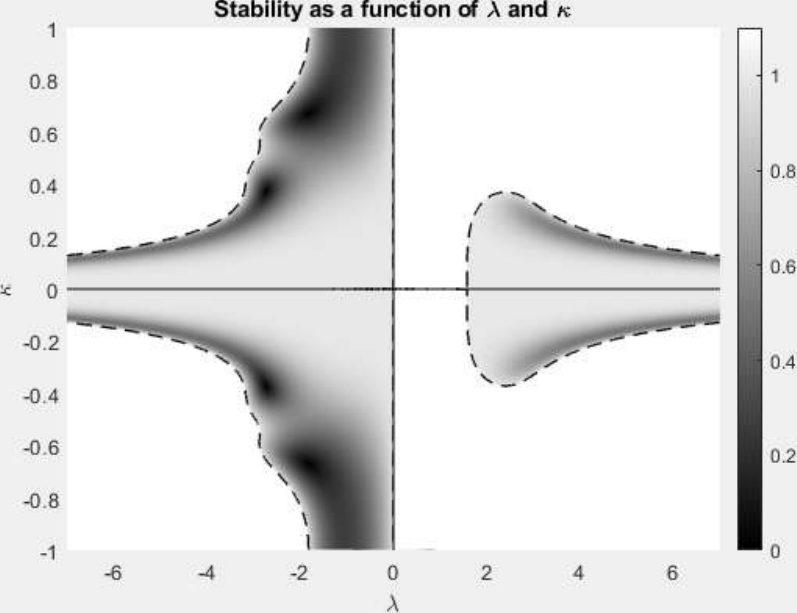}
\endminipage\hfill
\minipage{0.19\textwidth}%
  \includegraphics[width=\linewidth]{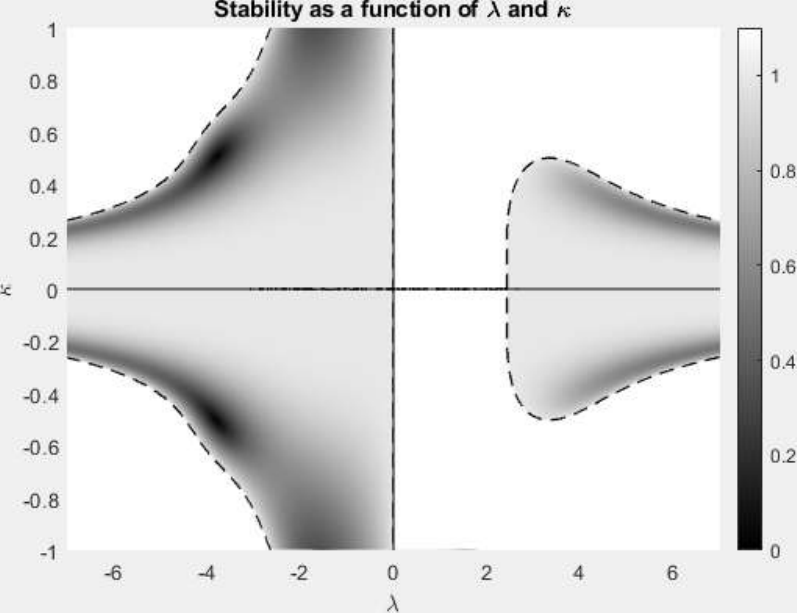}
\endminipage
\caption{RK4; left to right : stability regions for \textbf{Order 4} c, cs, f, zff and zffs.\label{rk44}}
\end{figure}

\begin{figure}[H]
\minipage{0.19\textwidth}
  \includegraphics[width=\linewidth]{blank.pdf}
\endminipage\hfill
\minipage{0.19\textwidth}
  \includegraphics[width=\linewidth]{blank.pdf}
\endminipage\hfill
\minipage{0.19\textwidth}
  \includegraphics[width=\linewidth]{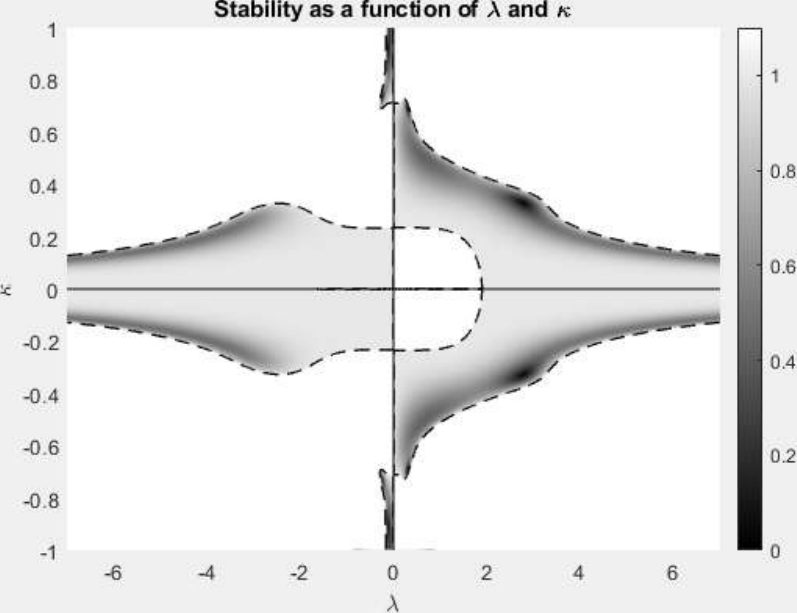}
\endminipage\hfill
\minipage{0.19\textwidth}
  \includegraphics[width=\linewidth]{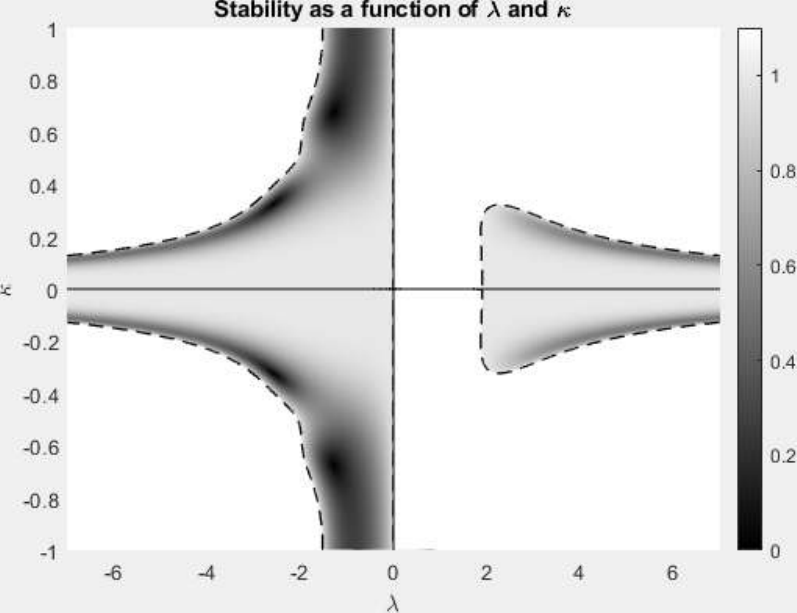}
\endminipage\hfill
\minipage{0.19\textwidth}%
  \includegraphics[width=\linewidth]{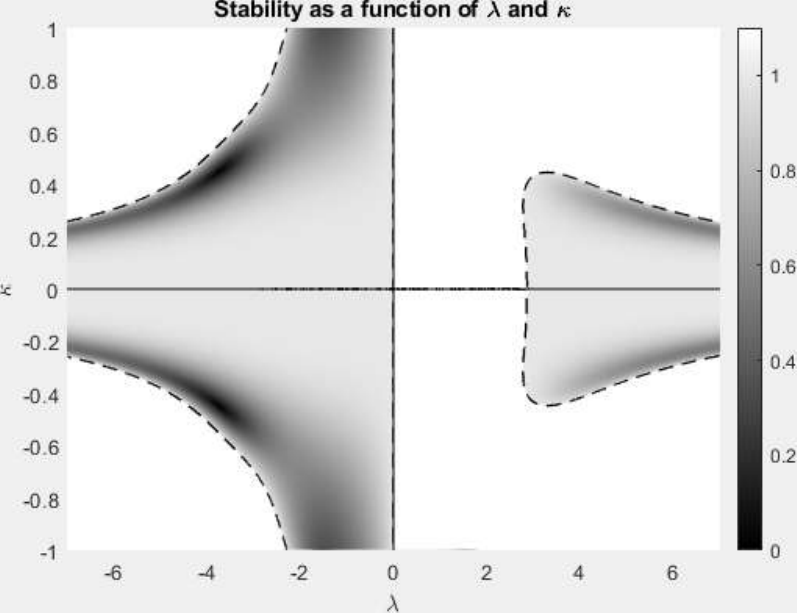}
\endminipage
\caption{RK4; left to right : stability regions for \textbf{Order 5} c, cs, f, zff and zffs.\label{rk45}}
\end{figure}

\begin{figure}[H]
\minipage{0.19\textwidth}
  \includegraphics[width=\linewidth]{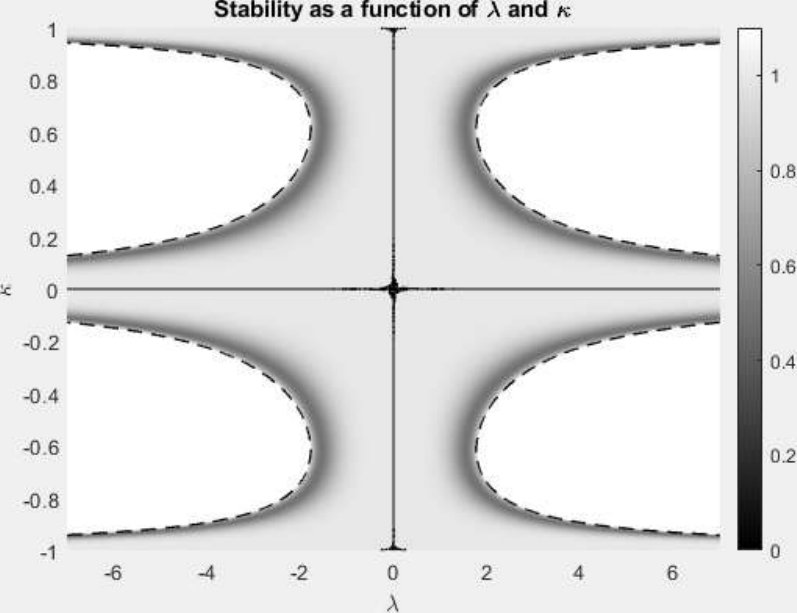}
\endminipage\hfill
\minipage{0.19\textwidth}
  \includegraphics[width=\linewidth]{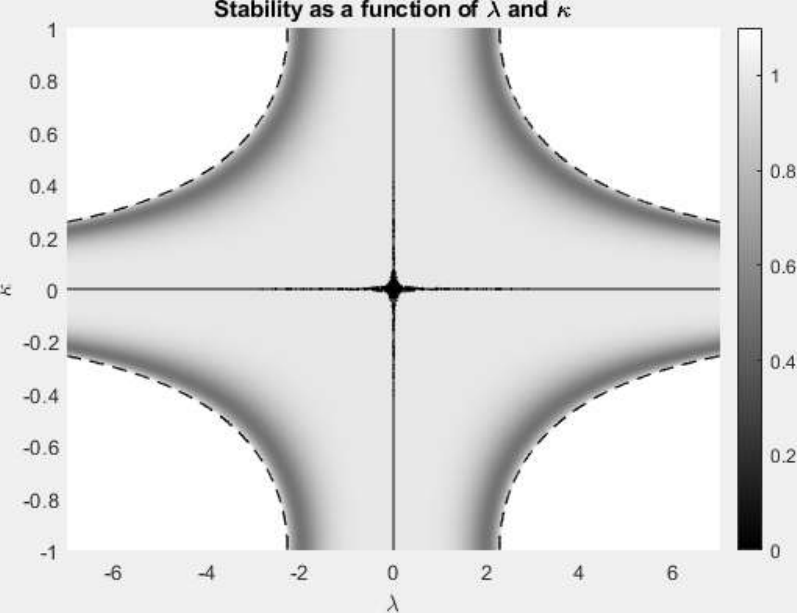}
\endminipage\hfill
\minipage{0.19\textwidth}
  \includegraphics[width=\linewidth]{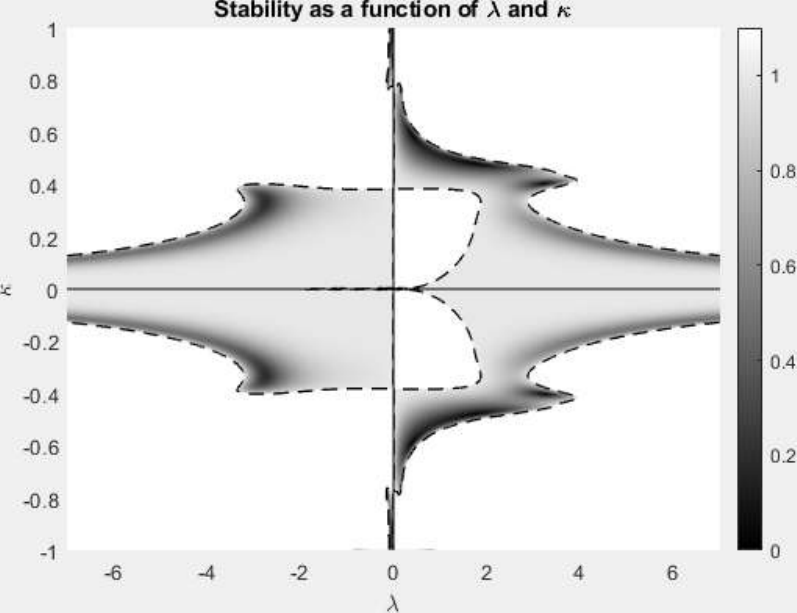}
\endminipage\hfill
\minipage{0.19\textwidth}
  \includegraphics[width=\linewidth]{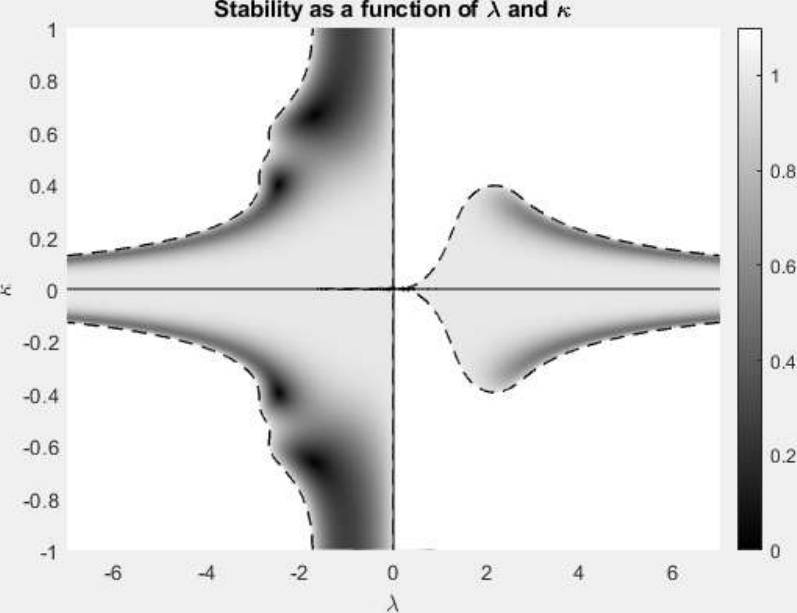}
\endminipage\hfill
\minipage{0.19\textwidth}%
  \includegraphics[width=\linewidth]{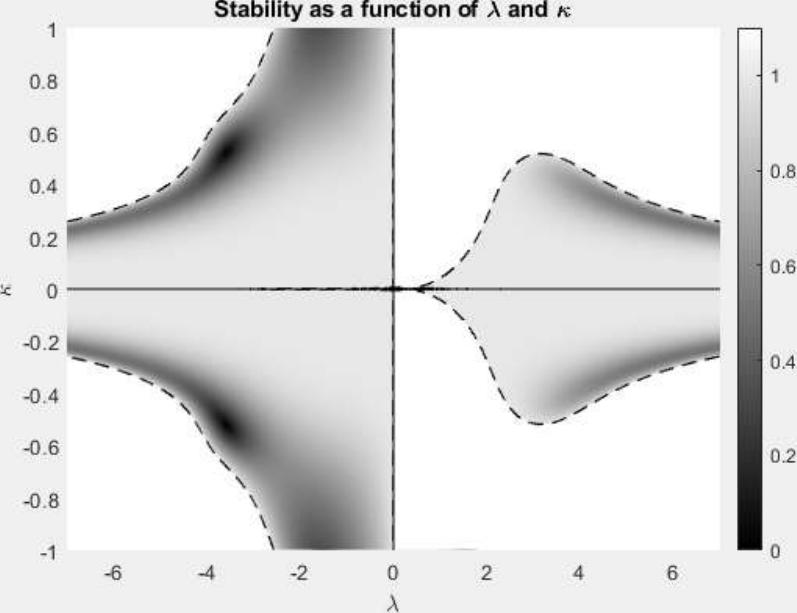}
\endminipage
\caption{RK4; left to right : stability regions for \textbf{Order 6} c, cs, f, zff and zffs.\label{rk46}}
\end{figure}

Fig. \ref{rk41} to \ref{rk46} show  the stability regions in the $\lambda$---$\kappa$ plane 
($\lambda$ and $\kappa$ respectively being the stability-number and the scaled wavenumber) for RK4 time integrator, 
for spatial orders $N=1$ to $6$, 
respectively (left to right) for centred (c), centred staggered (cs), forward (f), zigzag forward-first (zff) and zigzag forward-first staggered (zffs) schemes.
The greyed out areas 
represent the couples $(\lambda,\kappa)$ for which the scheme is stable 
(i.e. the amplification factor $|G| \leq 1$), while the dotted contours 
correspond to the critical case $|G| = 1$. A given scheme is conditionally
stable if and only if there exist a $\lambda_c \in \mathds{R}$ such that
the line of equation $\lambda = \lambda_c$ is included in the grey area.

\newpage
\section{Stability regions for RK5}
\begin{figure}[H]
\minipage{0.18\textwidth}
  \includegraphics[width=\linewidth]{blank.pdf}
\endminipage\hfill
\minipage{0.19\textwidth}
  \includegraphics[width=\linewidth]{blank.pdf}
\endminipage\hfill
\minipage{0.19\textwidth}
  \includegraphics[width=\linewidth]{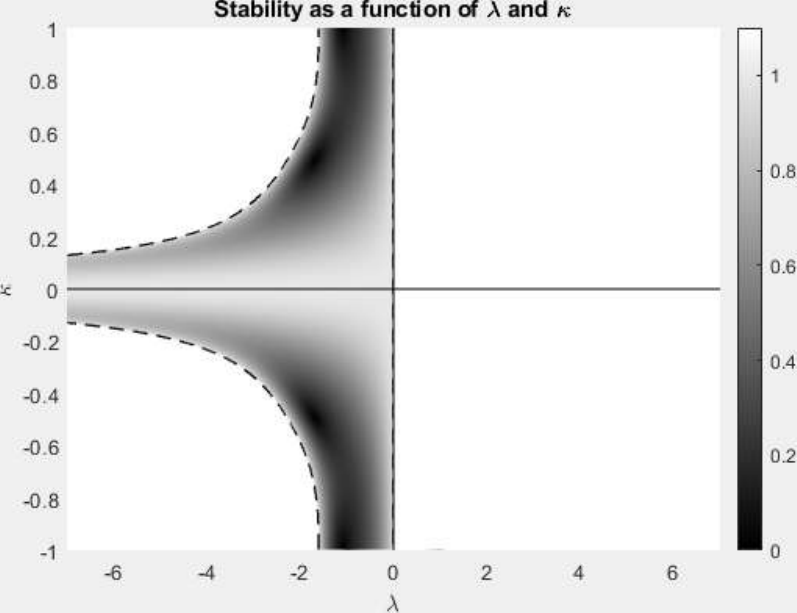}
\endminipage\hfill
\minipage{0.19\textwidth}
  \includegraphics[width=\linewidth]{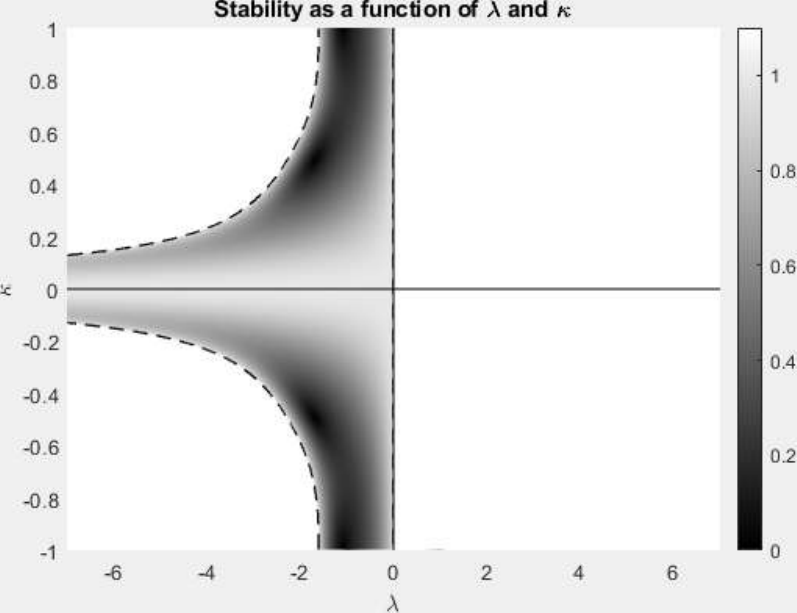}
\endminipage\hfill
\minipage{0.19\textwidth}%
  \includegraphics[width=\linewidth]{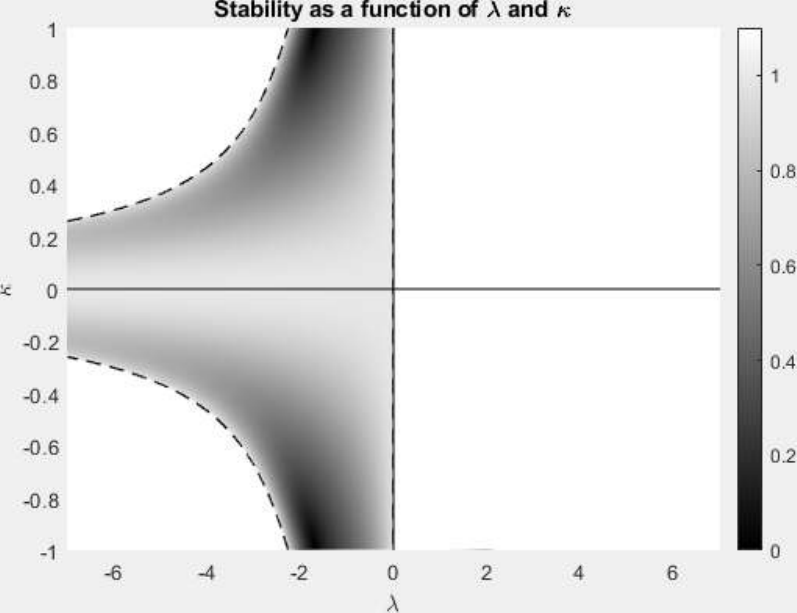}
\endminipage
\caption{RK5; left to right : stability regions for \textbf{Order 1} c, cs, f, zff and zffs.\label{rk51}}
\end{figure}

\begin{figure}[H]
\minipage{0.19\textwidth}
  \includegraphics[width=\linewidth]{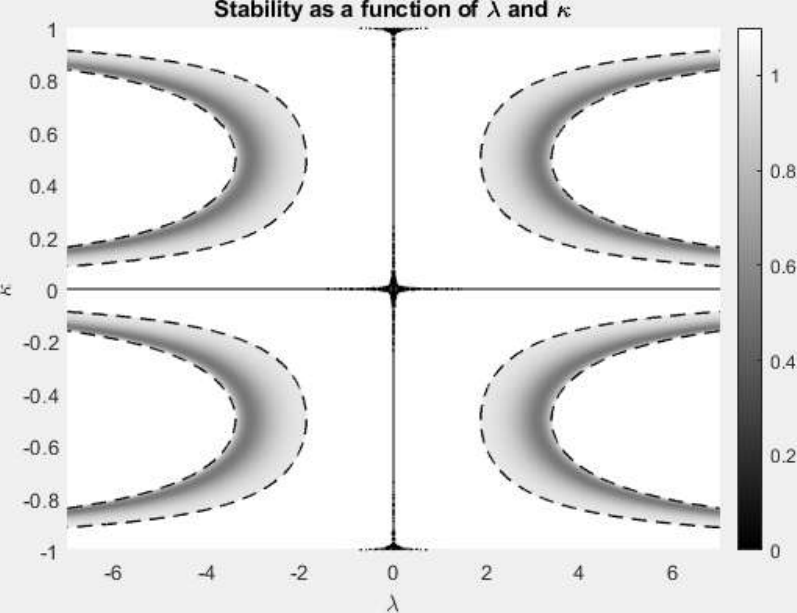}
\endminipage\hfill
\minipage{0.19\textwidth}
  \includegraphics[width=\linewidth]{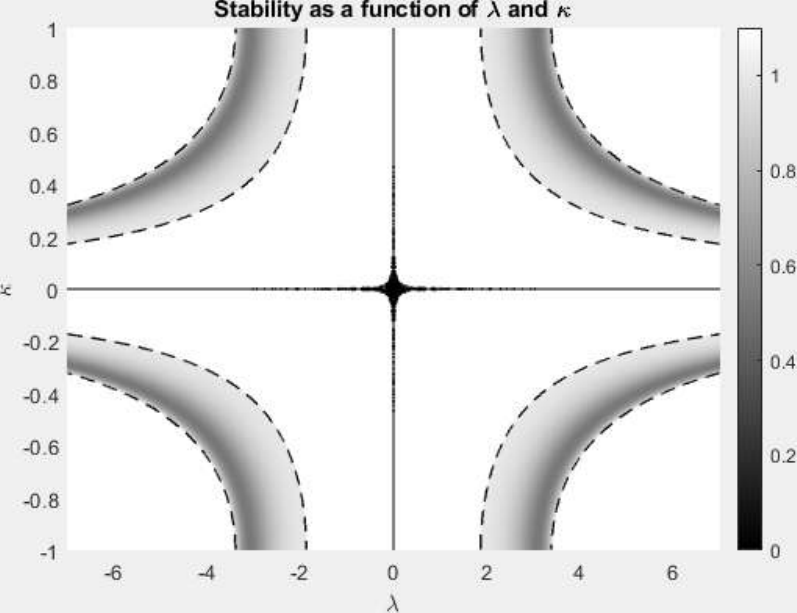}
\endminipage\hfill
\minipage{0.19\textwidth}
  \includegraphics[width=\linewidth]{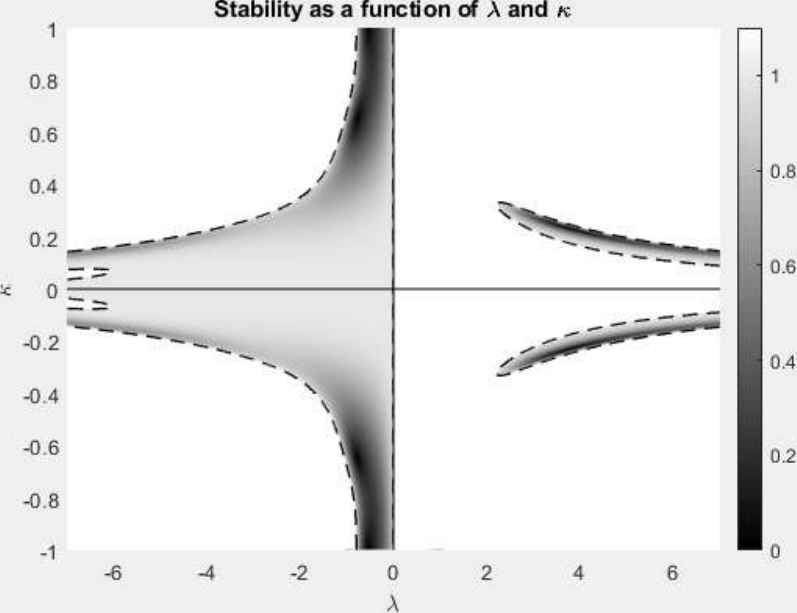}
\endminipage\hfill
\minipage{0.19\textwidth}
  \includegraphics[width=\linewidth]{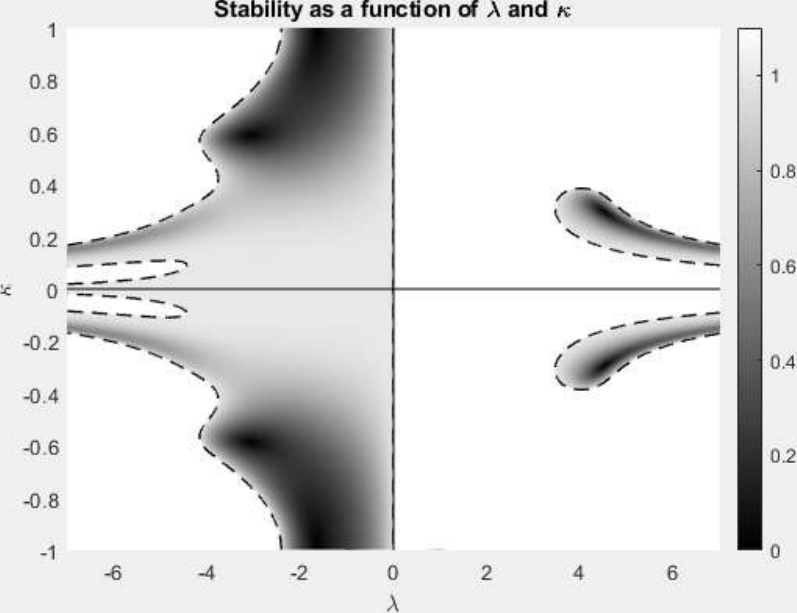}
\endminipage\hfill
\minipage{0.19\textwidth}%
  \includegraphics[width=\linewidth]{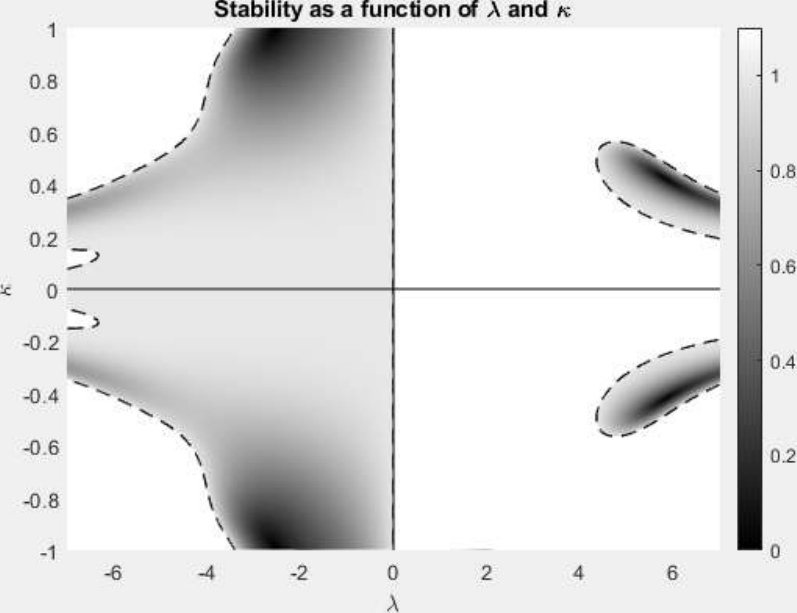}
\endminipage
\caption{RK5; left to right : stability regions for \textbf{Order 2} c, cs, f, zff and zffs.\label{rk52}}
\end{figure}

\begin{figure}[H]
\minipage{0.19\textwidth}
  \includegraphics[width=\linewidth]{blank.pdf}
\endminipage\hfill
\minipage{0.19\textwidth}
  \includegraphics[width=\linewidth]{blank.pdf}
\endminipage\hfill
\minipage{0.19\textwidth}
  \includegraphics[width=\linewidth]{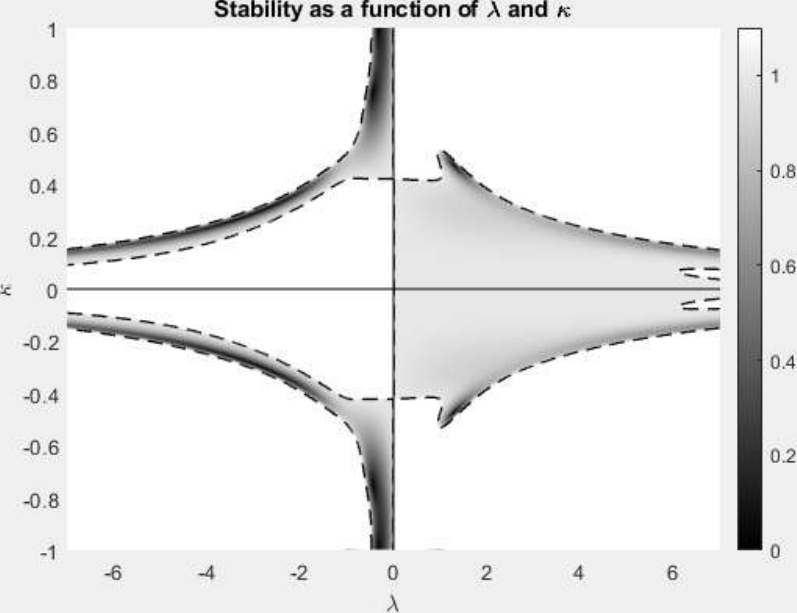}
\endminipage\hfill
\minipage{0.19\textwidth}
  \includegraphics[width=\linewidth]{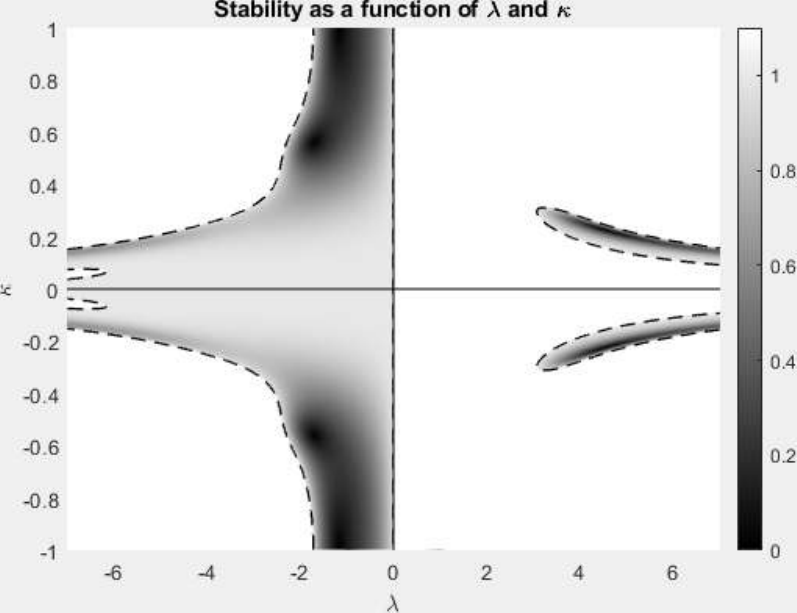}
\endminipage\hfill
\minipage{0.19\textwidth}%
  \includegraphics[width=\linewidth]{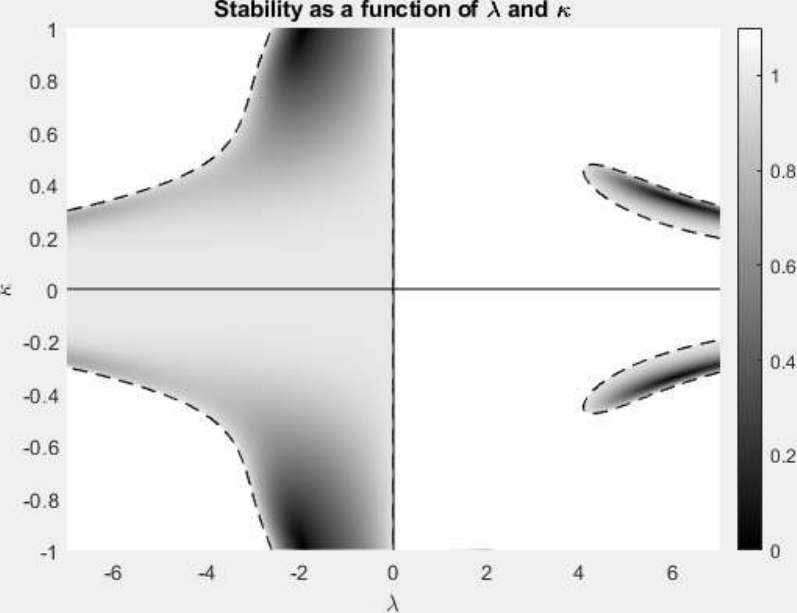}
\endminipage
\caption{RK5; left to right : stability regions for \textbf{Order 3} c, cs, f, zff and zffs.\label{rk53}}
\end{figure}

\begin{figure}[H]
\minipage{0.19\textwidth}
  \includegraphics[width=\linewidth]{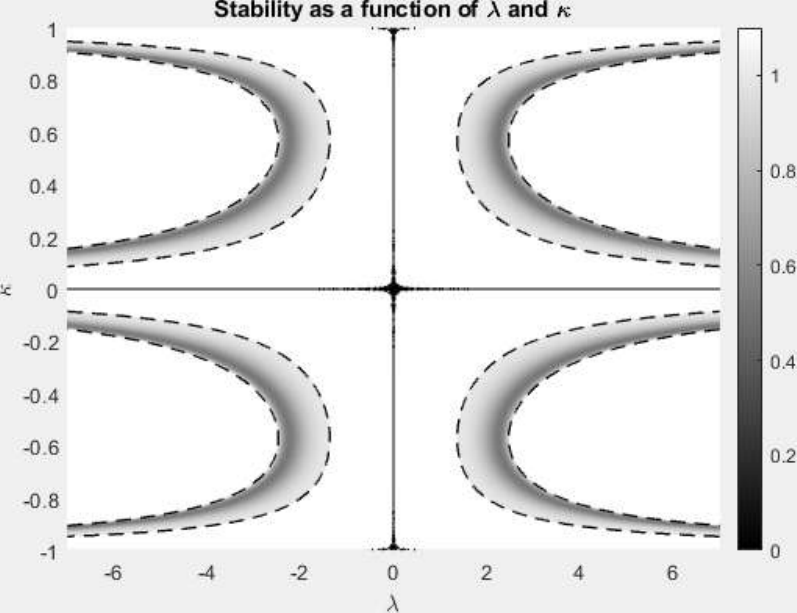}
\endminipage\hfill
\minipage{0.19\textwidth}
  \includegraphics[width=\linewidth]{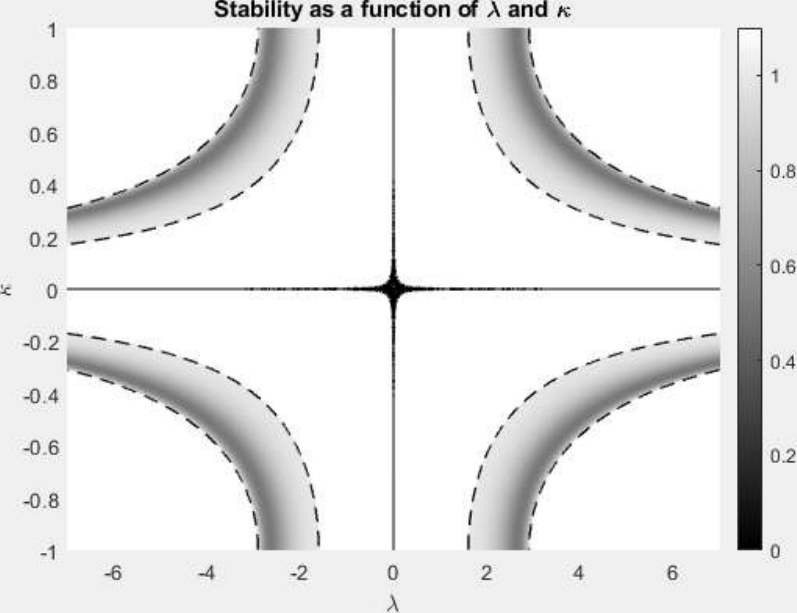}
\endminipage\hfill
\minipage{0.19\textwidth}
  \includegraphics[width=\linewidth]{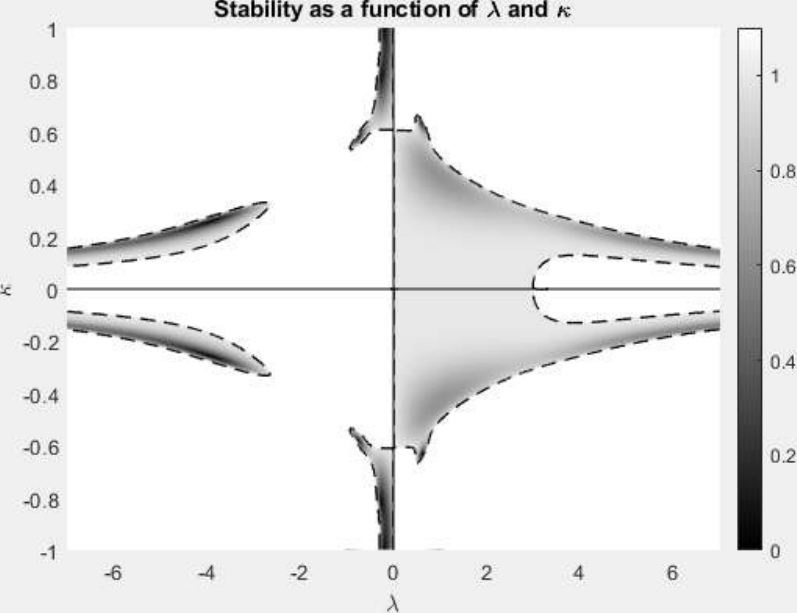}
\endminipage\hfill
\minipage{0.19\textwidth}
  \includegraphics[width=\linewidth]{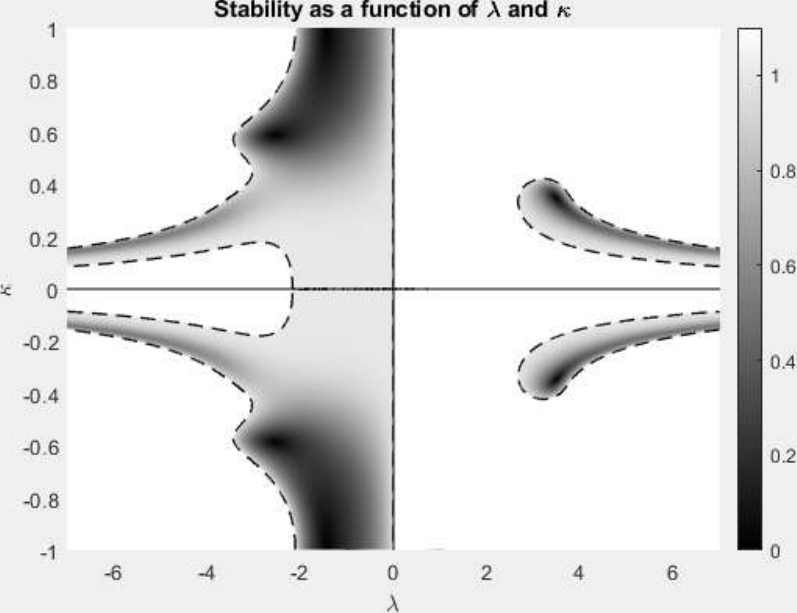}
\endminipage\hfill
\minipage{0.19\textwidth}%
  \includegraphics[width=\linewidth]{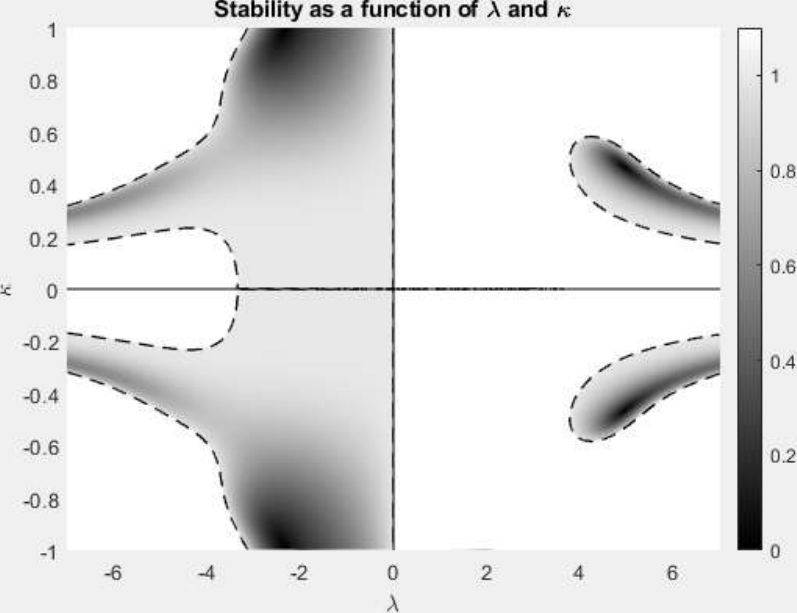}
\endminipage
\caption{RK5; left to right : stability regions for \textbf{Order 4} c, cs, f, zff and zffs.\label{rk54}}
\end{figure}

\begin{figure}[H]
\minipage{0.19\textwidth}
  \includegraphics[width=\linewidth]{blank.pdf}
\endminipage\hfill
\minipage{0.19\textwidth}
  \includegraphics[width=\linewidth]{blank.pdf}
\endminipage\hfill
\minipage{0.19\textwidth}
  \includegraphics[width=\linewidth]{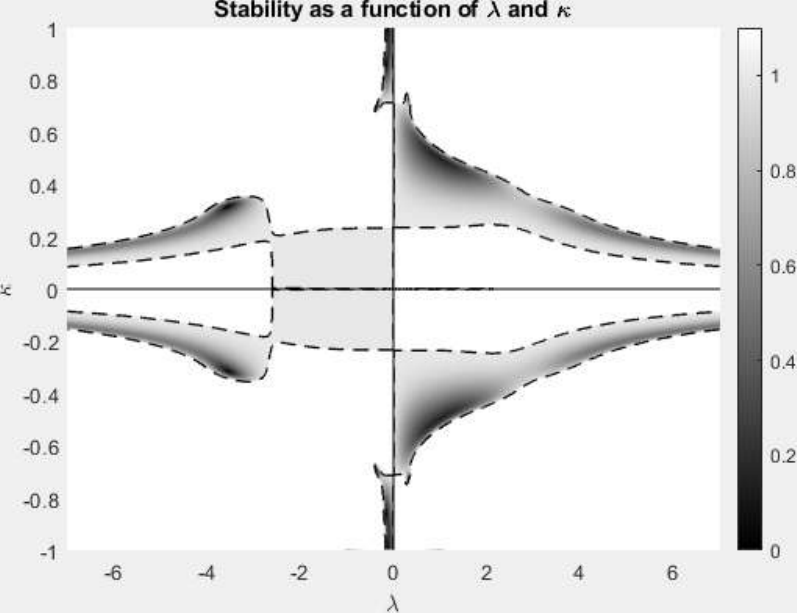}
\endminipage\hfill
\minipage{0.19\textwidth}
  \includegraphics[width=\linewidth]{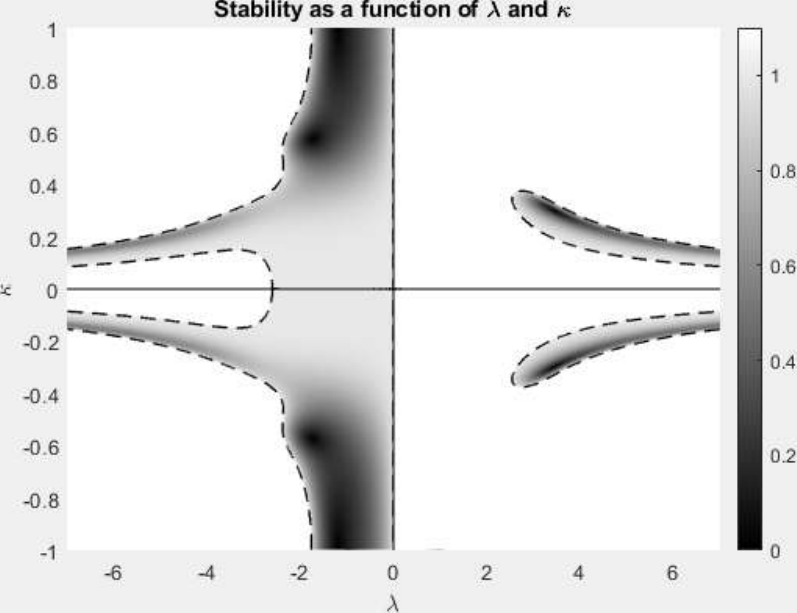}
\endminipage\hfill
\minipage{0.19\textwidth}%
  \includegraphics[width=\linewidth]{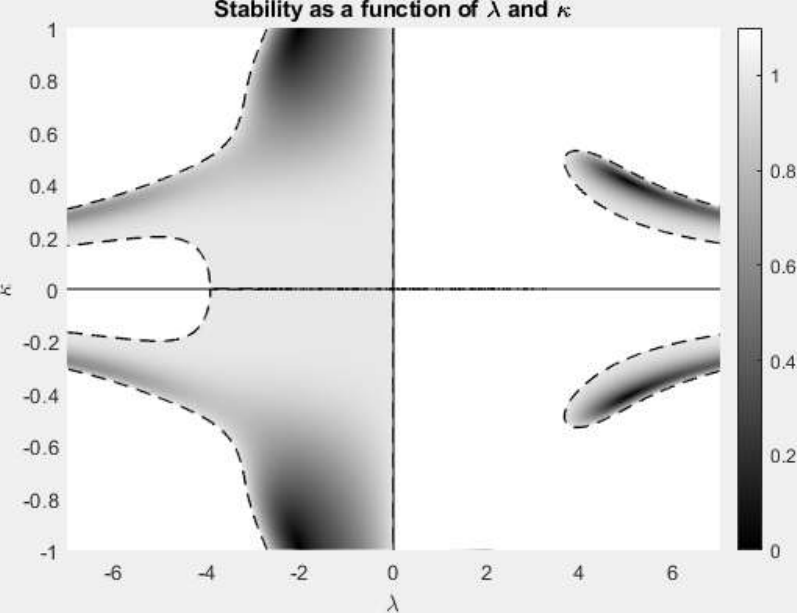}
\endminipage
\caption{RK5; left to right : stability regions for \textbf{Order 5} c, cs, f, zff and zffs.\label{rk55}}
\end{figure}

\begin{figure}[H]
\minipage{0.19\textwidth}
  \includegraphics[width=\linewidth]{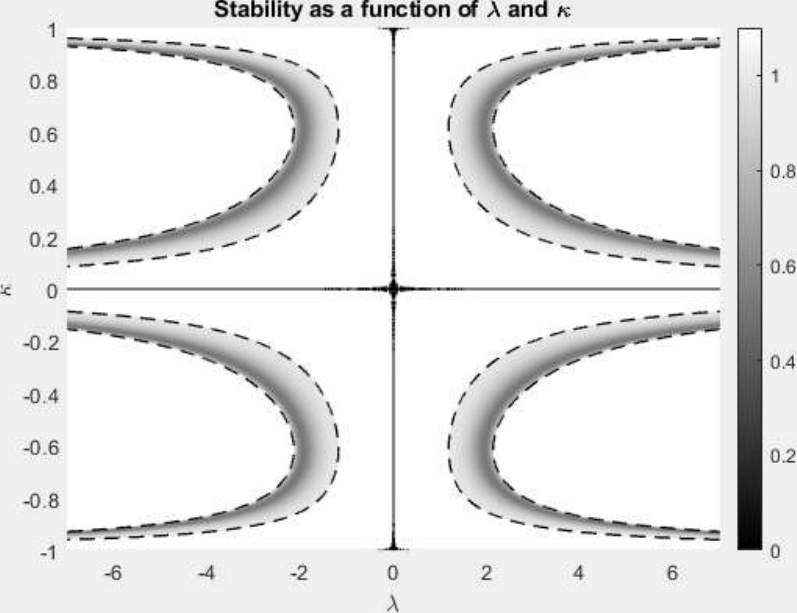}
\endminipage\hfill
\minipage{0.19\textwidth}
  \includegraphics[width=\linewidth]{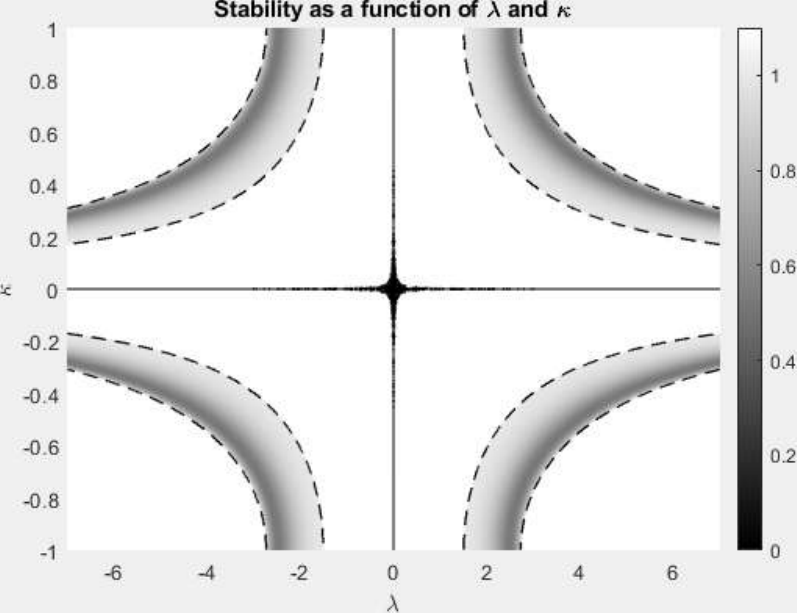}
\endminipage\hfill
\minipage{0.19\textwidth}
  \includegraphics[width=\linewidth]{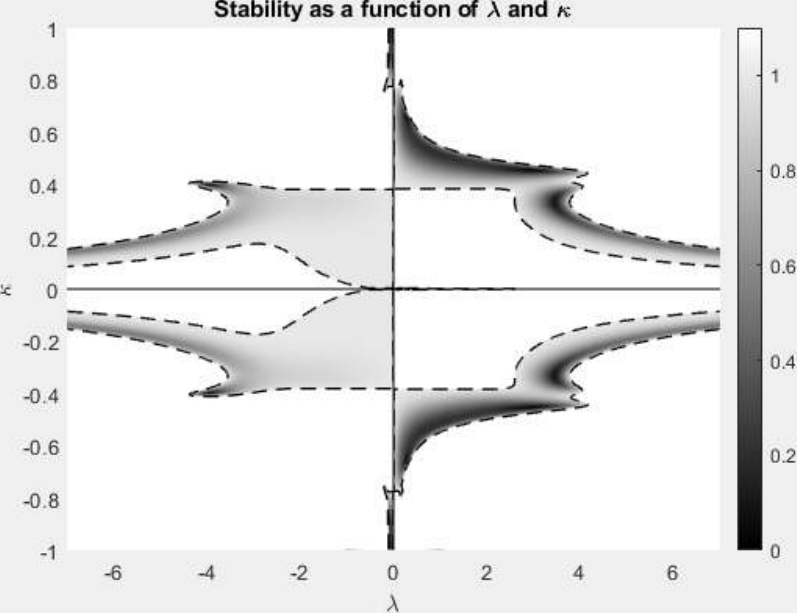}
\endminipage\hfill
\minipage{0.19\textwidth}
  \includegraphics[width=\linewidth]{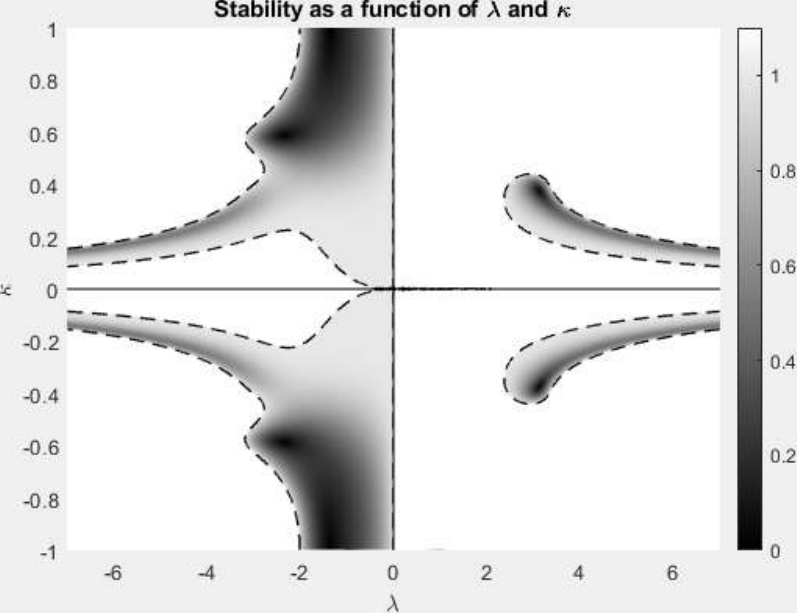}
\endminipage\hfill
\minipage{0.19\textwidth}%
  \includegraphics[width=\linewidth]{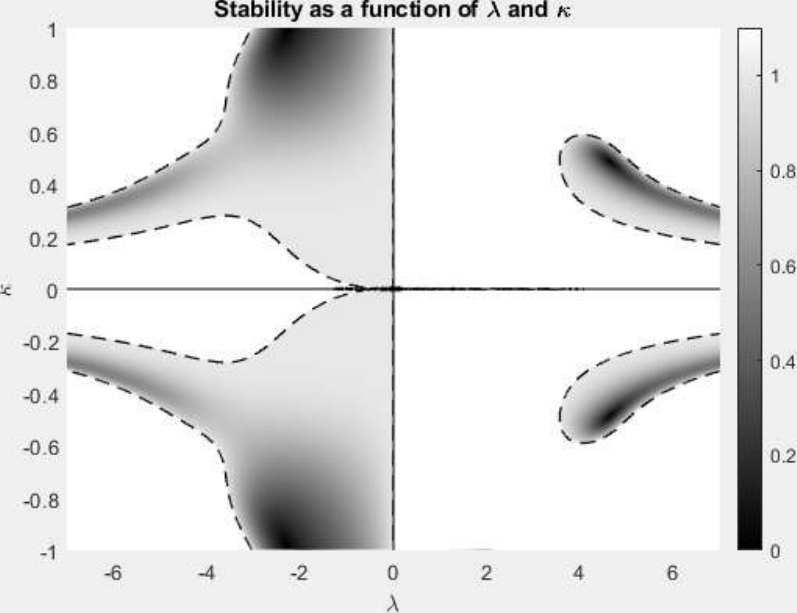}
\endminipage
\caption{RK5; left to right : stability regions for \textbf{Order 6} c, cs, f, zff and zffs.\label{rk56}}
\end{figure}

Fig. \ref{rk51} to \ref{rk56} show  the stability regions in the $\lambda$---$\kappa$ plane 
($\lambda$ and $\kappa$ respectively being the stability-number and the scaled wavenumber) for RK5 time integrator, 
for spatial orders $N=1$ to $6$, 
respectively (left to right) for centred (c), centred staggered (cs), forward (f), zigzag forward-first (zff) and zigzag forward-first staggered (zffs) schemes.
The greyed out areas 
represent the couples $(\lambda,\kappa)$ for which the scheme is stable 
(i.e. the amplification factor $|G| \leq 1$), while the dotted contours 
correspond to the critical case $|G| = 1$. A given scheme is conditionally
stable if and only if there exist a $\lambda_c \in \mathds{R}$ such that
the line of equation $\lambda = \lambda_c$ is included in the grey area.

\newpage
\section{Stability regions for RK6}
\begin{figure}[H]
\minipage{0.18\textwidth}
  \includegraphics[width=\linewidth]{blank.pdf}
\endminipage\hfill
\minipage{0.19\textwidth}
  \includegraphics[width=\linewidth]{blank.pdf}
\endminipage\hfill
\minipage{0.19\textwidth}
  \includegraphics[width=\linewidth]{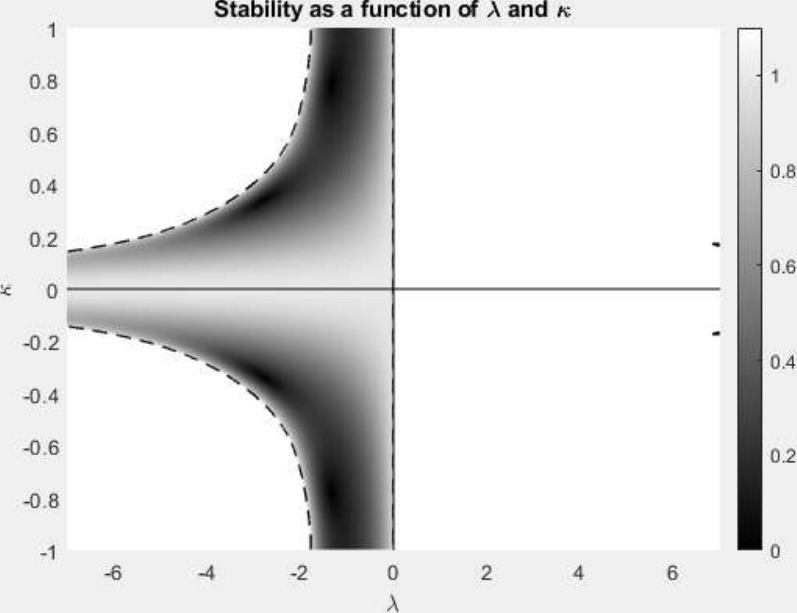}
\endminipage\hfill
\minipage{0.19\textwidth}
  \includegraphics[width=\linewidth]{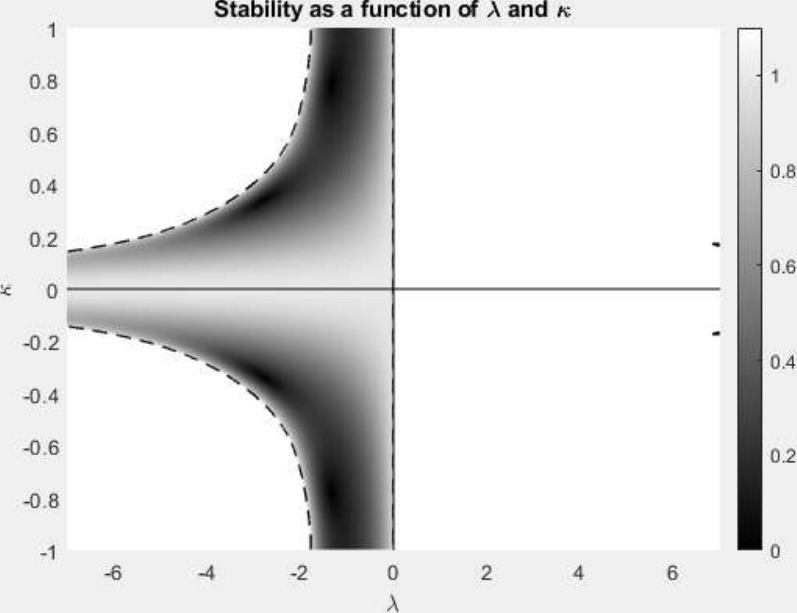}
\endminipage\hfill
\minipage{0.19\textwidth}%
  \includegraphics[width=\linewidth]{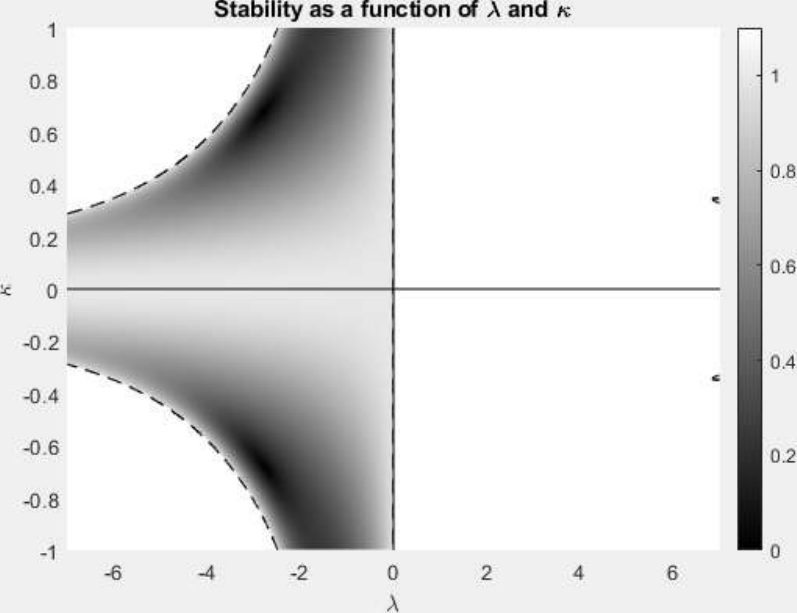}
\endminipage
\caption{RK6; left to right : stability regions for \textbf{Order 1} c, cs, f, zff and zffs.\label{rk61}}
\end{figure}

\begin{figure}[H]
\minipage{0.19\textwidth}
  \includegraphics[width=\linewidth]{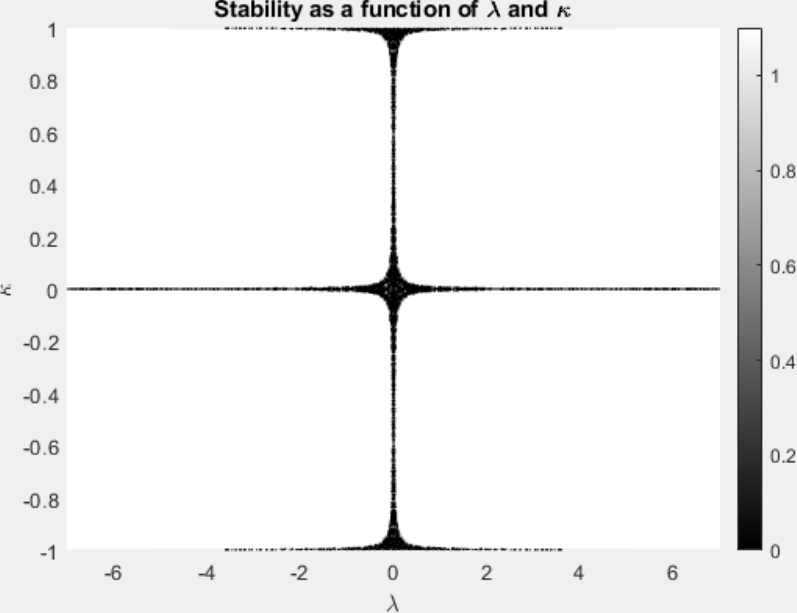}
\endminipage\hfill
\minipage{0.19\textwidth}
  \includegraphics[width=\linewidth]{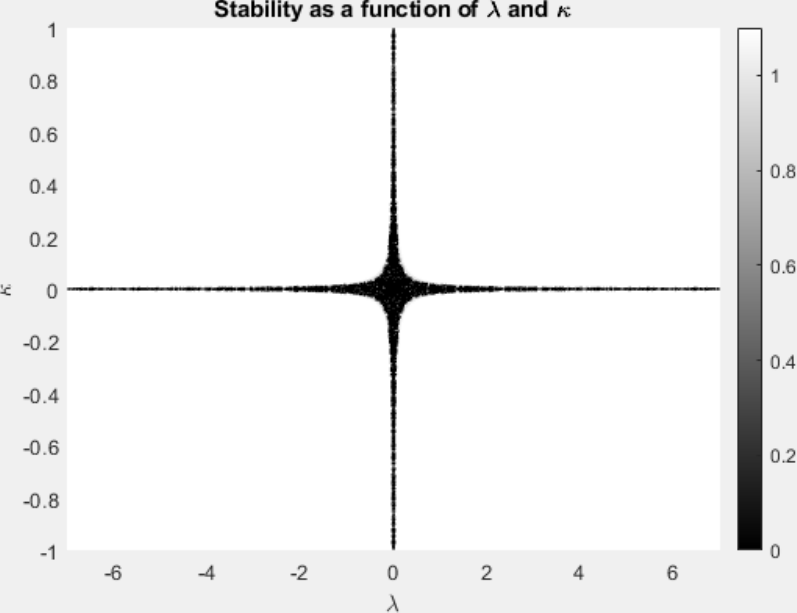}
\endminipage\hfill
\minipage{0.19\textwidth}
  \includegraphics[width=\linewidth]{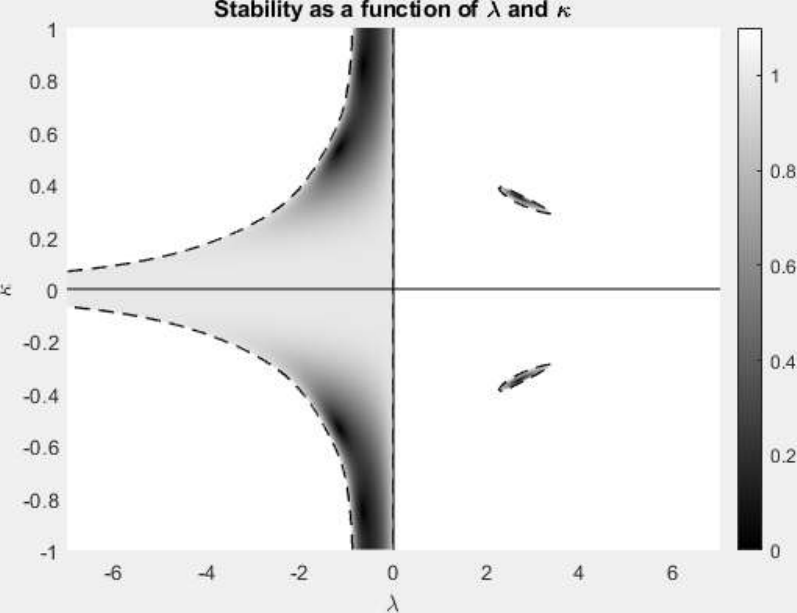}
\endminipage\hfill
\minipage{0.19\textwidth}
  \includegraphics[width=\linewidth]{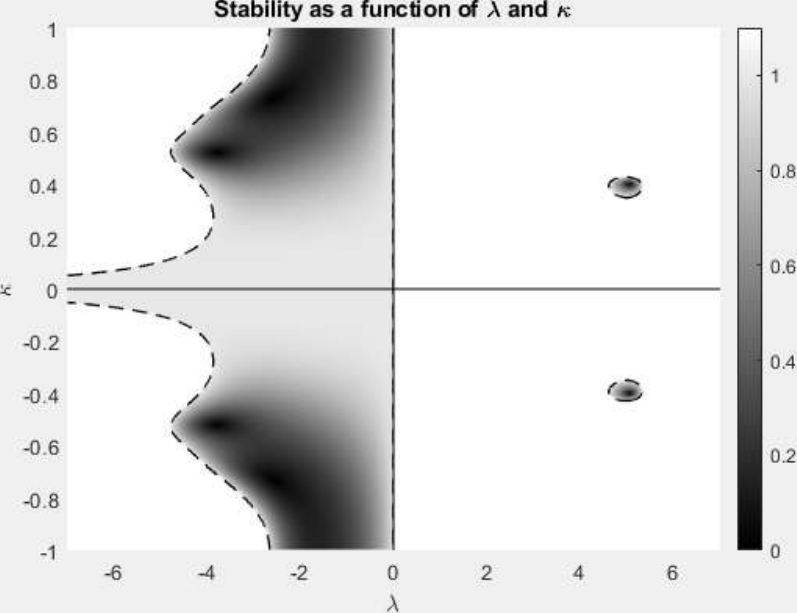}
\endminipage\hfill
\minipage{0.19\textwidth}%
  \includegraphics[width=\linewidth]{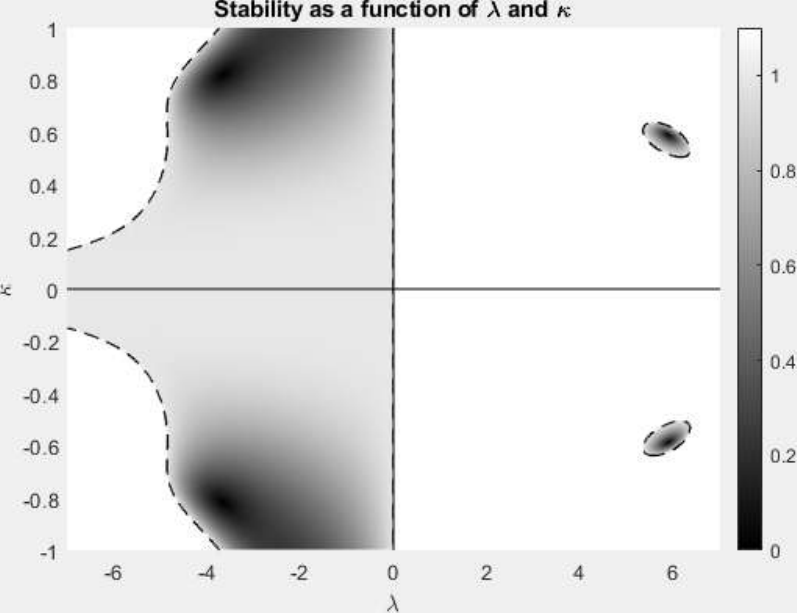}
\endminipage
\caption{RK6; left to right : stability regions for \textbf{Order 2} c, cs, f, zff and zffs.\label{rk62}}
\end{figure}

\begin{figure}[H]
\minipage{0.19\textwidth}
  \includegraphics[width=\linewidth]{blank.pdf}
\endminipage\hfill
\minipage{0.19\textwidth}
  \includegraphics[width=\linewidth]{blank.pdf}
\endminipage\hfill
\minipage{0.19\textwidth}
  \includegraphics[width=\linewidth]{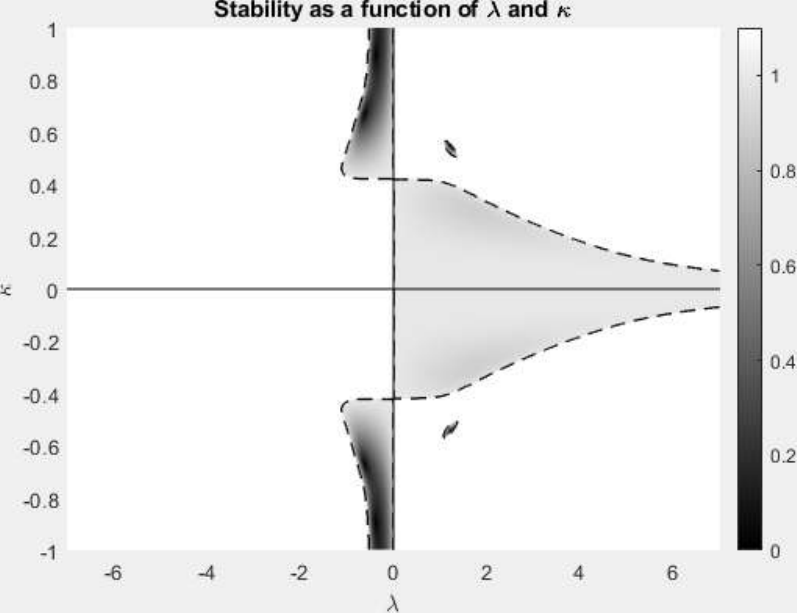}
\endminipage\hfill
\minipage{0.19\textwidth}
  \includegraphics[width=\linewidth]{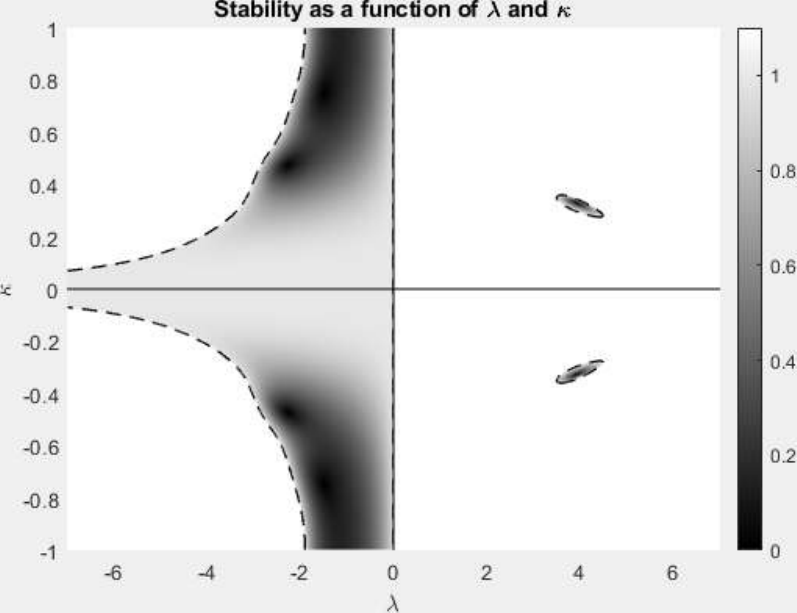}
\endminipage\hfill
\minipage{0.19\textwidth}%
  \includegraphics[width=\linewidth]{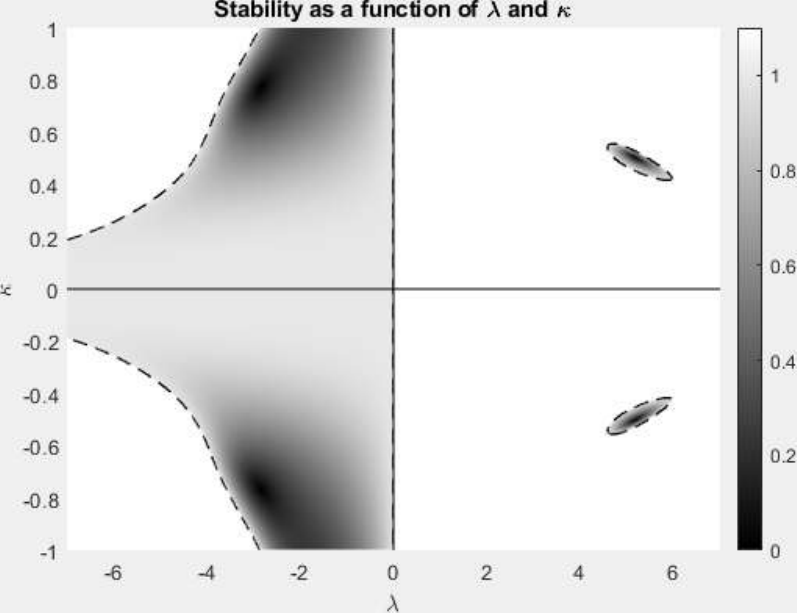}
\endminipage
\caption{RK6; left to right : stability regions for \textbf{Order 3} c, cs, f, zff and zffs.\label{rk63}}
\end{figure}

\begin{figure}[H]
\minipage{0.19\textwidth}
  \includegraphics[width=\linewidth]{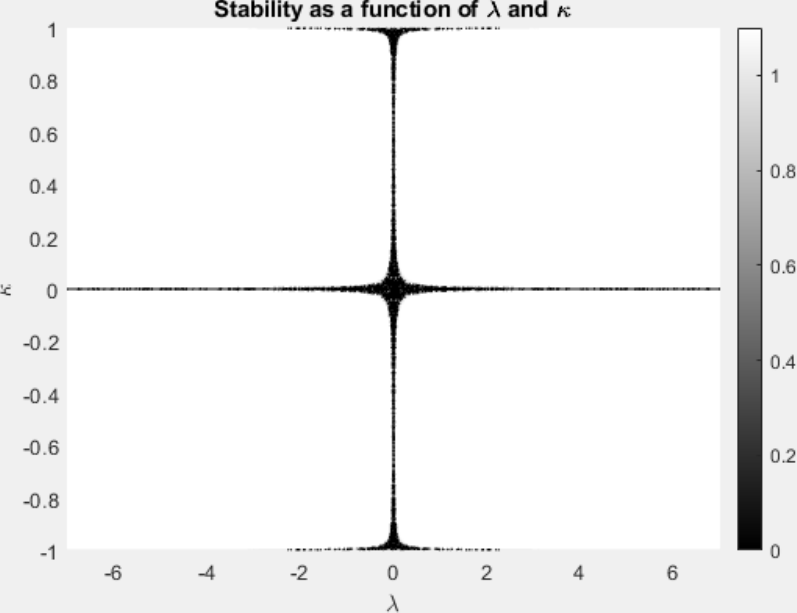}
\endminipage\hfill
\minipage{0.19\textwidth}
  \includegraphics[width=\linewidth]{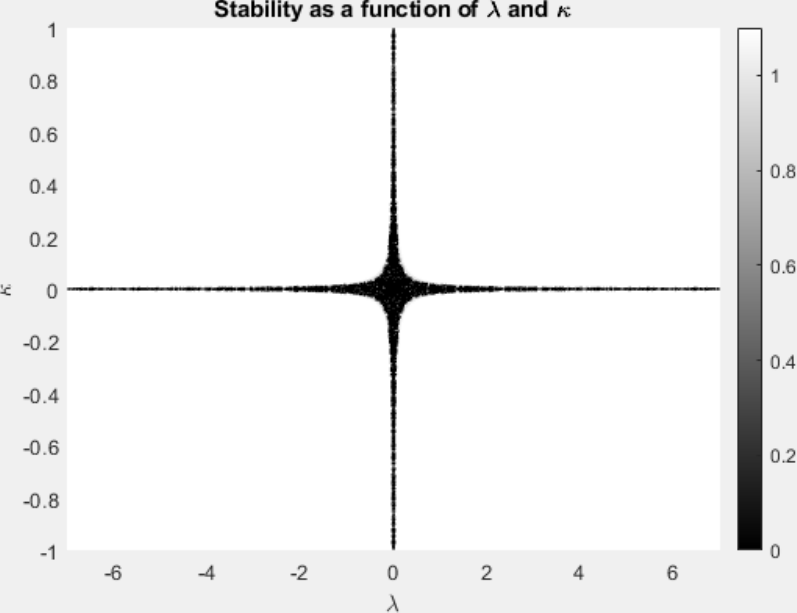}
\endminipage\hfill
\minipage{0.19\textwidth}
  \includegraphics[width=\linewidth]{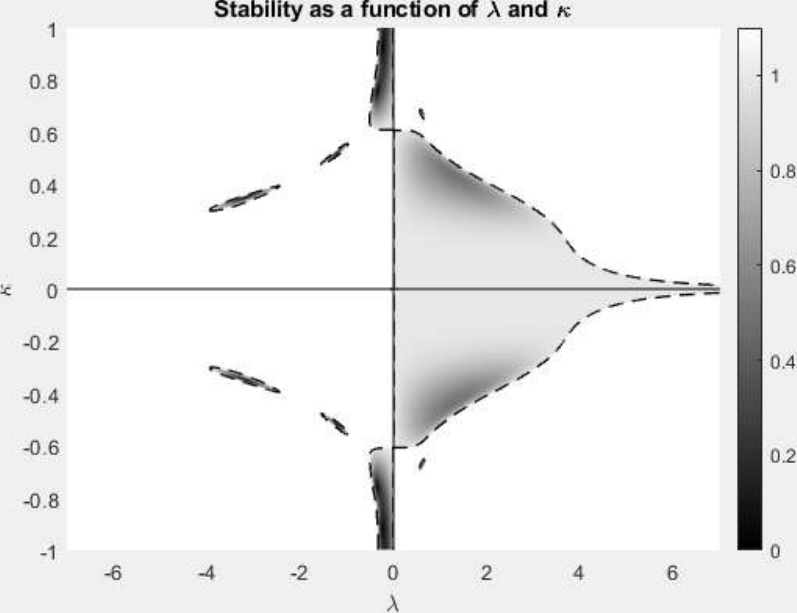}
\endminipage\hfill
\minipage{0.19\textwidth}
  \includegraphics[width=\linewidth]{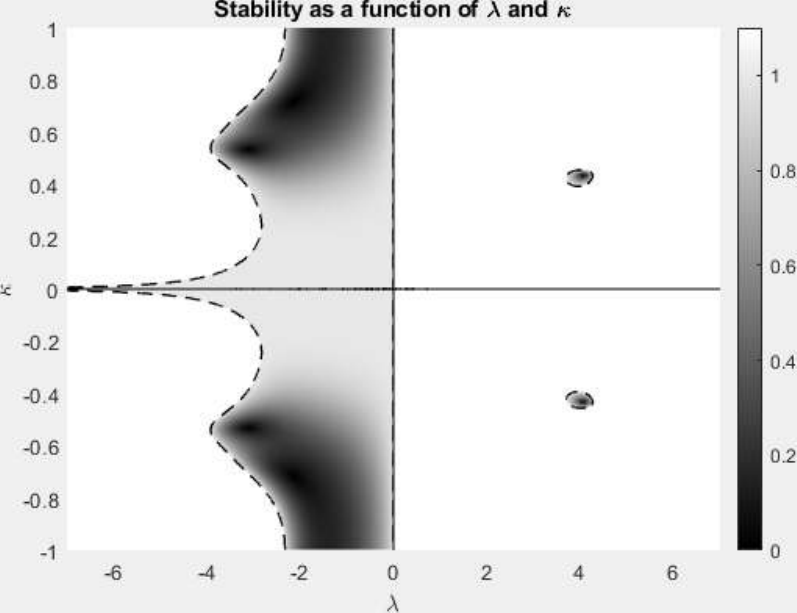}
\endminipage\hfill
\minipage{0.19\textwidth}%
  \includegraphics[width=\linewidth]{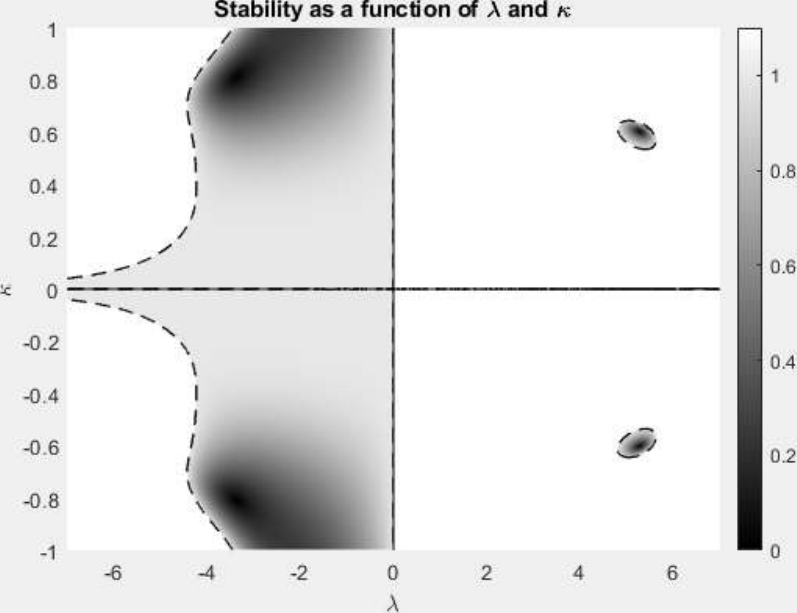}
\endminipage
\caption{RK6; left to right : stability regions for \textbf{Order 4} c, cs, f, zff and zffs.\label{rk64}}
\end{figure}

\begin{figure}[H]
\minipage{0.19\textwidth}
  \includegraphics[width=\linewidth]{blank.pdf}
\endminipage\hfill
\minipage{0.19\textwidth}
  \includegraphics[width=\linewidth]{blank.pdf}
\endminipage\hfill
\minipage{0.19\textwidth}
  \includegraphics[width=\linewidth]{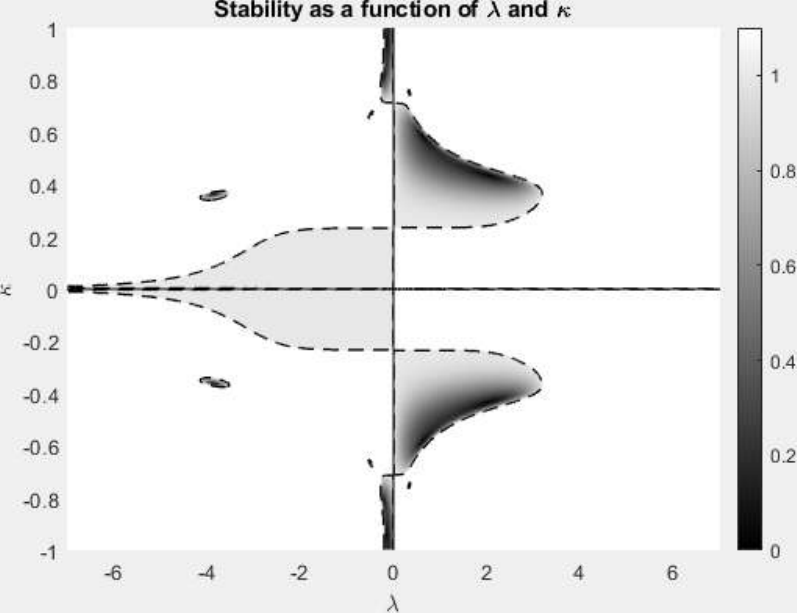}
\endminipage\hfill
\minipage{0.19\textwidth}
  \includegraphics[width=\linewidth]{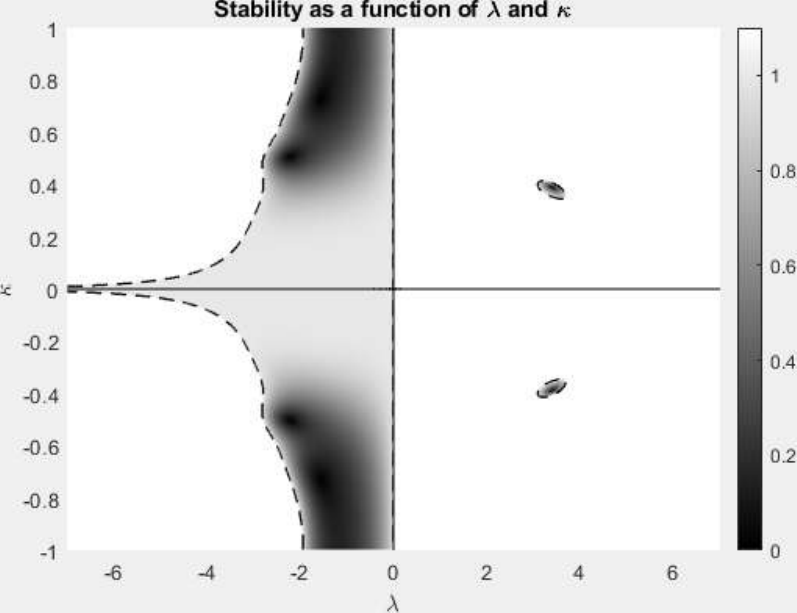}
\endminipage\hfill
\minipage{0.19\textwidth}%
  \includegraphics[width=\linewidth]{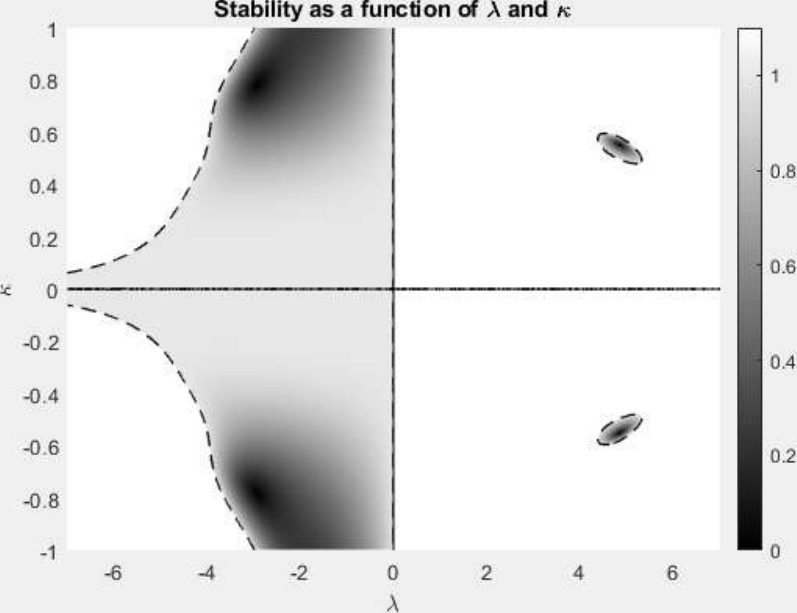}
\endminipage
\caption{RK6; left to right : stability regions for \textbf{Order 5} c, cs, f, zff and zffs.\label{rk65}}
\end{figure}

\begin{figure}[H]
\minipage{0.19\textwidth}
  \includegraphics[width=\linewidth]{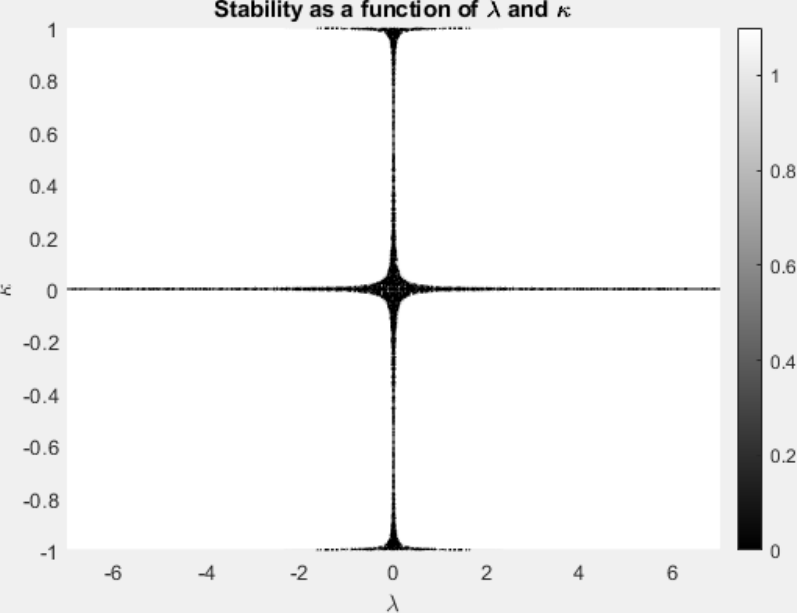}
\endminipage\hfill
\minipage{0.19\textwidth}
  \includegraphics[width=\linewidth]{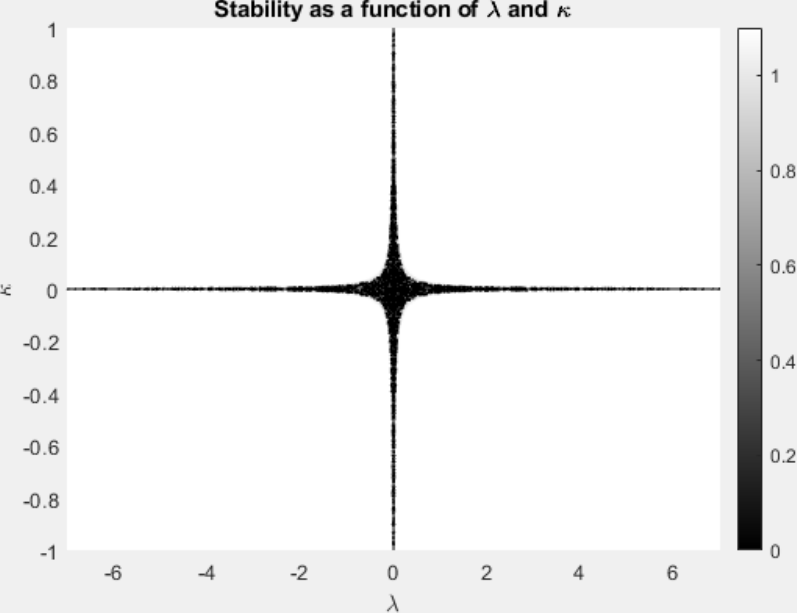}
\endminipage\hfill
\minipage{0.19\textwidth}
  \includegraphics[width=\linewidth]{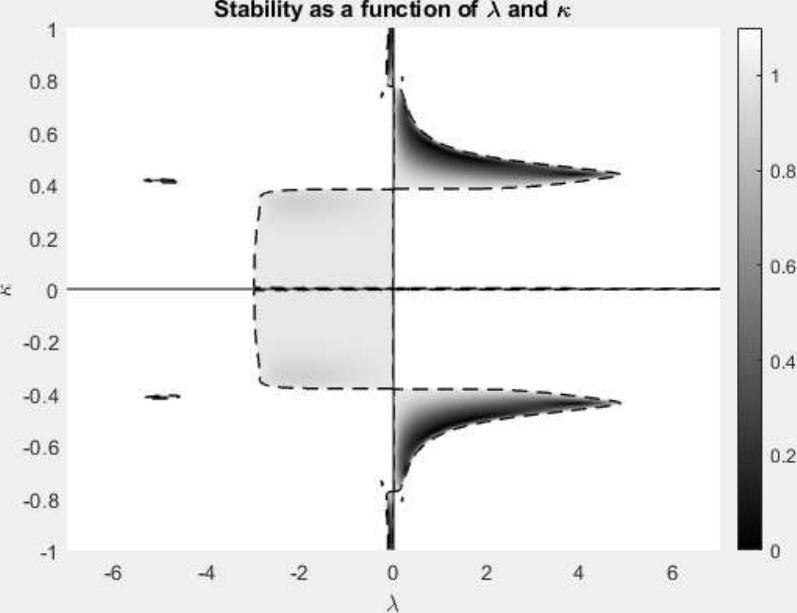}
\endminipage\hfill
\minipage{0.19\textwidth}
  \includegraphics[width=\linewidth]{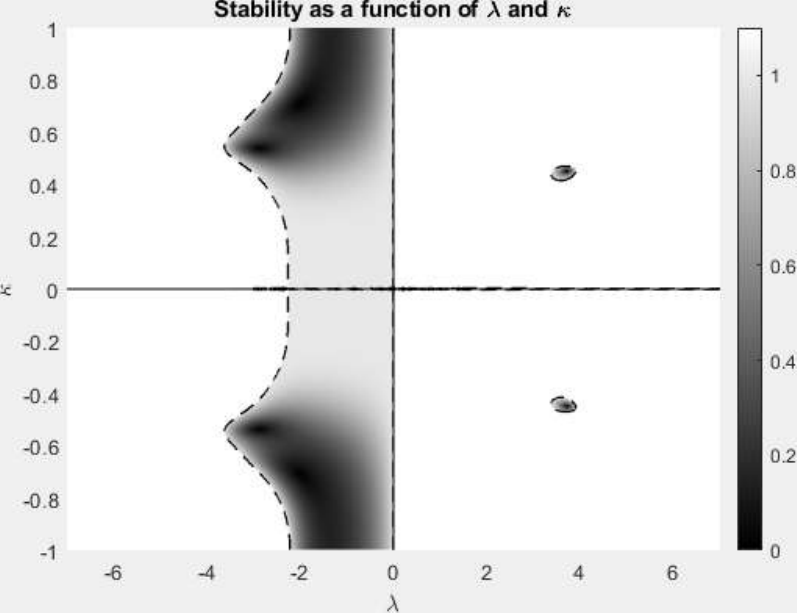}
\endminipage\hfill
\minipage{0.19\textwidth}%
  \includegraphics[width=\linewidth]{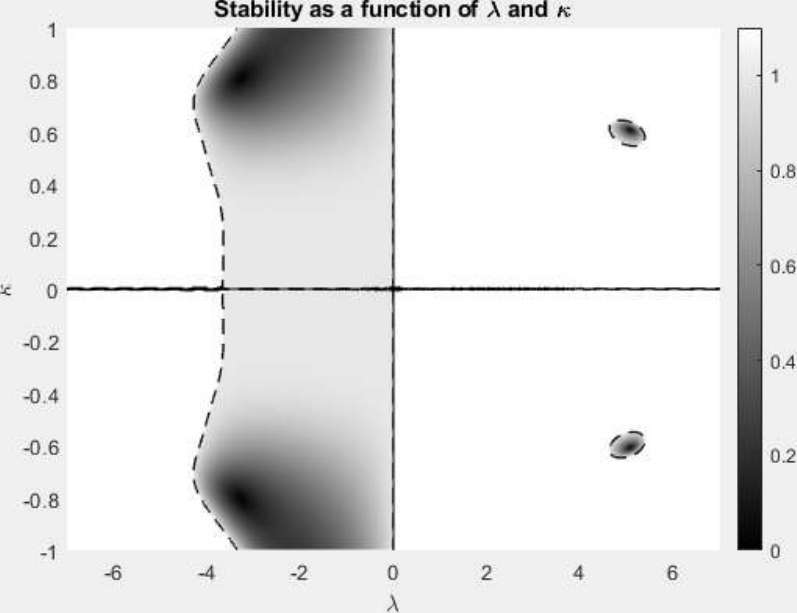}
\endminipage
\caption{RK6; left to right : stability regions for \textbf{Order 6} c, cs, f, zff and zffs.\label{rk66}}
\end{figure}

Fig. \ref{rk61} to \ref{rk66} show  the stability regions in the $\lambda$---$\kappa$ plane 
($\lambda$ and $\kappa$ respectively being the stability-number and the scaled wavenumber) for RK6 time integrator, 
for spatial orders $N=1$ to $6$, 
respectively (left to right) for centred (c), centred staggered (cs), forward (f), zigzag forward-first (zff) and zigzag forward-first staggered (zffs) schemes.
The greyed out areas 
represent the couples $(\lambda,\kappa)$ for which the scheme is stable 
(i.e. the amplification factor $|G| \leq 1$), while the dotted contours 
correspond to the critical case $|G| = 1$. A given scheme is conditionally
stable if and only if there exist a $\lambda_c \in \mathds{R}$ such that
the line of equation $\lambda = \lambda_c$ is included in the grey area.

\newpage
\section{Stability regions for RK7} \label{RK7Annex}
\begin{figure}[H]
\minipage{0.18\textwidth}
  \includegraphics[width=\linewidth]{blank.pdf}
\endminipage\hfill
\minipage{0.19\textwidth}
  \includegraphics[width=\linewidth]{blank.pdf}
\endminipage\hfill
\minipage{0.19\textwidth}
  \includegraphics[width=\linewidth]{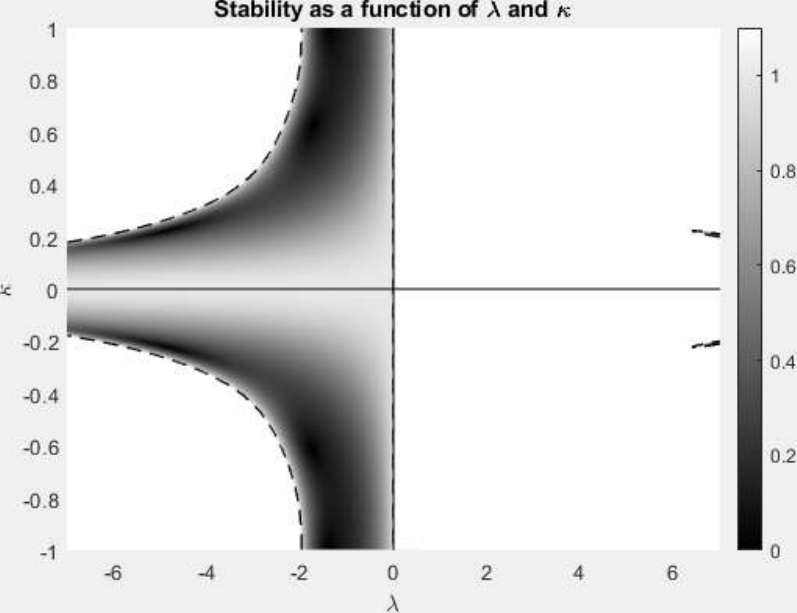}
\endminipage\hfill
\minipage{0.19\textwidth}
  \includegraphics[width=\linewidth]{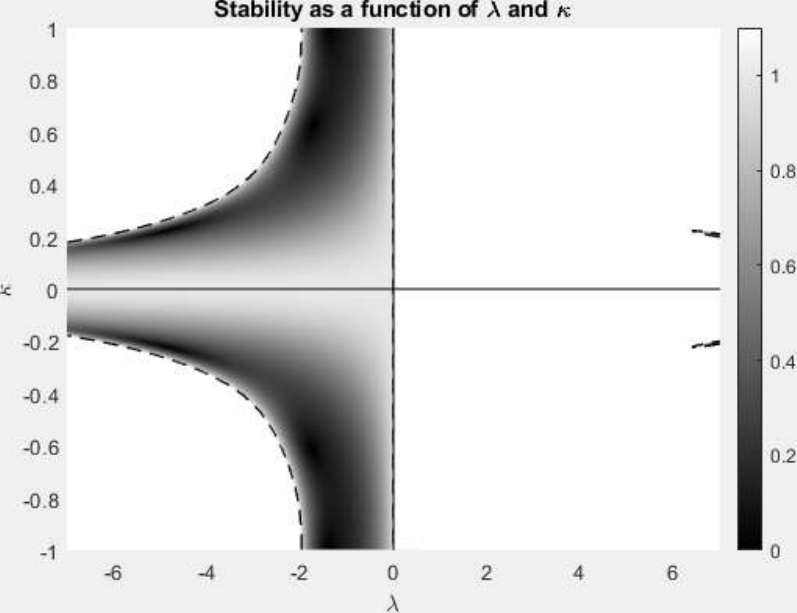}
\endminipage\hfill
\minipage{0.19\textwidth}%
  \includegraphics[width=\linewidth]{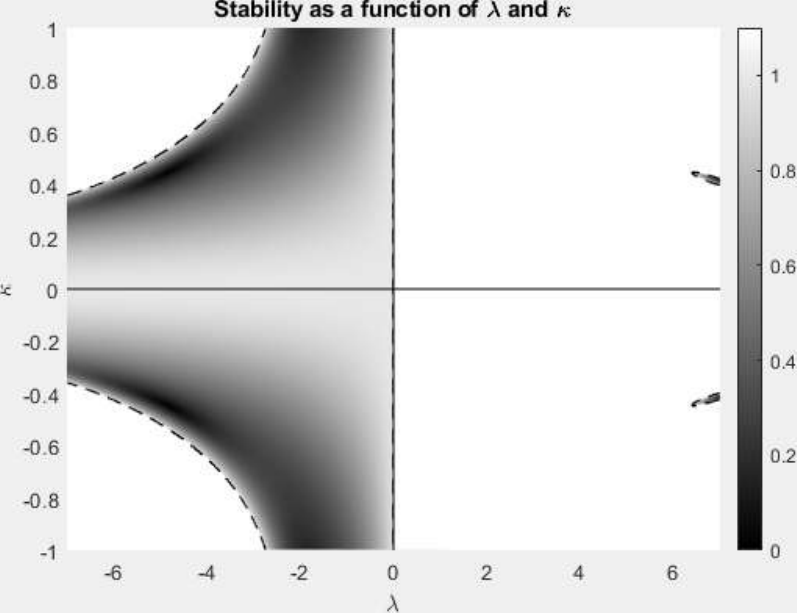}
\endminipage
\caption{RK7; left to right : stability regions for \textbf{Order 1} c, cs, f, zff and zffs.\label{rk71}}
\end{figure}

\begin{figure}[H]
\minipage{0.19\textwidth}
  \includegraphics[width=\linewidth]{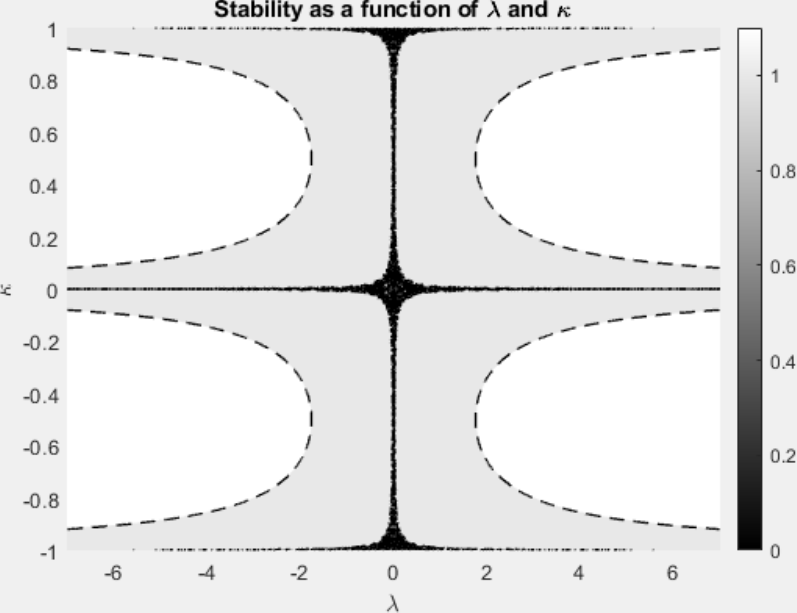}
\endminipage\hfill
\minipage{0.19\textwidth}
  \includegraphics[width=\linewidth]{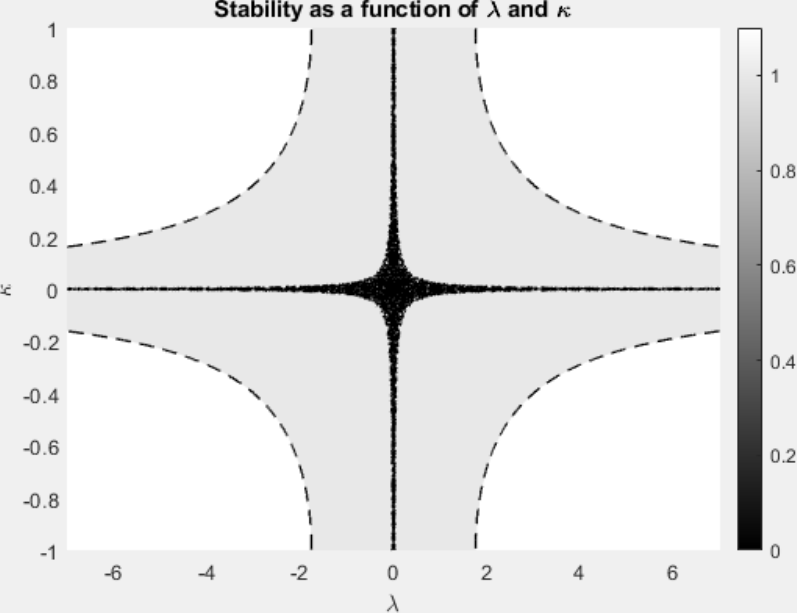}
\endminipage\hfill
\minipage{0.19\textwidth}
  \includegraphics[width=\linewidth]{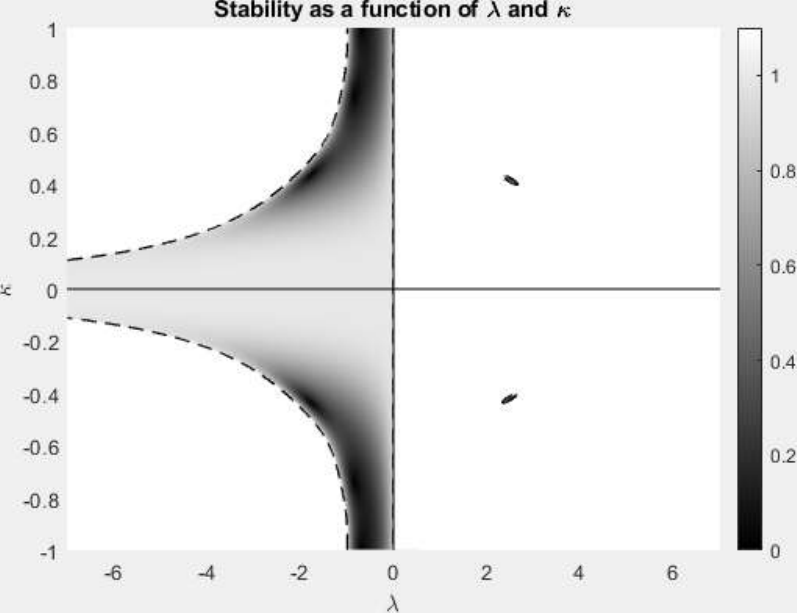}
\endminipage\hfill
\minipage{0.19\textwidth}
  \includegraphics[width=\linewidth]{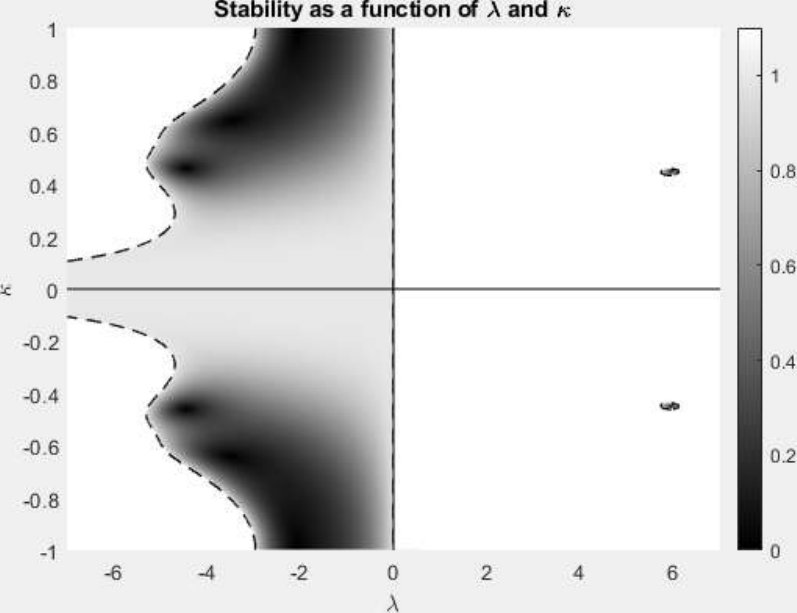}
\endminipage\hfill
\minipage{0.19\textwidth}%
  \includegraphics[width=\linewidth]{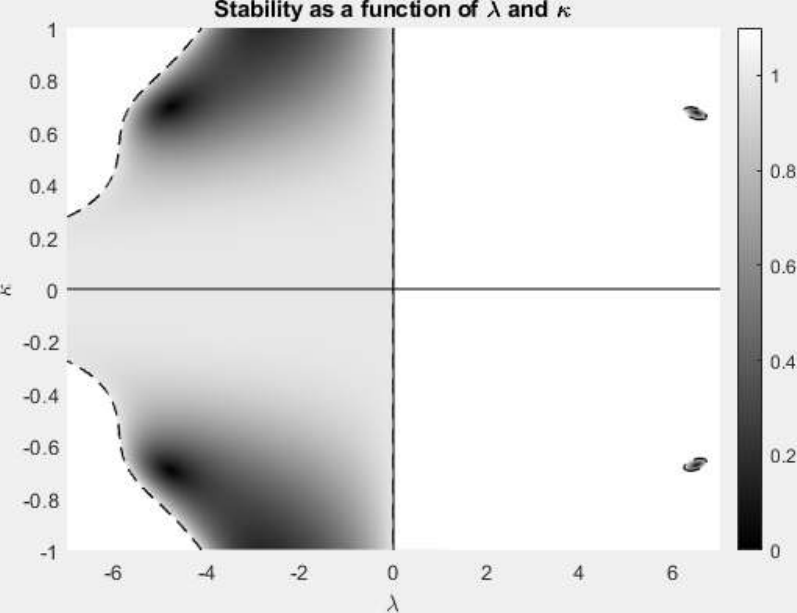}
\endminipage
\caption{RK7; left to right : stability regions for \textbf{Order 2} c, cs, f, zff and zffs.\label{rk72}}
\end{figure}

\begin{figure}[H]
\minipage{0.19\textwidth}
  \includegraphics[width=\linewidth]{blank.pdf}
\endminipage\hfill
\minipage{0.19\textwidth}
  \includegraphics[width=\linewidth]{blank.pdf}
\endminipage\hfill
\minipage{0.19\textwidth}
  \includegraphics[width=\linewidth]{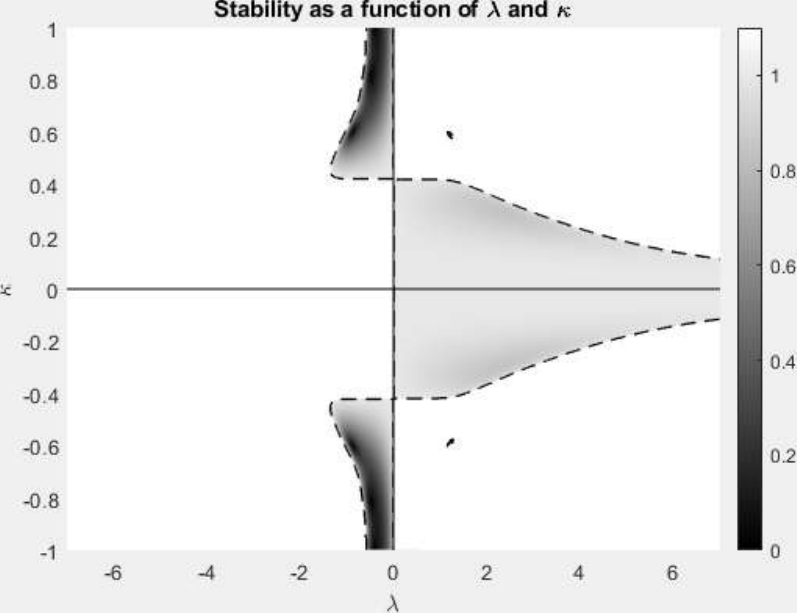}
\endminipage\hfill
\minipage{0.19\textwidth}
  \includegraphics[width=\linewidth]{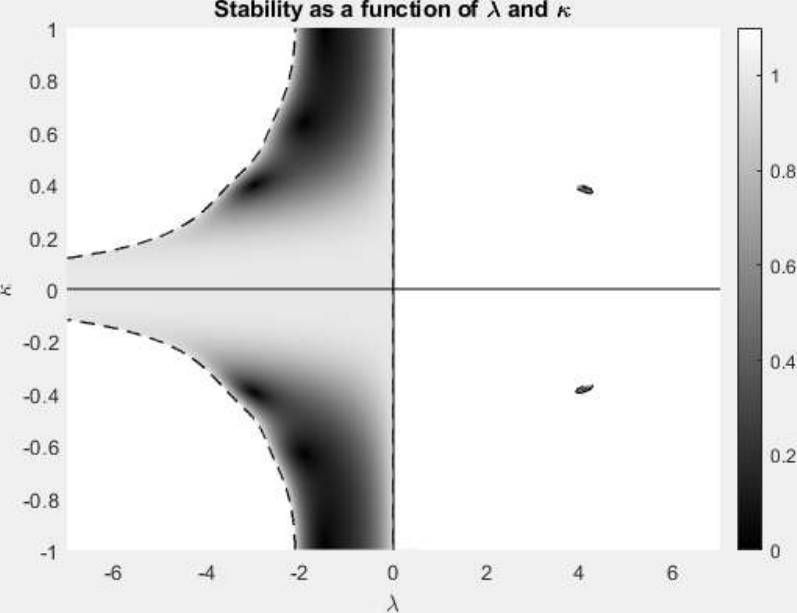}
\endminipage\hfill
\minipage{0.19\textwidth}%
  \includegraphics[width=\linewidth]{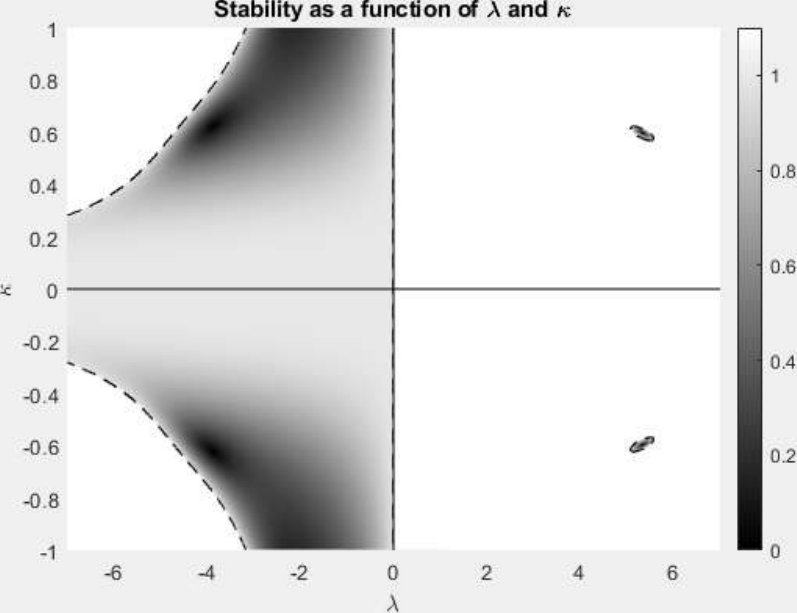}
\endminipage
\caption{RK7; left to right : stability regions for \textbf{Order 3} c, cs, f, zff and zffs.\label{rk73}}
\end{figure}

\begin{figure}[H]
\minipage{0.19\textwidth}
  \includegraphics[width=\linewidth]{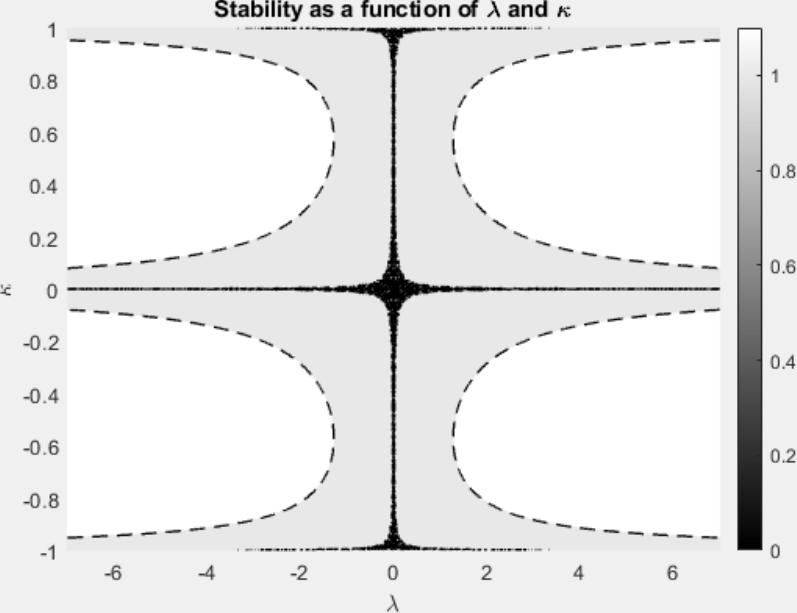}
\endminipage\hfill
\minipage{0.19\textwidth}
  \includegraphics[width=\linewidth]{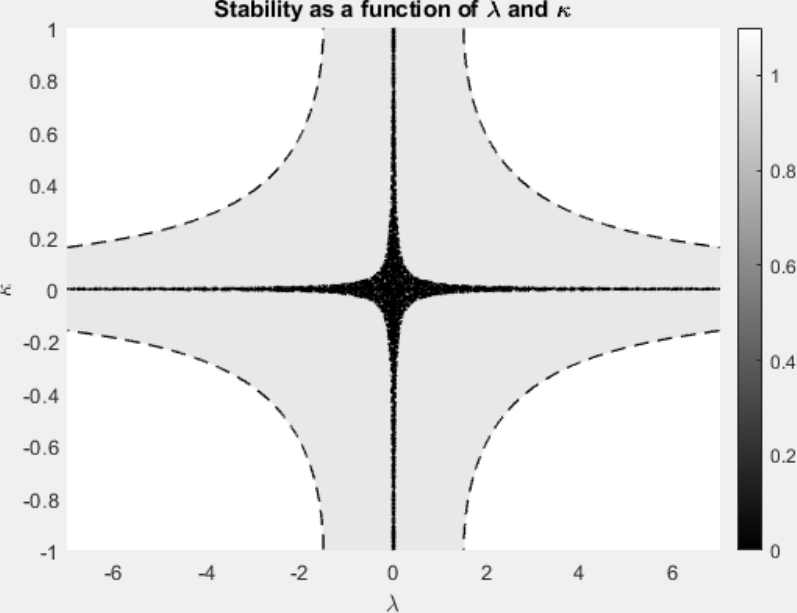}
\endminipage\hfill
\minipage{0.19\textwidth}
  \includegraphics[width=\linewidth]{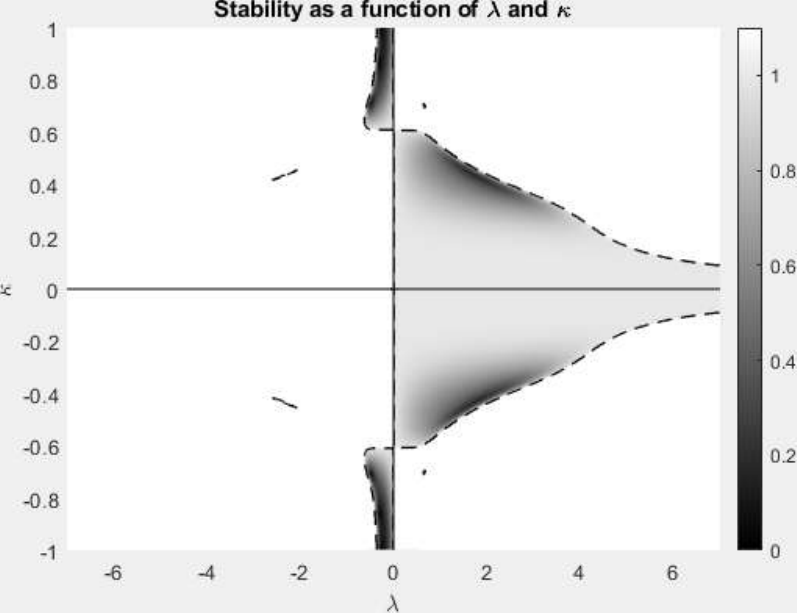}
\endminipage\hfill
\minipage{0.19\textwidth}
  \includegraphics[width=\linewidth]{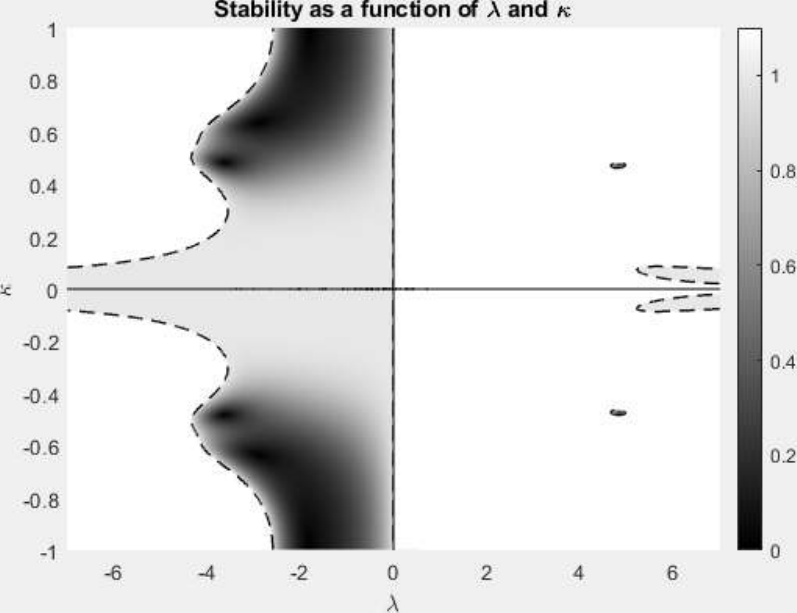}
\endminipage\hfill
\minipage{0.19\textwidth}%
  \includegraphics[width=\linewidth]{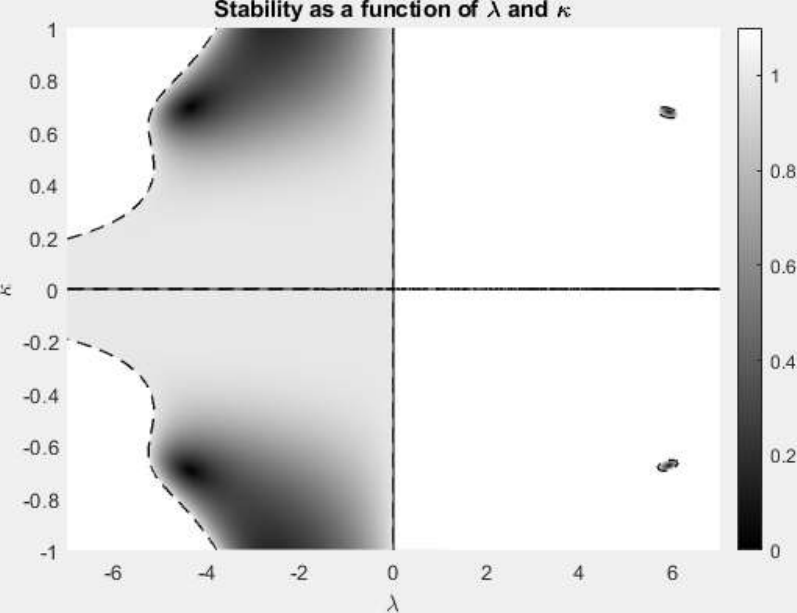}
\endminipage
\caption{RK7; left to right : stability regions for \textbf{Order 4} c, cs, f, zff and zffs.\label{rk74}}
\end{figure}

\begin{figure}[H]
\minipage{0.19\textwidth}
  \includegraphics[width=\linewidth]{blank.pdf}
\endminipage\hfill
\minipage{0.19\textwidth}
  \includegraphics[width=\linewidth]{blank.pdf}
\endminipage\hfill
\minipage{0.19\textwidth}
  \includegraphics[width=\linewidth]{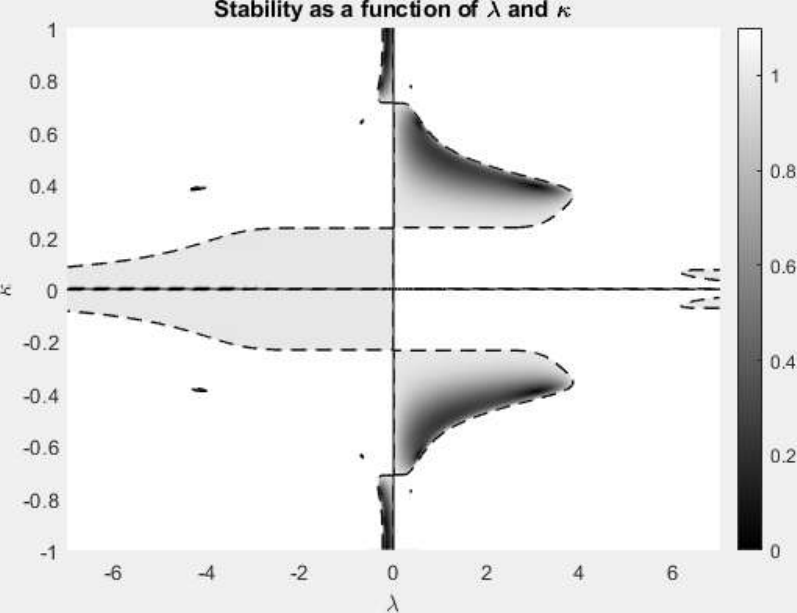}
\endminipage\hfill
\minipage{0.19\textwidth}
  \includegraphics[width=\linewidth]{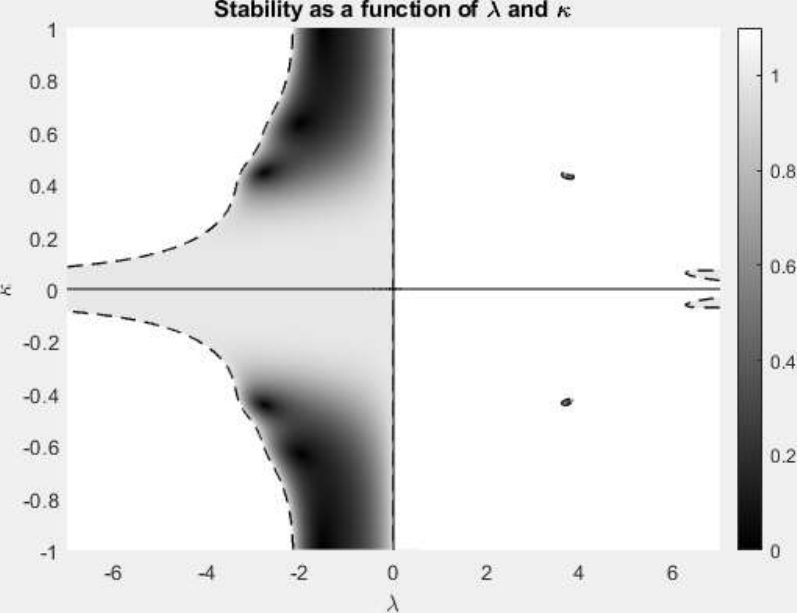}
\endminipage\hfill
\minipage{0.19\textwidth}%
  \includegraphics[width=\linewidth]{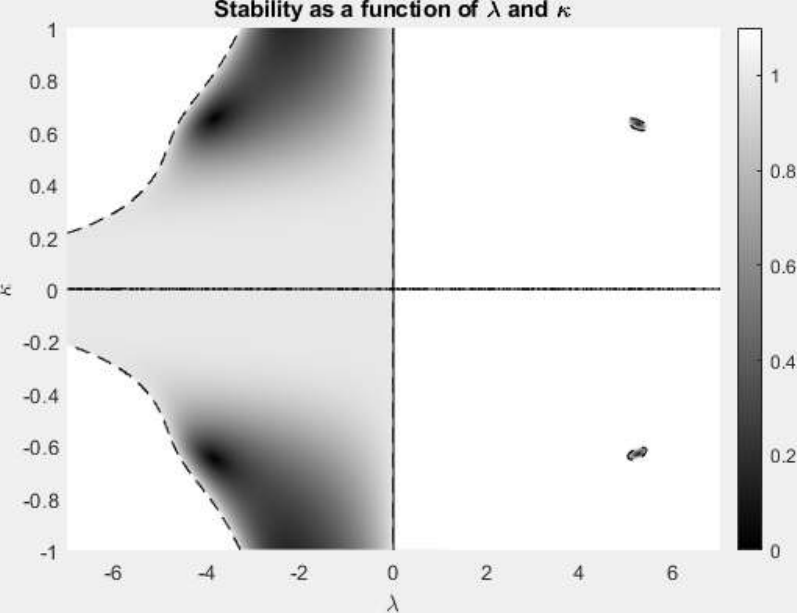}
\endminipage
\caption{RK7; left to right : stability regions for \textbf{Order 5} c, cs, f, zff and zffs.\label{rk75}}
\end{figure}

\begin{figure}[H]
\minipage{0.19\textwidth}
  \includegraphics[width=\linewidth]{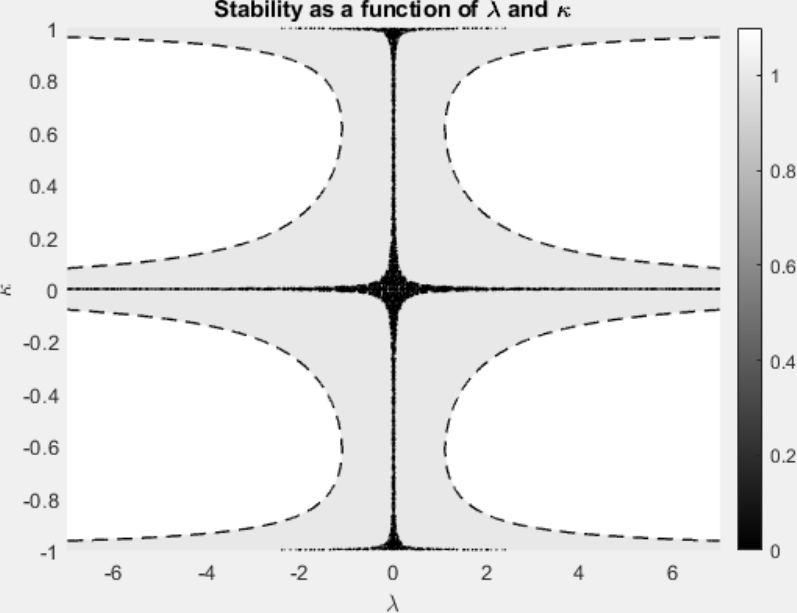}
\endminipage\hfill
\minipage{0.19\textwidth}
  \includegraphics[width=\linewidth]{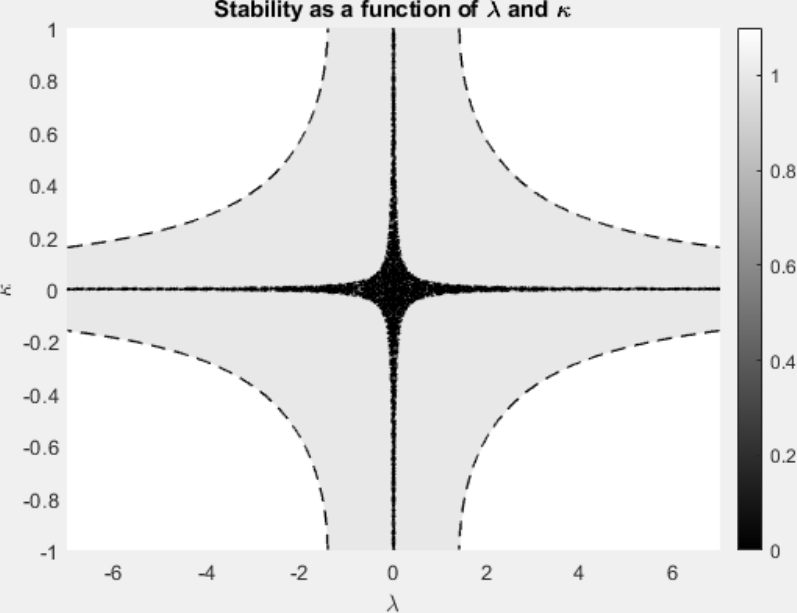}
\endminipage\hfill
\minipage{0.19\textwidth}
  \includegraphics[width=\linewidth]{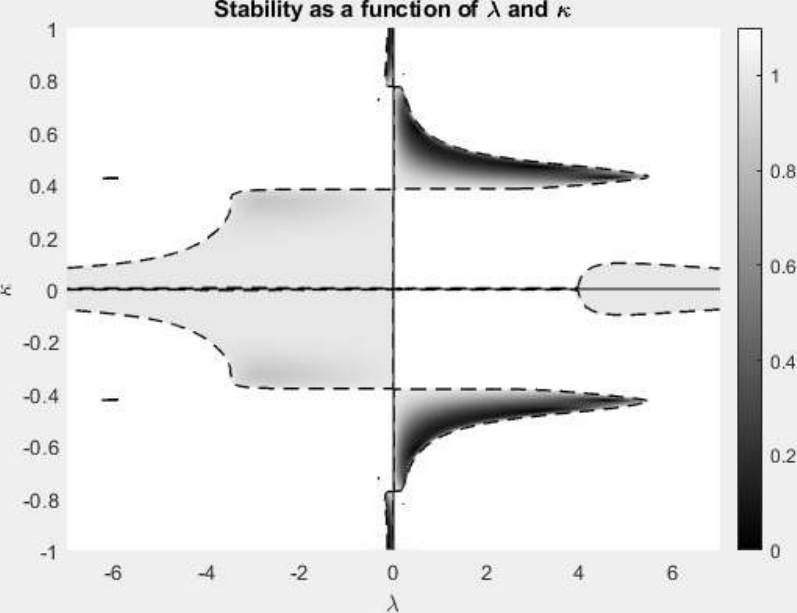}
\endminipage\hfill
\minipage{0.19\textwidth}
  \includegraphics[width=\linewidth]{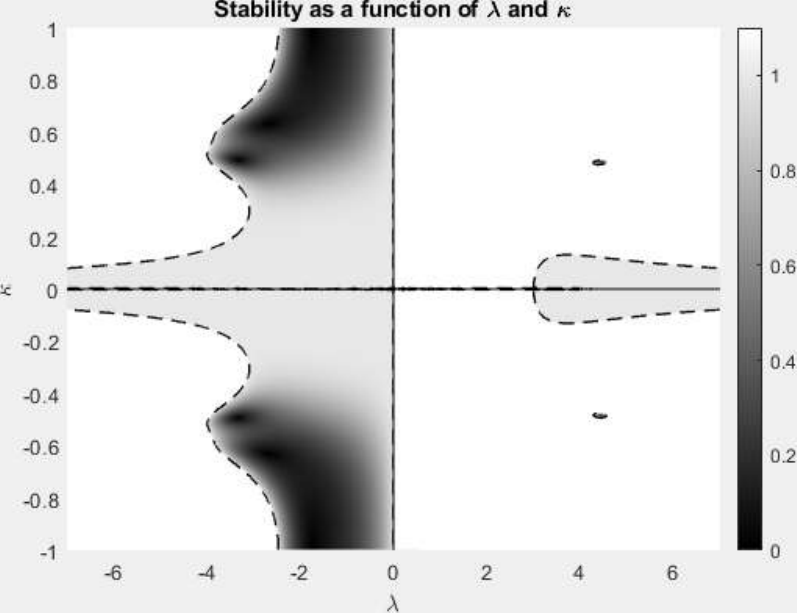}
\endminipage\hfill
\minipage{0.19\textwidth}%
  \includegraphics[width=\linewidth]{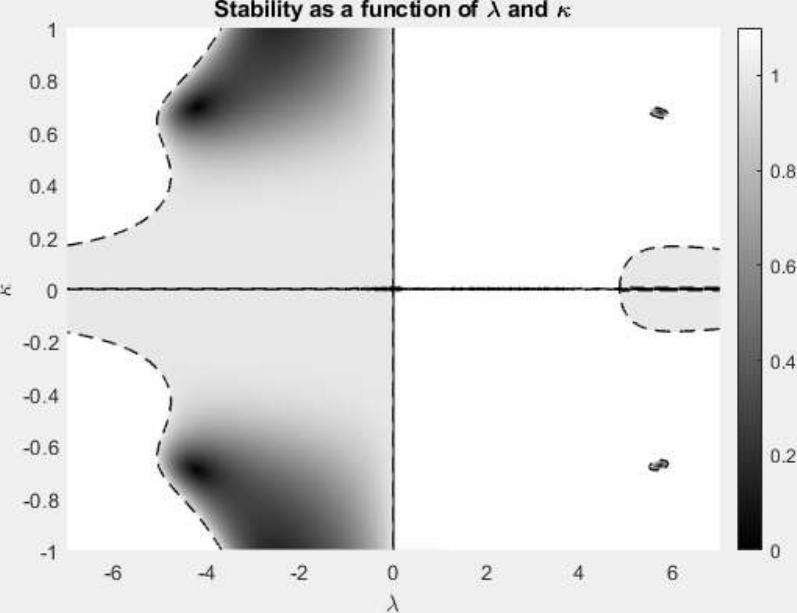}
\endminipage
\caption{RK7; left to right : stability regions for \textbf{Order 6} c, cs, f, zff and zffs.\label{rk76}}
\end{figure}

Fig. \ref{rk71} to \ref{rk76} show  the stability regions in the $\lambda$---$\kappa$ plane 
($\lambda$ and $\kappa$ respectively being the stability-number and the scaled wavenumber) for RK7 time integrator, 
for spatial orders $N=1$ to $6$, 
respectively (left to right) for centred (c), centred staggered (cs), forward (f), zigzag forward-first (zff) and zigzag forward-first staggered (zffs) schemes.
The greyed out areas 
represent the couples $(\lambda,\kappa)$ for which the scheme is stable 
(i.e. the amplification factor $|G| \leq 1$), while the dotted contours 
correspond to the critical case $|G| = 1$. A given scheme is conditionally
stable if and only if there exist a $\lambda_c \in \mathds{R}$ such that
the line of equation $\lambda = \lambda_c$ is included in the grey area.

\end{document}